\documentclass[a4paper, 10pt, reqno]{amsart}

\usepackage{amssymb,amsthm,amsmath,mathrsfs}
\usepackage{mathtools}
\usepackage{bbm}
\usepackage{stackengine}
\usepackage{enumitem}
\usepackage{soul}
\usepackage{thm-restate}
\usepackage{subcaption}

\usepackage[foot]{amsaddr}

\usepackage{caption,color,graphicx}
\usepackage[dvipsnames]{xcolor}

\usepackage{hyperref}
\hypersetup{colorlinks  = true,
            urlcolor    = blue,
            linkcolor   = red,
            citecolor   = blue}

\usepackage{tikz-cd}

\usepackage{verbatim}

\usepackage{environ}
\mathtoolsset{showonlyrefs,showmanualtags}

\usepackage[centering, top=3cm, bottom=3cm]{geometry}

\allowdisplaybreaks

\newcommand{\N}{\mathbb{N}}
\newcommand{\Z}{\mathbb{Z}}

\newcommand{\R}{\mathbb{R}}
\newcommand{\C}{\mathbb{C}}
\newcommand{\E}{\mathbb{E}}
\renewcommand{\P}{\mathbb{P}}
\newcommand{\1}{\mathbbm{1}}

\newcommand{\supp}{\mathrm{supp}\,}

\renewcommand{\d}{\mathrm{d}}
\newcommand{\dx}{\mathrm{d}x}

\newcommand{\Leb}{\mathrm{Leb}}

\newcommand{\calP}{\mathcal{P}}

\newcommand{\U}{\mathcal{U}}
\newcommand{\e}{\varepsilon}
\newcommand{\w}{\omega}

\newcommand{\dist}{\mathrm{dist}}
\newcommand{\Cp}{\mathbf{C}_+}
\newcommand{\Cm}{\mathbf{C}_-}
\newcommand{\Phip}{\Phi_+}
\newcommand{\Phim}{\Phi_-}

\newcommand{\vertiii}[1]{{\left\vert\kern-0.25ex\left\vert\kern-0.25ex\left\vert #1 \right\vert\kern-0.25ex\right\vert\kern-0.25ex\right\vert}}

\let\oldtocsection=\tocsection
\let\oldtocsubsection=\tocsubsection
\renewcommand{\tocsection}[2]{\hspace{0em}\oldtocsection{#1}{#2}}
\renewcommand{\tocsubsection}[2]{\hspace{1em}\oldtocsubsection{#1}{#2}}

\theoremstyle{plain}
\newtheorem*{theorem*}{Theorem}
\newtheorem{theorem}{Theorem}[section]
\newtheorem{lemma}[theorem]{Lemma}
\newtheorem{proposition}[theorem]{Proposition}
\newtheorem{corollary}[theorem]{Corollary}

\newtheorem{manualtheorem}{Theorem}

\newtheorem{manualhyp}{Hypothesis}
\newenvironment{hypothesis}[1][]{%
    \renewcommand\themanualhyp{#1}%
    \begin{manualhyp}
}{\end{manualhyp}}

\theoremstyle{remark}
\newtheorem{definition}[theorem]{Definition}
\newtheorem{remark}[theorem]{Remark}
\newtheorem{example}[theorem]{Example}
\newtheorem{notation}[theorem]{Notation}
  
\newtheorem{manualstep}{\bf Step}
\newenvironment{step}[1][]{%
    \renewcommand\themanualstep{#1}%
    \begin{manualstep}\itshape
}{\end{manualstep}}

\newenvironment{stepproof}[1][]{%
  \begin{proof}[#1]%
}{%
  \end{proof}%
}

\IfPackageLoadedTF{MnSymbol}{\newcommand{\coloneqq}{:=}}{}

\newlist{todolist}{itemize}{2}
\setlist[todolist]{label=$\square$}
\usepackage{pifont}

\usepackage{doi}

\usepackage[numbers]{natbib}

\author[Bernat Bassols Cornudella]{Bernat Bassols Cornudella$^{1}$}
\email{\href{mailto:bernat.bassols-cornudella20@imperial.ac.uk}{bernat.bassols-cornudella20@imperial.ac.uk}}

\author[Matheus M.~Castro]{Matheus M. Castro$^{2}$}
\email{\href{mailto:mmcastro@ime.unicamp.br}{mmcastro@ime.unicamp.br}}

\address{$^{1}$Department of Mathematics, Imperial College London, London SW7 2AZ, UK}

\address{$^{2}$Departamento de Matemática – Universidade Estadual de Campinas (UNICAMP), 13083-859, Campinas, SP, Brazil.}

\title[Conditioned stochastic stability for hyperbolic sets]{Conditioned stochastic stability of equilibrium states on uniformly hyperbolic sets} 
\date{\today}

\begin{document}

\begin{abstract}
    We establish the conditioned stochastic stability of equilibrium states for Hölder potentials on uniformly hyperbolic sets. While standard stochastic stability characterises measures on attractors, we analyse the statistics of transient dynamics on non-attracting sets by conditioning small random perturbations of the dynamics to not escape from our regions of interest. We prove that as the noise intensity vanishes, the quasi-ergodic measure of the $e^\phi$-weighted process generated by $\e$-small random perturbations of the deterministic dynamics converges to the unique equilibrium state associated with the potential $\phi - \log \left|\det \left. D T\right|_{E^u}\right|$. The results are obtained via perturbative spectral analysis of transfer operators acting on anisotropic Banach spaces and topological hyperbolic dynamics arguments. Furthermore, we extend this framework globally to Axiom A diffeomorphisms with multiple basic sets using dynamical filtrations. 
    This work provides a rigorous characterisation of natural measures on uniformly hyperbolic repellers, which are fundamental in the context of transient chaos. Our results have been formalised in Lean~4.
\end{abstract}

\maketitle
{\small\tableofcontents}

\newpage
\section{Introduction}\label{sec:intro}

The long-term statistical behaviour of trajectories for typical initial conditions in a dynamical system is well-known to be characterised by ergodic measures.
Given a map $T:M \to M$ and a $T$-invariant ergodic measure $\nu$, Birkhoff's ergodic theorem states that, for any bounded measurable observable ${\varphi:M \to \mathbb R}$, ``time averages equal space averages'' in the sense that
\begin{equation}\label{eq:BET}
    \lim_{n\to \infty}\frac{1}{n}\sum_{i = 0}^{n-1}\varphi \circ T^i (x) = \int \varphi\, \d\nu,\quad \nu\text{-a.s. on }x.
\end{equation}
Importantly, ``typical initial conditions'' refers to \eqref{eq:BET} being a ``$\nu$-almost sure'' statement and highlights the possible coexistence of many (even uncountably many) ergodic invariant measures.
In this context, the notion of \emph{stochastic stability} provides a strategy for identifying those ergodic invariant measures that persist under small random perturbations. Stochastic stability consists of: (i) proving the existence and uniqueness of a stationary measure $\eta_\e$ for the process $X_\e$ generated by $\e$-small random perturbations of $T$, and (ii) characterising the limit (in weak-$^*$) of $\eta_\e$ as $\e \to 0$. If $\eta_\e \to \eta$ as $\e \to 0$, then $\eta$ is an invariant measure selected by the prescribed family of random perturbations and, in many classical settings, is ergodic and uniquely selected. In this case, we say that $\eta$ is stochastically stable~\cite{Kifer1974, Alves2007}.
By construction, this strategy is fit to identify invariant measures sitting on attractors, near which $\eta_\e$ accumulates most mass.\footnote{This is not always the case; see for example \cite{Gora1985} or \cite{AraujoTahzibi2005} where the stochastically stable measure of LSV maps with $\alpha \geq 1$ is $\delta_0$, which sits on a neutral fixed
point.}

Repelling (or non-attracting) invariant sets may also support many ergodic invariant measures. Nevertheless, the stochastic stability of such measures cannot be characterised with the strategy presented above since the stationary measure of $X_\e$ will often attribute little-to-no mass near such sets and instead will highlight attractors~\cite{KavitaZeitouni1999}. To address this issue, we introduced the notion of ``conditioned stochastic stability'' in \cite{BassolsCastroLamb2024} where, instead of stationary measures, we proposed considering \emph{quasi-ergodic measures} of the process $X_\e$ conditioned upon remaining near the repeller.

In essence, conditioned stochastic stability consists of: (i) proving the existence and uniqueness of a quasi-ergodic measure $\nu_\e$ for the process $X_\e$ generated by $\e$-small random perturbations of $T$ and \emph{conditioned} upon remaining in a neighbourhood of a particular set $\Lambda$, and (ii) characterising the limit (in weak-$^*$) of $\nu_\e$ as $\e \to 0$. In \cite{BassolsCastroLamb2024}, we applied this strategy to uniformly expanding sets $\Lambda$, i.e.~repellers, and identified that quasi-ergodic measures converge to so-called \emph{equilibrium states} in the limit of $\e \to 0$, thus being \emph{conditionally stochastically stable}.

The study of conditioned stochastic stability provides a mathematically rigorous approach to transient chaos originating from the presence of a non-attracting chaotic set~\cite{TelGruiz2006, AltmannTel2009, LaiTel2011}. In particular, quasi-ergodic measures are the stochastic analogue of so-called ``natural measures'', which describe the statistical behaviour of the dynamics on such invariant non-attracting chaotic sets~\cite{Kantz1985}.
Natural measures are often computationally approximated via sampling techniques, such as the \emph{ensemble}~\cite{Kantz1985} and \emph{single-trajectory (PIM-triple)} methods~\cite{Nusse1989}. These methods aim to efficiently sample typical trajectories that remain near the non-attracting chaotic set for sufficiently long times, requiring repeated iterations of the underlying map. In this paper, we show that quasi-ergodic measures can be constructed by combining the dominant eigenfunctions of the random system's transfer operator and its dual (see Theorems~\ref{thm:qem-local} and \ref{thm:global-existence}), sidestepping cumbersome sampling strategies and offering a more robust procedure. In the process, we also obtain an approximation of the expected escape rate near a non-attracting chaotic set from the spectral radius of the transfer operator (see equation~\eqref{eq:growth-rate}), which aligns with well-established results~\cite{KadanoffTang1984, Pianigiani1979}.

Classical stochastic stability results for uniformly hyperbolic systems often concern equilibrium states associated with the geometric potential and show that stationary measures of small random perturbations converge to the corresponding SRB measures (see, for instance, \cite{Kifer1974, Kifer1986, Young1986, CowiesonYoung2005}).
Stability of other equilibrium states has been addressed in several dynamical settings and under different types of perturbation, such as small, possibly random, changes of the potential or the dynamics. We summarise some of these results below.

In~\cite{Kifer1990}, Kifer considered small diffusive perturbations of deterministic flows with a fixed potential and showed that, as the noise vanishes, the corresponding principal eigenvalue converges to the topological pressure of the deterministic system, establishing stochastic stability results for its equilibrium states.
In the context of discrete-time random dynamical systems, Bogenschütz~\cite{Bogenschutz1996}, Maume-Deschamps~\cite{MaumeDeschamps1997} and Baladi~\cite{Baladi1997} considered small random perturbations of expanding dynamics, obtained by successive random applications of maps close to the original one, and established stability of the corresponding equilibrium states as the perturbation vanishes. For non-uniformly expanding maps and Hölder potentials, Varandas and Viana~\cite{VarandasViana2010} proved stochastic stability of equilibrium states under random perturbations of the dynamics, while Castro and Varandas~\cite{CastroVarandas2013} later obtained stronger stability results allowing perturbations of both the dynamics and the potential. Statistical stability of equilibrium states under deterministic perturbations has also been studied in \cite{FreitasTodd2011, AlvesRamosSiqueira2019, BilbaoBioniLucena2024}. In a complementary and related direction, for Axiom A systems, stability of equilibrium states under perturbations of the dynamics and the potential has also been obtained through spectral stability of the corresponding weighted transfer operators by Gouëzel--Liverani~\cite{Gouezel-Liverani2008}, Baladi--Tsujii~\cite{Baladi-Tsujii2007} and Baladi \cite{Baladi2018-book}.

In contrast with the works above, the hyperbolic invariant set of interest in this paper may be non-attracting, so small random perturbations typically allow trajectories to escape from its neighbourhood. We therefore study quasi-ergodic measures of the corresponding process conditioned upon survival, instead of stationary measures, and show that they converge to the deterministic equilibrium state in the zero-noise limit.

The extension from uniformly expanding repellers in \cite{BassolsCastroLamb2024} to uniformly hyperbolic sets in this paper is conceptually natural but technically challenging.
The obstructions are of two main types: spectral and topological. 
From the spectral perspective, the Hilbert cone techniques used in \cite{BassolsCastroLamb2024} are no longer available in the current setting, and we instead rely on the anisotropic Banach spaces $B^{t,s}$ of Baladi and Tsujii \cite{Baladi-Tsujii2008, Baladi2018-book}. For these spaces, applying the Keller--Liverani perturbative argument is less immediate than in \cite{BassolsCastroLamb2024}: the annealed transfer operator naturally enjoys strong positivity and spectral properties on $\mathcal C^0(M)$, whereas the deterministic hyperbolic transfer operator is studied on $B^{t,s}$, and one must carefully relate the spectral information obtained on these two different functional settings.
From the topological perspective, proving that the support of the quasi-ergodic measure contains the underlying hyperbolic basic set is substantially more delicate.
Transitivity on the basic set does not extend to a neighbourhood of it. This issue comes from the fact that $W^u(\Lambda)$, i.e.~the unstable set of the hyperbolic basic set $\Lambda$, is not an open set but a fractal, having therefore a complicated geometry. Finally, we provide a more detailed analysis of the global conditioned dynamics and global quasi-ergodic measures for systems with multiple repellers (see the discussion in Section~\ref{sec:statements} and Theorem~\ref{thm:global-global} in particular). This introduces a further difficulty absent from the expanding case: one must understand how the different hyperbolic basic sets are dynamically connected and how these connections interact with their respective spectral contributions to determine the conditioned dynamics. We employ dynamical filtrations in the spirit of Conley's fundamental theorem of dynamical systems~\cite{Conley1978, Norton1995}.

\subsection{Conditioned Markov processes}\label{sec:conditionedMP}

Given a Markov chain on a metric space $(M, d)$, suppose that we are only interested in its behaviour as long as it remains in a compact subset $\overline{\mathcal U} \subset M$, where $\mathcal{U}\subset M$ is open. We thus identify $M\setminus \overline{\mathcal U}$ with a \emph{cemetery} state $\partial$ and
consider the absorbing Markov chain given by
\[X \coloneqq (\Omega, \{\mathcal F_n\}_{n \in\mathbb N_0}, \{X_n\}_{n \in\mathbb N_0}, \{\mathbf P^n\}_{n \in\mathbb N_0}, \{\mathbb P_x\}_{x\in \overline{\mathcal U} \sqcup \partial})\]
in the usual sense (see e.g.~\cite[Definition III.1.1]{Rogers1994}), with $\partial:= M \setminus \overline{\mathcal U}$. Since $\partial$ is a cemetery state, we impose $\mathbf P^n(\partial, \partial) = 1$ for all $n \in \mathbb N_0$, i.e.~$X_n$ is absorbed in $\partial$.

Given a non-positive continuous \emph{potential} function $\phi : M \to (-\infty, 0]$, we lift $X$ to the new process $X^{\phi}$ given by (see also \cite[Section~8.1.2]{Gouezel-Liverani2008} and \cite[Remark~5.23]{Baladi2018-book})
\begin{equation}\label{eq:weighted}
X_{n+1}^{\phi} = \begin{cases}
    X_{n+1},&\text{with prob. }e^{\phi(X_n)}\\
    \partial,&\text{with prob. }1 - e^{\phi(X_n)}.
\end{cases}
\end{equation}
The process $X^{\phi}$ is again a Markov chain, with transition kernel
$\mathbf P_\phi$ satisfying
$$
\mathbf P_\phi(x,B)
\coloneqq
\begin{cases}
e^{\phi(x)}\mathbf P(x,B),
& B\subset\overline{\mathcal U}\ \text{measurable},\\[1mm]
1-e^{\phi(x)}\mathbf P(x,\overline{\mathcal U}),
& B=\partial,
\end{cases}
\quad x\in\overline{\mathcal U}.
$$
and $\mathbf P_\phi(\partial,\{\partial\})=1$. We denote by $\mathbb P_x^\phi$ the path measure of $X^\phi$ starting
from $x$ and by $\E_x^\phi$ the corresponding expectation. Similarly,
$\E_x$ denotes expectation with respect to $\mathbb P_x$. If
$\phi\equiv0$, then $\mathbf P_\phi=\mathbf P$, and hence
$\mathbb P_x^\phi=\mathbb P_x$ and $X_n^\phi=X_n$. 

Here, we are interested in the behaviour of $X_n^\phi$ before it is absorbed by $\partial$ or, in other words, before the stopping time 
\[\tau(\omega, x) \coloneqq \inf\{n > 0 :\, X_n^\phi \in \partial\} = \inf\{n > 0:X_n^\phi \in M \setminus \overline{\mathcal U}\}, \quad (\w, x )\in \Omega \times M\]
occurs. As mentioned above, quasi-ergodic measures are defined to capture the statistical behaviour of the process $X_n^\phi$ conditioned upon not being absorbed, i.e.~conditioned upon $\tau >n$. To be more precise:

\begin{definition}[Quasi-ergodic measure on $\overline{\U}\subset M$]  
    We say that the probability measure $\nu$ on $M$ is a \emph{quasi-ergodic measure on $\overline{\U} \subset M$} for the $e^\phi$-weighted process $X^\phi$ if for any bounded and measurable observable $\varphi:\overline{\U} \to \mathbb R$ it holds that 
\begin{equation}\label{eq:qem-cond}
\begin{split}
\E_x^\phi\left[\frac{1}{n}\sum_{i= 0}^{n-1}\varphi \circ X_i^\phi \, \bigg|\, \tau >n \right]
\xrightarrow[]{n \to \infty} \int \varphi \d\nu, \quad \nu\text{-a.s. on }x.
\end{split}
\end{equation}
The conditional expectation on the left-hand side is defined as 
\begin{equation}\label{eq:qem}
\frac{1}{\P_x^\phi\left[\tau > n\right]} \E_x^\phi\left[\1_{\{\tau > n\}}\frac{1}{n}\sum_{i= 0}^{n-1} \varphi \circ X_i^\phi\right] = \frac{1}{\E_x[e^{S_n\phi} \1_{\overline{\U}}(X_n)]}\E_x\left[e^{S_n\phi} \1_{\overline{\U}}(X_n) \frac{1}{n} \sum_{i = 0}^{n-1}\varphi \circ X_i \right],
\end{equation}
where $S_n \phi = \sum_{i=0}^{n-1}\phi \circ X_i$ denotes the Birkhoff sum.
\end{definition}

\begin{remark}\label{rmk:positive}
Observe that this definition is well-posed for any continuous $\phi$ and not only non-positive weights. For a general $\phi$, we interpret $X^\phi$ as follows: the process evolves in time and carries a certain amount of mass. If the process enters (dynamically) the region $\partial$, then it is absorbed and all mass is lost. At every time step $n$ and in every position $X_n$, the process will increase (or decrease) its mass by a factor $e^{\phi(X_n)}$. Conditioning upon $\tau > n$ as in \eqref{eq:qem-cond} must then be interpreted as the right-hand side in \eqref{eq:qem}.

More generally, if $\varphi:M\to[0,\infty)$ is a continuous weight which is allowed to vanish, we use the notation $X^{\log \varphi}$ through the non-negative weight $\varphi$, with the convention $\log0=-\infty$. In particular, $e^{S_n\log \varphi}\coloneqq\prod_{j=0}^{n-1}\varphi(X_j),$ which is well defined even when $\varphi$ vanishes. When $0\leq \varphi\leq1$, this has the usual absorbing interpretation: if $\varphi(x)=0$, the process is absorbed with probability one at the next transition.
\end{remark}

The study of Markov processes conditioned upon never entering a cemetery state has long been present in the literature \cite{Yaglom1947, Pinsky1985, Breyer1999, Collet2013} and has more recently become of increasing interest \cite{Champagnat2016, Champagnat2018, Colonius2021, Castro2024,Castro2022}. To prove the existence and the uniqueness of the quasi-ergodic measure for a process $X_\e$, it is often required for the chain to be strong Feller and transitive (see e.g.~\cite[Appendix]{BassolsCastroLamb2024}). Nevertheless, this assumption is no longer valid when considering small random perturbations of hyperbolic (non-expanding) deterministic maps, as done herein. To address this shortcoming, we provide new results guaranteeing the existence and uniqueness of quasi-ergodic measures in the small-noise regime, whose support contains relevant dynamical information (see Theorems~\ref{thm:qem-local} and~\ref{thm:global-existence} guaranteeing $\Lambda \subset \supp \nu$).

\subsection{Conditioned stochastic stability}

The notion of conditioned stochastic stability presented in \cite{BassolsCastroLamb2024} provides a natural candidate to highlight the statistics of the transient behaviour of a map $T$.
As mentioned above, we may consider the quasi-ergodic measure of the Markov process $X_\e$ generated by small, $\e$-bounded, random perturbations of a map $T:M \to M$ and conditioned upon remaining in a particular region of interest, e.g.~where the transient evolves. As $\e \to 0$, if the (weak-$^*$) limit\footnote{While we only consider weak-$^*$ convergence, one may also ask for the convergence in \emph{total variation} of the quasi-ergodic measure, which would refer to \emph{strong} conditioned stochastic stability \cite{Shen-vanStrien2013}.} of this quasi-ergodic measure exists and is unique, we have effectively identified a $T$-invariant measure \cite{Khasminskii1963} that is robust under small random perturbations.

We distinguish the complement of two regions in state space with the potential to highlight dynamical transients: (i) $\partial = M \setminus \overline{V}$, where $V$ is the (isolating) neighbourhood of an invariant set (see Definition~\ref{def:hyp-basic-set}), and (ii) $\partial = \mathcal A,$ an open neighbourhood of a forward invariant set for $T$ which we may think of as an attractor or trapping region. We refer to these settings as the local and global problems, respectively. Conditioned stochastic stability, therefore, concerns understanding the limiting behaviour of the quasi-ergodic measure $\nu_\e$ for the process $X_\e$ conditioned upon not entering $\partial$ as $\e \to 0$.

In this paper, we show that the quasi-ergodic measure of the process $X_\e$ approximates, as $\e \to 0$, the \emph{equilibrium state}
associated with the geometric potential $\psi = - \log |\det DT\vert_{E^u}|$ (see Section~\ref{sec:TDF-natural_measures}), where $E^u$ denotes the unstable expanding direction of $T$. This holds both locally and globally. This extends the previous results in \cite{BassolsCastroLamb2024}, where repellers are assumed to be uniformly expanding. More generally, we show that the quasi-ergodic measure for the process $X_\e^\phi$ approximates the equilibrium state associated with the potential $\psi = \phi - \log|\det DT\vert_{E^u}|$ when restricted to $\Lambda$.\\

The rest of the manuscript is organised as follows. In Section~\ref{sec:TDF-natural_measures}, we briefly introduce equilibrium states, as they feature in our results when $\e \to 0$, and relate quasi-ergodic to so-called \emph{natural} measures, which lie at the heart of the transient chaos literature. In Section~\ref{sec:setup} we formally introduce the systems considered (Section~\ref{sec:deterministic}) and their random perturbations (Section~\ref{sec:random_perturb}) to state the main results of the paper in Section~\ref{sec:statements} (Theorems~\ref{thm:qem-local}, \ref{thm:local} for the local setting and Theorems~\ref{thm:global-existence}, \ref{thm:global} for the global setting). We also provide some examples in Section~\ref{sec:examples} to illustrate the novelty of these results. Section~\ref{sec:BanachSpaces} introduces the (anisotropic) Banach spaces $B^{t,s}$ of \cite[Chapter~4]{Baladi2018-book}, which we employ to establish spectral properties of the transfer operators considered. We dedicate Sections~\ref{sec:QEMs-vanishing} and~\ref{sec:QEMs-non-vanishing} to the local problem, where the random process is conditioned upon staying in a neighbourhood $V$ of a hyperbolic basic set $\Lambda$. For weights vanishing on $\partial V$, Section~\ref{sec:QEMs-vanishing} provides the existence of the quasi-ergodic measure as well as conditioned stochastic stability in the local context. These results are extended to non-vanishing weights in Section~\ref{sec:QEMs-non-vanishing}. Finally, in Section~\ref{sec:proofs}, we prove the uniqueness of the quasi-ergodic measure supporting the hyperbolic set $\Lambda$ (Section~\ref{sec:proof_local}) and promote all these results to the global setting (Section~\ref{sec:global-proof}) by means of dynamical filtrations \cite{Conley1978}. All results have been formalised in Lean~4.35.0 as explained in Section~\ref{sec:lean-formalisation}. We dedicate the Appendix~\ref{sec:appendix} to the statement of two theorems exploited in Section~\ref{sec:QEMs-vanishing}.

\subsection{Equilibrium states and natural measures}\label{sec:TDF-natural_measures}

The theory of \emph{thermodynamic formalism} concerns the study of so-called \emph{equilibrium states}. For a map $T: \Lambda \to \Lambda$ and a continuous weight $g:\Lambda \to \R$, an equilibrium state is a $T$-invariant (ergodic) Borel probability measure $\nu$ on $\Lambda$ realising the supremum of the \emph{metric pressure}
\begin{equation*}
        P_\mu(T, \Lambda, \psi) = h_\mu(T,\Lambda) + \int \psi \, \d\mu
\end{equation*}
over all such $T$-invariant ergodic measures, where $h_\mu(T, \Lambda)$ denotes the Kolmogorov-Sinai (metric) entropy \cite{Kolmogorov1958, Sinai1959} and $\psi = \log g - \log |\det DT|_{E^u}|$. These measures were introduced in the foundational work of Sinai, Ruelle and Bowen \cite{Sinai1972, Bowen1975-book, BowenRuelle1975, Ruelle1976, Ruelle1978} and are motivated by ideas from the theory of statistical mechanics as they are set to minimise a theoretical analogue of the system's free energy.
Amongst equilibrium states, we distinguish those associated with the ``geometric'' potential $\psi = - \log |\det DT\vert_{E^u}|$, where $E^u$ denotes the unstable expanding direction of $T$, i.e.~$g = 1$. It is well-known that these equilibrium states correspond to the Sinai-Ruelle-Bowen (SRB) measure for hyperbolic attractors, which characterises the (natural) distribution of Lebesgue typical orbits on $\Lambda$ \cite{Young2002}. 

From a different point of view, empirical results in transient chaos show that the transient dynamics of a map $T$ near a hyperbolic invariant set $\Lambda$ are governed by the statistics of nearby points in $\Lambda$ \cite{GrebogiOttYorke1983, Lai2011}. In this context, the so-called \emph{natural measure} on $\Lambda$ may be computed from the histogram of trajectories that remain close to $\Lambda$ for sufficiently long times \cite{Kantz1985, Nusse1989, Bassols2023}. Quasi-ergodic measures, as defined above, and their limit as $\e \to 0$, aim to formalise these objects. In the process, we provide a new technique to approximate natural measures as the limit of quasi-ergodic measures, which, in turn, may be approximated from the dominant eigenfunctions of the (annealed) transfer operator studied herein. 

\section{Setup, results and examples} \label{sec:setup}
\subsection{The deterministic dynamics}\label{sec:deterministic}
Throughout this paper, unless stated otherwise, the map $T:M \to M$ is a $\mathcal C^r(M, M)$ diffeomorphism, $r \geq 1$, on an orientable compact Riemannian manifold $M$ and $g \in \mathcal C^{r-1}_0(M, \mathbb R^+_0)$ denotes a weight function, $\mathbb R^+_0 \coloneqq [0,\infty)$. 
Denote by {$\Leb$} a Borel measure on $M$ induced by a smooth volume form compatible with this metric, which we may think of as Lebesgue.

We consider diffeomorphisms $T$ for which a locally maximal hyperbolic set $\Lambda$ exists.

\begin{definition}[Locally maximal hyperbolic set, isolating neighbourhood] A $T$-invariant compact set $\Lambda$ is called \emph{hyperbolic} if there exists a continuous invariant decomposition $T_\Lambda M = E^u \oplus E^s$ of the tangent bundle over $\Lambda$ into two $DT$-invariant sub-bundles, and there exist constants $C>0$ and $0<\rho<1$ such that
$$\|\left.D T_x^m\right|_{E_x^s}\|\leq C \rho^m,\qquad \|\left.D T_x^{-m}\right|_{E_x^u}\|\leq C \rho^m,$$ 
for every $m\geq 0$ and $x\in\Lambda$.
The hyperbolic set $\Lambda$ is called \emph{locally maximal} (or isolated) if it admits an open neighbourhood $V$ such that $\Lambda = \cap_{m\in\mathbb Z} T^m(\overline{V}).$ The set $V$ is called an \emph{isolating neighbourhood} of $\Lambda$.
\end{definition}

\begin{definition}[{Hyperbolic basic set~\cite[Definition~10.1.3, adapted]{HasselblattKatok2003}}]\label{def:hyp-basic-set} A hyperbolic set $\Lambda$ is called transitive if $T$ has a dense orbit in $\Lambda$. A \emph{hyperbolic basic set} for $T$ is a transitive locally maximal hyperbolic set for $T$.  
\end{definition}

In particular, transitivity and continuity of the fibre bundles imply that the dimensions of $E^u(x)$ and $E^s(x)$ are constant on a hyperbolic basic set $\Lambda$~\cite{AraujoViana2012}. We denote these by $d_u$ and $d_s$, respectively. If $d_s=0$, then $\Lambda$ is a periodic repelling orbit, and this uniformly expanding case is already covered in~\cite{BassolsCastroLamb2024}. Since our interest here is in non-attracting hyperbolic sets, we exclude the case $d_u=0$, which corresponds to a periodic attracting orbit. Hence, in the construction of the anisotropic spaces below, we may assume that $d_s,d_u\geq1$.

\begin{remark}\label{rmk:non-transitive}
If $\Lambda$ is a locally maximal hyperbolic set but non-transitive, then it is well-known that there exists a partition of $R=  \overline{\mathrm{Per}(T\vert_\Lambda)} \coloneqq \overline{\{p \in \Lambda:\, p \text{ is $T$-periodic}\}}$ into finitely many non-empty compact subsets $R^{i,j}, 1 \leq i \leq k$, $1\leq j \leq m(i)$ such that for every $i,j$ \cite[Theorem~10.3.6]{HasselblattKatok2003}:
\begin{enumerate}[label = (\roman*)]
    \item $R^{i} = \cup_{j=1}^{m(i)} R^{i,j}$ is $T$-invariant,
    \item $T(R^{i,j}) = R^{i,j+1\ (\mathrm{mod}\ m(i))}$,
    \item $T:R^i\to R^i$ is uniformly hyperbolic and topologically transitive, and
    \item each $T^{m(i)}:R^{i,j}\to R^{i,j}$ is uniformly hyperbolic and topologically mixing.
\end{enumerate}
Thus, we may always reduce to the hyperbolic basic setting when dealing with the local problem by taking a small enough isolating neighbourhood around each $R^i$ (see e.g.~\cite[Section~18.3.1]{Katok1995}).
\end{remark}

Equilibrium states for H\"older weights on hyperbolic basic sets are well known to exist and to be unique~\cite{Bowen1975-book, Baladi-Tsujii2007, Gouezel-Liverani2008}. This result is essential for understanding the limit of quasi-ergodic measures as $\e \to 0$, as we do not develop a thermodynamic formalism.

The so-called global problem in this paper concerns Axiom A diffeomorphisms. Let us recall some basic definitions for such maps.

\begin{definition}[{Set of non-wandering points, $NW(T)$}]\label{def:NW}
    A point $x \in M$ is said to be \emph{non-wandering} if for every neighbourhood $U \ni x$ there is $n \geq 1 $ such that $T^n(U) \cap U \neq \emptyset$. We denote by $NW(T)$ the \emph{set of non-wandering points}.
\end{definition}

\begin{definition}[Axiom A]\label{def:AxiomA}
    We say that a $\mathcal C^r(M)$ diffeomorphism $T:M \to M$, $r> 1$, satisfies \emph{Axiom A} if:
    \begin{enumerate}[label = (\roman*)]
        \item $NW(T)$ is hyperbolic, and
        \item the set of $T$-periodic points is dense in $NW(T)$, i.e.~$NW(T) = \overline{\mathrm{Per}(T)}$.
    \end{enumerate}
\end{definition}

If $T$ satisfies Axiom A, then observe that the density of periodic points in $NW(T)$ implies that $NW(T)$ is a locally maximal hyperbolic set, so we may take $\Lambda = NW(T)$. In this setting, we further set $\Lambda_i \coloneqq R^i$ of the dynamical or spectral decomposition in Remark~\ref{rmk:non-transitive} if $T$ is not transitive.

\subsection{The random perturbations}\label{sec:random_perturb}
For $\e >0$, consider the random perturbations of $T$ given by $F_\e:[-\e, \e]^m\times M \to M$, where $F_\e$ is jointly $\mathcal C^1$ in $(\omega,x)$, $F_\e(\omega, \cdot) \in \mathcal C^r(M)$ and $\partial_\w F_\e(\w, x)$ is surjective for all $\w \in [-\e, \e]^m$. We assume that $\dist_{\mathcal C^r}(F_\e(\w), T) \leq C \|\w\|$ for some uniform $C>0$, where $\dist_{\mathcal C^r}$ denotes the $\mathcal C^r$-Whitney topology \cite[Chapter~1.2]{Palis1982}. Note that $m \geq \dim M$ by surjectivity of $\partial_\w F_\e(\w, x)$. Denote by $\Omega_\e \coloneqq ([-\e, \e]^m)^{\mathbb Z}$ the space of bi-infinite sequences of elements in $[-\e, \e]^m$ endowed with the probability measure $\P_\e \coloneqq \upsilon^{\otimes \mathbb Z}$, where $\upsilon$ is an absolutely continuous probability measure on $[-\e, \e]^m$ of full support whose density with respect to Lebesgue measure is bounded. For $\w\in\Omega_\e$, we define
$T_\w(x)\coloneqq T_{\w_0}(x)\coloneqq F_\e(\w_0,x)$ and, for $n\in\mathbb N$, write
$$
T_\w^n(x)\coloneqq T_{\w_{n-1}}\circ\cdots\circ T_{\w_0}(x),
$$
with $T_\w^0=\mathrm{id}$. Moreover, we denote by $\theta:\Omega_\e\to\Omega_\e$ the two-sided shift, $(\theta\w)_i=\w_{i+1}$ for $i\in\mathbb Z$, and, for $n\geq1$, define
$$
T_\w^{-n}\coloneqq\left(T_{\theta^{-n}\w}^{n}\right)^{-1}
=T_{\w_{-n}}^{-1}\circ\cdots\circ T_{\w_{-1}}^{-1}.
$$
With this convention, $T_\w^{n+m}=T_{\theta^n\w}^{m}\circ T_\w^n$ for every $m,n\in\mathbb Z$. We use $\Theta:\Omega_\e\times M\to\Omega_\e\times M$ to denote the skew-product map\footnote{The symbol $\Theta$ is also used in Section~\ref{sec:BanachSpaces} to denote cone systems but this will be clear by context.} $
\Theta(\w,x)=(\theta\w,T_\w(x)).$

Let $X_\e$ be the Markov process generated by $F_\e$, where $\Omega = \Omega_\e$, $X^\e_n(\w, x) \coloneqq T^n_\w(x)$ for $n \in \mathbb N$ and $X^\e_0(\w, x) = x$ for every $x \in M$, $\mathbb P_\e$-almost surely. Given a suitable $\phi:M \to \mathbb R$, we denote by $X^\phi_\e$ the $e^\phi$-weighted Markov process as described in Section~\ref{sec:conditionedMP}, where $\partial$ is an open subset of $M$. 
We remark that for such random perturbations, the annealed transfer operator associated with $X^\phi_\e$ it has the strong Feller property (see Lemma~\ref{lem:strong_feller}).

\begin{remark}\label{rmk:abuse}
Observe that $X_n^\phi = X_n^\phi(\omega,x) = T^n_\w(x)$ for all $n < \tau := \inf\{n > 0:X_n^\phi \in \partial\}$ and $X_n^\phi \in \partial$ for all $n \geq \tau$. In this paper, we are interested in the dynamics up until the absorbing time $\tau$. Thus, without loss of generality, we may redefine $T_\omega^n$ so that $T_\omega^n (x) = \partial$ if $n\geq\tau(\omega,x)$.
In particular, this allows us to write $\1_{A}(T^N_\w(x))$ instead of $\1_{A}(T^N_\w(x))\prod_{i=0}^{N-1}\1_{M\setminus \partial}(T^i_\w(x))$ for sets $A \subset M \setminus \partial$.
\end{remark}

\subsection{Main results}\label{sec:statements}
We are ready to state the main assumptions considered in this paper and the four theorems concerning the existence and the uniqueness of quasi-ergodic measures (Theorem~\ref{thm:qem-local} local, Theorem~\ref{thm:global-existence} global) and conditioned stochastic stability (Theorem~\ref{thm:local} local, Theorem~\ref{thm:global} global). 

In the local setting, we shall work under the following hypothesis:

\begin{hypothesis}[HL]\label{hyp:local}
    We say that a triple $(T,g,\Lambda)$ satisfies Hypothesis~\ref{hyp:local} if: 
\begin{enumerate}[label = (\roman*)]
    \item $T:M\to M$ is a $\mathcal C^r$ diffeomorphism with $r>1$,
    \item $\Lambda$ is a basic hyperbolic set with isolating neighbourhood $V$,
    \item $g:M \to \mathbb R$ is a $\mathcal C^{r-1}$ weight function supported in $V$, with {$g|_{\Lambda} > 0$}, and
    \item $T:\Lambda\to\Lambda$ is topologically mixing.
\end{enumerate}
\end{hypothesis}

Assuming Hypothesis~\ref{hyp:local}, there exists a unique equilibrium state $\nu^g$ on $\Lambda$ associated with $g$ which realises the topological pressure of the system (see \cite[Theorem~7.5]{Baladi2018-book} or \cite{Gouezel-Liverani2008})
    $$
    P_{\mathrm{top}}(T,\Lambda, \psi) = h_{\nu^g}(T) + \int \psi \, \d\nu^g,
    $$
with $\psi = \log g - \log|\det DT\vert_{E^u}|$.
A central object in the study of equilibrium states is the Ruelle transfer operator $\mathcal P_g$ on a suitable function space defined as 
\[\mathcal P_g \varphi(x) = g(x)(\varphi \circ T(x)).\]
Analogously, for the randomly perturbed and $e^\phi$-weighted process $X_\e^\phi$, we consider the (local) annealed transfer operator on $V$ instead given by
\[\calP_\e \varphi(x) = e^{\phi(x)}\E_\e\left[\varphi \circ T_\w(x) \cdot \1_{\overline{V}}\circ T_\w(x)\right].\]
The following result provides the existence and uniqueness of quasi-ergodic measures around each hyperbolic basic set.

\begin{theorem}\label{thm:qem-local}
    Assume that $(T,g,\Lambda)$ satisfies Hypothesis~\ref{hyp:local} and let $V$ be an isolating neighbourhood of $\Lambda$. Let $e^\phi\in\mathcal C^{r-1}(\overline V)$ and assume that $e^\phi=g$ on $\Lambda$.
 Then for $\e>0$ small enough, $X_\e^{\phi}$ admits a unique quasi-ergodic measure $\nu_\e^\phi$ on $\overline{V}$ such that $\Lambda \subset \supp \nu_\e^\phi$. Moreover, $\nu_\e^\phi(\varphi) = \mu_\e(\varphi \cdot g_\e)$ for any observable $\varphi\in \mathcal C^0(\overline{V})$, where $g_\e\in\mathcal C^0(\overline V)$ and $\mu_\e\in\mathcal M(\overline V)$ satisfy
$
\calP_\e g_\e=\lambda_\e g_\e$, $\calP_\e^*\mu_\e=\lambda_\e\mu_\e$, $\mu_\e(g_\e)=1,$ with $\lambda_\e=r(\calP_\e|_{\mathcal C^0(\overline V)})$.
\end{theorem}

Local conditioned stochastic stability is then given by: 

\begin{theorem}[Local conditioned stochastic stability, $\partial = M \setminus \overline{V}$]\label{thm:local}
    For each $\e>0$ small enough, let $\nu_\e$ be the quasi-ergodic measure on $\overline{V}$ from Theorem \ref{thm:qem-local}. Then $\nu_\e \to \nu^g$ in {weak-$^*$} as $\e \to 0$, where $\nu^g$ is the unique equilibrium state of $T$ on $\Lambda$ associated with the potential $\phi - \log |\left.\det DT\right|_{E^u}|.$
\end{theorem}

We prove Theorem~\ref{thm:qem-local} and Theorem~\ref{thm:local} in Section~\ref{sec:proof_local}.

In data-driven applications and the practical analysis of transient chaos, it may not be possible to identify the locus of a hyperbolic basic set $\Lambda$ nor a suitable isolating neighbourhood $V$. Following Theorem~\ref{thm:qem-local}, this information is crucial to approximate the transfer operator $\calP_\e$ and compute its dominant left and right eigenfunctions, which, combined, make up the quasi-ergodic measure and approximate the underlying equilibrium state.

On the other hand, provided enough time has elapsed, a system will often settle in a trapping region $\mathcal A$ on which the dynamics are no longer transient. While we may not be interested in the behaviour within $\mathcal A$, data provides us with access to its position in the state space. Therefore, it is well-justified to ask what the statistical behaviour of the process is before reaching this global trapping region.

\begin{remark}
    We may also choose the trapping region $\mathcal A$ to be the empty set, in which case the absorbing Markov process generated by small random perturbations is only subject to absorption by the action of the potential.
\end{remark}

In Theorems~\ref{thm:global-existence} and \ref{thm:global}, we address this question and show that the longest transient behaviour of the system is driven by the dynamics near the most dominant hyperbolic basic set, i.e.~that of the largest topological pressure and thus of smallest escape rate~\cite{DemersYoung2006}. In this context, we extend the previous Hypothesis~\ref{hyp:local} to distinguish the different basic sets $\Lambda_i, i \in \{1, \dots, k\}$ and ensure there are no cycles between them:

\begin{hypothesis}[HG1]\label{hyp:global}
 We say that $(T,g,\mathcal A)$ satisfies Hypothesis~\ref{hyp:global} on $M$ if:
    \begin{enumerate}[label = (\roman*)]
        \item $T:M \to M$ is an Axiom A diffeomorphism (see Definition~\ref{def:AxiomA}), and $T(\overline{\mathcal A}) \subset \mathcal A$ with $\mathcal A$ being an open set, \label{item:gh1i}
        \item there exists $r>1$ such that $T$ is a $\mathcal C^r(M)$ function and
$g\in\mathcal C^{r-1}(M,\mathbb R_0^+)$ is a non-negative weight function
such that $g|_{\Lambda_i}>0$ for every basic set
$\Lambda_i\subset M\setminus\mathcal A$,\footnote{The weight $g$ is global and need not be supported on a particular isolating neighbourhood $V$. Positivity is only required on the basic sets contained in $M\setminus\mathcal A$.}
\label{item:gh1ii}
        \item $\Lambda = \cup_{i = 1}^{k}\Lambda_i = NW(T) = \overline{\mathrm{Per}(T)}$ admits no cycles (see Definition~\ref{def:order})\label{item:gh1iii},
        \item there exists a unique basic set $\Lambda_i\subset M\setminus\mathcal A$ on which the topological pressure is maximal among the basic sets contained in $M\setminus\mathcal A$, and this pressure is realised by the unique equilibrium
state $\nu^g$ supported on $\Lambda_i$,
\label{item:gh1iv}
        \item $T:\Lambda_i\to \Lambda_i$ is topologically mixing.\label{item:ghv}
    \end{enumerate}
\end{hypothesis}

For this problem, we consider the (global) annealed transfer operator on $M \setminus \mathcal A$ for the $e^\phi$-weighted process $X_\e^\phi$ given by
\[\widehat{\calP}_\e \varphi(x) = e^{\phi(x)}\E_\e\left[\varphi \circ T_\w(x) \cdot \1_{M\setminus \mathcal{A}}\circ T_\w(x)\right].\]
The following theorem provides the existence of quasi-ergodic measures globally:

\begin{theorem}\label{thm:global-existence}
   Assume that $(T,g,\mathcal A)$ satisfies Hypothesis~\ref{hyp:global}. Let $\phi:M\to\mathbb R$ be such that $e^\phi$ is of class $\mathcal C^{r-1}$ on an open neighbourhood of every basic set $\Lambda_i\subset M\setminus\mathcal A$. Assume moreover that $e^\phi=g$ on each such basic set. Then, for all sufficiently small $\e>0$, the following holds:
   \begin{enumerate}[label = (\roman*)]
    \item \label{item:gex1} the operator $\widehat{\calP}_\e$ admits a unique (up to scaling) dominant eigenfunction $g_\e\in\mathcal C^0(M\setminus\mathcal A)$ such that $\widehat{\calP}_\e g_\e=\lambda_\e g_\e$,
    \item \label{item:gex2} the operator $\widehat{\calP}_\e^*$ admits a unique (up to scaling) dominant eigenmeasure $\mu_\e\in\mathcal M(M\setminus\mathcal A)$ such that $\widehat{\calP}_\e^*\mu_\e=\lambda_\e\mu_\e$, where $\lambda_\e=r(\widehat{\calP}_\e|_{\mathcal C^0(M\setminus\mathcal A)})$, and we normalise $g_\e$ and $\mu_\e$ so that $\mu_\e(g_\e)=1$,
    \item \label{item:gex3} $\nu_\e^\phi(\,\cdot\,)=\mu_\e(\,\cdot\,g_\e)$ is the unique quasi-ergodic measure for $X_\e^\phi$ on $M\setminus \mathcal A$ such that $\Lambda_i\subset \supp \nu_\e^\phi$, with $\Lambda_i$ as in~\ref{hyp:global}~\ref{item:ghv}.
\end{enumerate}
\end{theorem}

\noindent Global conditioned stochastic stability is then given by:
\begin{theorem}[Global conditioned stochastic stability, $\partial = \mathcal A$]\label{thm:global}
    For each $\e>0$ let $\nu^\phi_\e$ be the quasi-ergodic measure from Theorem~\ref{thm:global-existence}. Then $\nu_\e^\phi$ converges to $\nu^g$ of Hypothesis~\ref{hyp:global}~\ref{item:gh1iv} in {weak-$^*$} as $\e \to 0$.
\end{theorem}

This quasi-ergodic measure $\nu_\e^\phi$ is not necessarily unique. In the presence of finitely many basic sets $\Lambda_1, \dots, \Lambda_k \subset M$, the trajectory of an initial condition $x$ may never spend enough time sufficiently close to the most dominant hyperbolic set. In this case we would expect that $x \not \in \{g_\e > 0\}$. Based on the classical notion of \emph{dynamical filtrations}~\cite{Conley1978}, and under some additional assumptions (see Hypothesis~\ref{hyp:global-global}), we prove a more detailed version of Theorems~\ref{thm:global-existence} and~\ref{thm:global} in Section~\ref{sec:global-proof} to deal with such cases. Theorem~\ref{thm:global-global} provides, in particular, $t \leq k$ different quasi-ergodic measures $\nu_{\e, 1}^{\phi}, \dots,\nu_{\e, t}^{\phi},$ and for each $x \in W^s(\Lambda_j)$ for some $j \in \{1, \dots, k\}$, the quasi-ergodic limit in~\eqref{eq:qem-cond} holds for some particular $\nu_{\e, \ell}^{\phi},$ $\ell \in \{1, \dots, t\}$. Conditioned stochastic stability of these quasi-ergodic measures $\nu_{\e,\ell}^{\phi},$ $\ell \in \{1, \dots, t\}$ is also established in Theorem~\ref{thm:global-global}~\ref{item:global-global-global-stoch-stab}.

\subsection{Examples}\label{sec:examples}

\subsubsection{The H\'enon repeller}

Consider the paradigmatic family of H\'enon maps~\cite{Henon1976}
\begin{align*}
    H_{a,b}:\mathbb R^2&\to \mathbb R^2,\\
    (x,y)&\mapsto (1+y-a x^2,bx),
\end{align*}
with real parameters $a,b>0$. It is well-known~\cite{Feit1978,DevaneyNitecki1979} that if $b\in(0,1)$ and $a>(5+2\sqrt5)(1+b)^2/4$, then
$$
\Lambda_{a,b}
=
\left\{
(x,y)\in\mathbb R^2:\,
\sup_{n\in\mathbb Z}\|H_{a,b}^n(x,y)\|<\infty
\right\}
$$
is compact and uniformly hyperbolic and coincides with the non-wandering set of $H_{a,b}$. Moreover, $H_{a,b}:\Lambda_{a,b}\to\Lambda_{a,b}$ is topologically conjugate to the two-sided $2$-shift~\cite[Theorem, item~iv, p.~138]{DevaneyNitecki1979}.

Let $V\Subset\mathbb R^2$ be an isolating neighbourhood of $\Lambda_{a,b}$ and identify $\mathbb R^2$ with $\mathbb S^2\setminus\{p_\infty\}$. We may choose a $\mathcal C^r$ diffeomorphism $\widetilde H_{a,b}:\mathbb S^2\to\mathbb S^2$ which coincides with $H_{a,b}$ on a neighbourhood of $\overline V$. For $\omega\in[-\e,\e]^2$, let $\tau_\omega:\mathbb S^2\to\mathbb S^2$ be a smooth diffeomorphism which, in the chart $\mathbb S^2\setminus\{p_\infty\}\simeq\mathbb R^2$, coincides with the translation $x\mapsto x+\omega$ on a neighbourhood of $\widetilde H_{a,b}(\overline V)$, and set $T_\omega=\tau_\omega\circ\widetilde H_{a,b}.$
Thus, as long as the trajectory remains in $\overline V$,
$$
T_\omega(x)=H_{a,b}(x)+\omega.
$$
Hence the corresponding Markov process is precisely the additive random perturbation
$$
X_{n+1}^{\e}=H_{a,b}(X_n^\e)+\omega_n^\e,
\qquad
\omega_n^\e\sim_{\mathrm{i.i.d.}}\mathrm{Uniform}([-\e,\e]^2),
$$
conditioned upon remaining in $\overline V$. For $\e>0$ sufficiently small, Theorems~\ref{thm:qem-local} and~\ref{thm:local} apply and yield a quasi-ergodic measure $\nu_\e^\phi$ such that $\Lambda_{a,b}\subset\supp\nu_\e^\phi$. Moreover, $\nu_\e^\phi$ converges in the weak-$^*$ topology to the equilibrium state associated with the potential $\phi-\log|\det DH_{a,b}|_{E^u}|$ on $\Lambda_{a,b}$.

\subsubsection{Arnold's cat map}
Let $T$ denote Arnold's \emph{cat} map given by \cite{ArnoldAvez1967}
\[
\begin{split}
T: \mathbb T^2 &\to \mathbb T^2\\
x &\mapsto Ax \mod 1,
\end{split} \qquad A = \begin{pmatrix} 2&1\\1&1 \end{pmatrix}.
\]
$T$ is an Anosov diffeomorphism, i.e.~a diffeomorphism with hyperbolic structure on the tangent bundle and a special case of an Axiom A map, hence satisfies Hypotheses~\ref{hyp:local} and \ref{hyp:global} (with $k = 1$ in item~\ref{item:gh1iii}), as well as \ref{hyp:global-global}, together with any given H\"older potential, $\Lambda = M = \mathbb T^2$ and $\mathcal A = \emptyset$. Notice, again, that existence of equilibrium states $\nu^g$ for H\"older continuous potentials follows from~\cite[Theorem 4.1]{Bowen1975-book} (see also \cite[Theorem 7.5]{Baladi2018-book}).

Consider the Markov process $X^{\e,\phi}_n$ generated by small random perturbations of $T$ and weighted by $e^\phi$, as outlined in Section~\ref{sec:random_perturb}. For instance, we may set
\[X_{n+1}^\e = T(X_n^\e) + \w_n^\e \mod  1, \qquad n \geq 0,\]
with $\w_n^\e~\sim_{\text{i.i.d.}} \text{Uniform} ([-\e, \e]^2)$, $\e > 0$, and $X_0 \in \mathbb T^2$.

Alternatively, using the Blaschke factors from~\cite{SlipantschukBandtlow2017}, let $\omega=(\omega_0,\omega_1)\in\mathbb C^2$ and define
$$
T_\omega(z_1,z_2)=\left(\frac{z_1+\omega_0}{1+\overline{\omega_0}z_1}z_1z_2,\frac{z_1+\omega_1}{1+\overline{\omega_1}z_1}z_2\right).
$$
Then $T_{(0,0)}=T$, and for $\omega$ sufficiently small the maps $T_\omega$ are Anosov diffeomorphisms close to $T$. Moreover, since $\omega_0$ and $\omega_1$ act independently on the two circle coordinates, $\partial_\omega T_\omega$ is surjective. Hence, taking $\omega_n$ i.i.d.\ with full support in a sufficiently small neighbourhood of $0\in\mathbb C^2\simeq\mathbb R^4$, we obtain a random perturbation satisfying the assumptions of Section~\ref{sec:random_perturb}.

In both cases, for $\e$ small enough, Theorems~\ref{thm:qem-local} and \ref{thm:global-existence} provide the existence and uniqueness of quasi-ergodic measures $\nu_\e$ for $X_n^{\e, \phi}$, and Theorems~\ref{thm:local} and \ref{thm:global} establish that $\nu_\e$ converge to the unique equilibrium state $\nu^g$ in weak-$^*$ as $\e \to 0$.

\section{The anisotropic Banach space \texorpdfstring{$B^{t,s}$}{Bᵗˢ}}\label{sec:BanachSpaces}

This section contains a brief summary of Chapters 4 and 5 in the book ``Dynamical Zeta Functions and Dynamical Determinants for Hyperbolic Maps'' by Baladi~\cite{Baladi2018-book}. The construction of the anisotropic Banach spaces $B^{t,s}$ is essential to obtain bounds on the spectral and essential spectral radii of the deterministic transfer operator, which can be extended to random settings by well-established perturbation theory results of Keller and Liverani \cite{KellerLiverani1999}. Importantly, the spaces $B^{t,s}$ allow us to interpret the eigenfunctions of the dual transfer operator as (quasi-stationary) measures, which we leverage in the construction of the quasi-ergodic measure. In the process, we show a compact embedding property of these spaces (see Theorem~\ref{thm:compact_embedding}), which is leveraged in the proof of Theorem~\ref{thm:qem-local} and Theorem~\ref{thm:local} for vanishing Hölder weights (in Section~\ref{sec:QEMs-vanishing}) and in Lemma~\ref{lem:LY-Baladi}.

\begin{notation} We mostly follow the definitions in the underlying reference and the standard nomenclature used in the literature, where possible.
\begin{enumerate}[label = (\roman*)]
    \item $K\subset \mathbb R^d$ is a closed subset.
   \item $\mathcal C^\alpha(K)$ denotes the $\alpha$-H\"older space on $K$, while
$\mathcal C^\alpha_0(K)\coloneqq
\left\{\varphi\in\mathcal C^\alpha(\mathbb R^d):\operatorname{supp}\varphi\subset K\right\},$
endowed with the ambient $\mathcal C^\alpha(\mathbb R^d)$-norm. For a closed subset $K\subset M$, we use the same notation $\mathcal C^\alpha_0(K)$ for the space $
\left\{\varphi\in\mathcal C^\alpha(M):\operatorname{supp}\varphi\subset K\right\},$
endowed with the ambient $\mathcal C^\alpha(M)$-norm.
    \item $\mathcal F_b(Y) \coloneqq \{h:Y\to \R; \ h \text{ is bounded and measurable}\}$ for any domain $Y \subset M$.
    \item $\mathcal F$ denotes the Fourier transform and $\mathcal{F}^{-1}$ the inverse Fourier transform.
    \item $H^{t}_p(\R^d)$, with $t \in \R, 1<p<\infty$ denotes the (isotropic) Sobolev space defined by
    \[H^{t}_p(\R^d) \coloneqq \{\varphi \in \mathcal{S}' \,:\, \|\varphi\|_{H^{t}_p(\R^d)} \coloneqq \|\mathcal F^{-1}\left((1 +\|\xi\|^2)^{t/2}\mathcal F(\varphi)(\xi)\right)\|_{L_p(\R^d)} < \infty\},\]
    where $\|\cdot\|_{L_p(\R^d)}$ is the usual norm and $\mathcal S'$ the (Schwartz) space of tempered distributions. Recall that if $\ell \in \N$, then \[H^{\ell}_p(\R^d) \equiv \{\varphi \in L_p(\R^d)\,: \,  \sum_{|\beta|  \leq \ell} \|\partial^\beta \varphi\|_{L_p(\R^d)} < \infty\}\] for all $1<p<\infty$.
    \item Given $t>0$, we write $H^t$ for the $L_2$ Sobolev space of order $t$, i.e.~$H^t = H^t_2$ (see \cite[Chapter~2, Section~2.2.1]{Baladi2018-book}). Recall that its dual is $H^{-t}$. If $\varphi \in H^{-t}$ and $f\in H^t$ we define $\langle \varphi,f\rangle := \varphi (f).$
    \end{enumerate}
\end{notation}

\subsection{Charts and cone systems adapted to \texorpdfstring{$(T,V)$}{(T,V)}} \label{sec:charts-cone-systems}

A cone in $\R^d$ is an invariant subset under scalar multiplication. For two cones $\mathbf C$ and $\mathbf C'$ in $\mathbb R^d$, we write 
$$\mathbf C \Subset \mathbf C'\ \text{if}\ \overline{\mathbf C}\subset \mathrm{Int}(\mathbf C ')\cup \{0\}.$$

We say that a cone $\mathbf C$ is $d'$-dimensional if $d'\geq 1$  is the maximal dimension of a linear subset of $\mathbf C.$

\begin{definition}[{Cone system, $\Theta<\Theta'$, \cite[Definition 4.10]{Baladi2018-book}}] Let $\Cp$ and $\Cm$ be closed cones in $\R^d$ with non-empty interiors, of respective dimensions $d_s$ and $d_u$, and such that $\Cp \cap \Cm =\{0\}$ (i.e.~the cones are transversal). Let $\Phip:\mathbb S^{d-1}\to [0,1]$ be a $\mathcal C^\infty$ function on the unit sphere $\mathbb{S}^{d-1}$ in $\R^d$ satisfying
\[ \Phi_{+}(\xi)=
\begin{cases}
1,&\xi\in \mathbb S^{d-1}\cap \Cp \\ 
0,& \xi\in \mathbb S^{d-1}\cap \Cm.
\end{cases} \]
Define $\Phim:\mathbb S^{d-1}\to [0,1]$ by $\Phi_{-}(\xi) = 1- \Phi_{+}(\xi).$ A quadruple $\Theta = (\Cp,\Cm, \Phip,\Phim)$ is called a \emph{cone system}. For another such quadruple $\Theta' = (\Cp', \Cm',\Phip', \Phim'),$ we write $\Theta <\Theta'$ if  $\R ^d \setminus \mathbf C'_+ \Subset \Cm.$ In particular, this implies $\Cp \Subset \Cp'$ and $\Cm'\Subset \Cm$. 
\end{definition}

\begin{definition}[{Cone-hyperbolic diffeomorphism, \cite[Definition 4.11]{Baladi2018-book}}]\label{def:cone-hyp}
    Let $U$ be an open and bounded set in $\R^d$, let $\Theta = (\Cp,\Cm, \Phip,\Phim)$ and $\Theta' = (\Cp', \Cm',\Phip', \Phim')$ be two cone systems. A $\mathcal C^r$ diffeomorphism $F:U\to \mathbb R^d$ onto its image is \emph{cone-hyperbolic} from $\Theta$ to $\Theta'$ if $F$ extends to a bilipschitz $\mathcal C^1$ diffeomorphism of $\R^d$ such that
$$D F^{\intercal}_x (\R^d \setminus \Cp) \Subset \mathbf C'_-\ \text{for every }x\in \R^d,$$
where $A^\intercal$ denotes the transpose of the matrix $A$.
\end{definition}

It is important to remark that we view the cones in the cotangent bundle $T^*\mathbb R^d$ so that $F$ acts on these with the transpose of $DF$. In some sense, the definition of cone-hyperbolicity ``from $\Theta$ to $\Theta'$'' concerns the direction for the cotangent dynamics, which are in the inverse direction of the dynamics for $F$. 

\begin{definition}[{Charts and partition of unity adapted to \texorpdfstring{$(T,V)$}{(T,V)}, \cite[Definition 4.14]{Baladi2018-book}}]\label{def:charts}Fix a finite system of $\mathcal C^\infty$ local charts $\{(V_i,\kappa_i)\}_{i\in I},$ with open subsets $V_i\subset M$ and maps $\kappa_i: V_i \to \R^d$ such that $V \subset \cup_{i} V_i$ and
\begin{enumerate}[label = (\arabic*)]
    \item $\mathcal V =\{V_i\}_{i\in I}$ is a generating cover\footnote{Recall that a cover is
    \emph{generating} for $T$ if the diameter of all sets $\{\cap_{k = -m}^{m-1} T^{-k}(V_{\omega_k}): \, \omega \in I^{2m}\}$ tends to zero as $m \to \infty$.}
    of $V$, and there is no strict sub-cover.
    \item $\kappa_i(V_i)$ is a bounded open subset of $\R^d$ with smooth boundary for each $i\in I$.
\end{enumerate}

Finally, let $\{\theta_i\}_{i\in I}$ be a $\mathcal C^\infty$ finite partition of unity for $V$ subordinate to the cover $\mathcal V$, that is, the support of each $\theta_i: M\to [0,1]$ is contained in $V_i$, and we have $\sum_{i\in I}\theta_i(x) = 1$ for every $x\in V.$
\end{definition}

\begin{definition}[{Cone systems adapted to $(T,V)$, \cite[Definition 4.15]{Baladi2018-book}}]\label{def:cone-adapted} If $\Lambda$ is a hyperbolic basic set for $T$ with isolating neighbourhood $V$, we may choose\footnote{Such a choice can always be made; we refer to \cite[Section 4 and Definition 4.15]{Baladi2018-book} for details.} a finite family of cone systems $\{\Theta_i = (\mathbf C_{i,+},\mathbf C_{i,-}, \Phi_{i,+},\Phi_{i,-})\}_{i\in I}$, where $I$ is the index set for the local charts from Definition~\ref{def:charts}, so that  the following conditions hold:
\begin{enumerate}[label = (\arabic*)]
    \item[(3)] If $x\in V_i\cap \Lambda$, the cone $(D\kappa_i)^*_x(\mathbf C_{i,+})$ contains the ($d_s$-dimensional) normal subspace of $E^u(x)$, and the cone  $(D\kappa_i)^*_x(\mathbf C_{i,-})$ contains the ($d_u$-dimensional) normal subspace of $E^s(x).$
    \item[(4)] If $V_{ij}\coloneqq T^{-1}(V_j)\cap V_i\neq \emptyset$, the map in charts
    $$ F = T_{ij} = \kappa_{j}\circ T\circ \kappa_{i}^{-1}:\kappa_{i}(V_{ij})\to \kappa_{j}(V_j)$$is a $\mathcal C^r$ cone-hyperbolic diffeomorphism from $\Theta_j$ to $\Theta_i$.
\end{enumerate}
\end{definition}

\subsection{The space \texorpdfstring{$B^{\Theta,t,s}$}{Bᵗˢ} in \texorpdfstring{$\R^d$}{Rᵈ} } \label{sec:Bts-in-Rd}

We follow the presentation given in \cite[Section~2.4.1]{Baladi2018-book}. Let $\chi:\R \to [0,1]$ be a symmetric $\mathcal C^\infty$ bump-function such that 
$$\chi(x) =\begin{cases}
    1,&|x|\leq 1,\\
    0, & |x| \geq 2.
\end{cases} $$
Define the functions $\psi_n: \R^d \to \R, n \in \N_0$, such that $\psi_0(\xi) = \chi(\|\xi\|)$ and
$$\psi_n(\xi) = \chi(2^{-n}\|\xi\|) - \chi(2^{-n+1}\|\xi\|), \text{ for } n \geq 1.$$ 
Given a cone system $\Theta = (\Cp,\Cm,\Phip,\Phim)$ in $\R^d$, let $(n,\sigma) \in (\N_0)\times \{-,+\}$. We define for $\xi \in \mathbb R^d:$
$$
\psi_{\Theta,0,\sigma}(\xi) \coloneqq \frac{\chi(\|\xi\|)}{2}, \quad\text{ and }\quad \psi_{\Theta,n,\sigma}(\xi) \coloneqq \Phi_\sigma\left(\frac{\xi}{\|\xi\|}\right) \psi_n(\xi), \text{ for }\ n\geq 1.
$$
Observe that $$\sum_{\sigma\in\{+,-\}}\sum_{n=0}^{\infty} \psi_{\Theta,n,\sigma} = 1.$$

\begin{notation} Given a $\mathcal C^\infty$ function $\varphi:\mathbb R^d \to \mathbb C$ with compact support, we define $\varphi_{\Theta,n,\sigma}$ to be
$$\varphi_{\Theta,n,\sigma} \coloneqq \widehat{\psi}_{\Theta,n,\sigma}\ast \varphi \coloneqq \int_{\R^d}\widehat{\psi}_{\Theta,n,\sigma}(x-y)\varphi(y)\d y,$$
where 
$$\widehat{\psi}_{\Theta,n,\sigma}(x) \coloneqq \mathcal F^{-1} {\psi}_{\Theta,n,\sigma} (x)\coloneqq \frac{1}{(2\pi)^d}\int_{\R^d} e^{i x \xi} \psi_{\Theta,n,\sigma}(\xi)\d \xi,$$
with $x\xi$ denoting the scalar product. In particular, we have that $\sup_{(n, \sigma)} \|\widehat{\psi}_{\Theta, n , \sigma}\|_{L^1(\mathbb R^d)} < \infty.$
\end{notation}

\begin{definition}[{Fake unstable leaves $\mathcal F(\Theta)$ and the $L_1(\mathcal F)$-norm, \cite[Definition 5.10]{Baladi2018-book}}] Let $\Theta = (\Cp,\Cm, \Phip,\Phim)$ be a cone system. We define the set of unstable leaves $\mathcal F= \mathcal F(\Theta)$ to be the set of $\mathcal C^1$-submanifolds $\Gamma \subset \R^d$, of dimension $d_u$, such that the straight line connecting any two points in $\Gamma$ is normal to a $d_s$-dimensional subspace contained in $\Cp$, i.e.
$$\sup_{\substack{x,y\in \Gamma\\x\neq y}}\inf_{\substack{C\subset \Cp\\d_s-\text{subspace}}}\sup_{c\in C\setminus \{0\}}\left\langle  \frac{y-x}{\|y-x\|} ,\frac{c}{\|c\|} \right\rangle  = 0.$$ 
For $\varphi \in \mathcal C^\infty (\R^d)$ and $\mathcal F = \mathcal F(\Theta)$ we introduce the $L_1(\mathcal F)$-norm
$$\|\varphi\|_{L_1(\mathcal F)} = \sup_{\Gamma\in\mathcal F}\|\varphi\|_{L_1(\mu_\Gamma)}, $$
where $\mu_\Gamma$ is the Riemannian volume on $\Gamma$ induced by the standard metric on $\R^d$.
\end{definition}

\begin{definition}[{The local space $B^{\Theta, t,s}(K)$, \cite[Definition 5.12]{Baladi2018-book}}]
\label{def:localbts}
Given $t,s \in \mathbb R$, a cone system $\Theta$, and a compact set $K \subset \R^d$, we define 
\[
\|\varphi\|_{B^{\Theta, t,s}} \coloneqq \max\left\{ \sup_{n\geq 0} 2^{tn}\|\varphi_{\Theta,n,+}\|_{L_1(\mathcal F)}, \sup_{n\geq 0} 2^{sn}\|\varphi_{\Theta,n,-}\|_{L_1(\mathcal F)}\right\}
\]
and set $B^{\Theta, t,s}$ to be the completion of $\mathcal C^u_0(K)$, for any fixed $u>t$, with respect to the $B^{\Theta,t,s}$-norm, $\|\cdot\|_{B^{\Theta,t,s}}$.
\end{definition}

\begin{proposition}[{\cite[Lemma 5.14]{Baladi2018-book}}]\label{prop:inj} Assume $s\leq t$. For any $u>\max\{0,t\}$, there exists a constant $C = C(u,K)$ such that $\|\varphi\|_{B^{\Theta, t, s}}\leq C \|\varphi\|_{\mathcal C^u_0}$ for every $\varphi \in \mathcal C^{\infty}(\R^d ; \C)$ supported in $K$. Moreover, for any $v > \max\{-s, 0\}$, the space $B^{\Theta, t, s}(K)$ is contained in the space of distributions of order $v$ supported in $K$. 

In other words, $\mathcal C^u_0(K) \subset B^{\Theta, t,s} \subset (\mathcal C^v_0(K))^*$ if $u> \max\{0,t\}$ and $v > \max\{-s, 0\}$.
\end{proposition}

\begin{notation}
Let $D\subset\mathbb R^d$ be a compact domain with smooth boundary. For an integer $\ell\geq0$, we denote by $
H_0^\ell(D)\coloneqq\overline{\mathcal C_c^\infty(\operatorname{int}D)}^{\,H^\ell(D)}$
the zero-boundary Sobolev space and set $H^{-\ell}(D)\coloneqq\bigl(H_0^\ell(D)\bigr)^*.$
\end{notation}

\begin{corollary}\label{cor:compemb}
Assume $s\leq t$. Let $K\subset\mathbb R^d$ be compact and let $\widetilde K\subset\mathbb R^d$ be a compact domain with smooth boundary such that
$
K\subset\operatorname{int}(\widetilde K).
$
Then, for some integer $\ell>0$, there is a compact embedding
$
B^{\Theta,t,s}(K)\hookrightarrow H^{-\ell}(\widetilde K).
$
\end{corollary}

\begin{proof}
Choose $v>\max\{-s,0\}$. By Proposition~\ref{prop:inj}, every element of $B^{\Theta,t,s}(K)$ defines a distribution of order $v$ supported in $K$. Since $
K\subset\operatorname{int}(\widetilde K),
$ we obtain a continuous embedding $B^{\Theta,t,s}(K)\hookrightarrow
\bigl(\mathcal C_0^v(\widetilde K)\bigr)^*.$

Choose $\eta>0$ and an integer $\ell>v+\eta+d/2$. By the general Sobolev inequalities
\cite[Section~5.6.3, Theorem~6]{Evans2010}, $
H_0^\ell(\widetilde K)\hookrightarrow
\mathcal C_0^{v+\eta}(\widetilde K)$ continuously. Since $\widetilde K$ is compact and $\eta>0$, the embedding $
\mathcal C_0^{v+\eta}(\widetilde K)\hookrightarrow
\mathcal C_0^v(\widetilde K)$
is compact by the Arzel\`a--Ascoli theorem. Hence, $
H_0^\ell(\widetilde K)\hookrightarrow
\mathcal C_0^v(\widetilde K)$
is a compact embedding.

Taking duals and using the fact that the adjoint of a compact operator is compact, we obtain a compact linear operator
$$
\bigl(\mathcal C_0^v(\widetilde K)\bigr)^*
\longrightarrow
\bigl(H_0^\ell(\widetilde K)\bigr)^*
=
H^{-\ell}(\widetilde K).
$$
Its composition with the continuous inclusion of $B^{\Theta,t,s}(K)$ is injective. Indeed, if an element of $B^{\Theta,t,s}(K)$ vanishes in $H^{-\ell}(\widetilde K)$, then it vanishes on $\mathcal C_c^\infty(\operatorname{int}\widetilde K)$ and hence is the zero distribution. Therefore
$
B^{\Theta,t,s}(K)\hookrightarrow H^{-\ell}(\widetilde K)
$
is a compact embedding.
\end{proof}

\subsection{The space \texorpdfstring{$B^{t,s}$}{Bᵗˢ} on \texorpdfstring{$M$}{M}}\label{sec:Bts-on-M} The ideas developed so far can be adapted to $M$ and the map $T$ by means of suitable local charts. 
\begin{definition}[{Regular cone-hyperbolicity, \cite[Definition 5.10]{Baladi2018-book}}]\label{def:regular_conehyp}
    Let $U$ be a bounded open subset in $\R^d$, and let $\Theta = (\Cp,\Cm, \Phip, \Phim)$ and $\Theta' = (\Cp',\Cm', \Phip', \Phim')$ be two cone systems. A $\mathcal C^r$ diffeomorphism onto its image $F:U \to \R^d$ is \emph{regular cone-hyperbolic from $\Theta$ to $\Theta'$} if $F$ is cone-hyperbolic in the sense of Definition~\ref{def:cone-hyp} and, in addition, there exists, for each $x, y \in U$, a linear transformation $\mathbb L_{xy}$ satisfying
    \[\mathbb L _{x,y}^{\intercal}(\R^d \setminus \Cp)\Subset \Cm' \text{ and } \mathbb L _{xy}(x-y) = F(x) - F(y).\]
\end{definition}

Notice that if $F$ is regular cone-hyperbolic, then the extension of $F$ to $\R^d$ maps each element of $\mathcal F(\Theta')$ to an element of $\mathcal F(\Theta)$.

\begin{definition}[Anisotropic spaces $B^{t,s}$ on $M$]
Fix real numbers $s$ and $t$. The Banach space $B^{t,s}=B^{\Theta,t,s}(T,V)$ is the completion of $\mathcal C_c^\infty(V)$ for the norm
$$
\|\varphi\|_{B^{t,s}(T,V)}
\coloneqq
\max_{i\in I}
\|(\theta_i\cdot\varphi)\circ\kappa_i^{-1}\|_{B^{\Theta_i,t,s}},
$$
where $\{\kappa_i:V_i\to\mathbb R^d\}_{i\in I}$ is a finite $\mathcal C^\infty$ atlas of $M$, $\{\theta_i\}_{i\in I}$ is a partition of unity and $\Theta=\{\Theta_i\}_{i\in I}$ is a family of cone systems satisfying the requirements of Definitions~\ref{def:charts} and~\ref{def:cone-adapted}, with item~(4) of Definition~\ref{def:cone-adapted} strengthened to require regular cone-hyperbolicity as in Definition~\ref{def:regular_conehyp}.
\end{definition}

Observe that the symbol $\Theta$ refers to a single cone system when considering $B^{\Theta, t, s}(K)$ for some compact $K \subset \R^d$ but refers to a finite family of cone systems, satisfying the regular cone-hyperbolic conditions of Definitions~\ref{def:cone-hyp} and \ref{def:regular_conehyp} when considering $B^{\Theta, t, s}(T,V)$. Moreover, we may drop the symbol $\Theta$ and simply write $B^{t,s}$ when a statement holds for any family of cone systems. 

 \begin{proposition}[Leibniz bounds on $B^{t,s}(T,V)$] \label{prop:items}
 Consider $r>1$ and $t-(r-1)<s<t$, a generating open cover $\mathcal V = \{V_i\}_{i \in I}$ of $M$, a finite smooth atlas $\{\kappa_i:V_i \to \R^d\}_{i\in I}$ of $M$, a partition of unity $\{\theta_i\}_{i\in I}$, and two finite families of cone systems in $\R^d$, $\Theta = \{\Theta_i\}_{i\in I}$ and $\Theta'= \{\Theta_i'\}_{i\in I}$, such that $\Theta_i' < \Theta_i$ for every $i\in I$, all satisfying the conditions of Definition~\ref{def:charts}. Then there exists $C= C(\Theta,\Theta')$ such that
 \[\|\varphi\|_{B^{\Theta',t,s}(T,V)} \leq C \|\varphi\|_{B^{\Theta,t,s}(T,V)}\]
 for every $\varphi\in B^{\Theta,t,s}(T,V).$
 \end{proposition}
\begin{proof}
For each $i\in I$, choose a $\mathcal C^\infty$ function $\widetilde{\theta}_i:M\to [0,1]$ such that 
$$\widetilde{\theta}_i(x) = \begin{cases}
    1,&x\in \mathrm{supp}(\theta_i) \Subset V_i,\\
    0,&x\in M\setminus V_i.
\end{cases}$$
From \cite[Corollary~5.19]{Baladi2018-book} we obtain that for every $i\in I$, there exists $C_i= C_i(\Theta,\Theta')$ such that 
\begin{align*}
   \|\left(\theta_i \cdot \varphi\right)\circ \kappa_i^{-1} \|_{B^{\Theta_i',t,s}(\overline{V_i} )}
   &=\|(\widetilde{\theta}_i \circ \kappa_i^{-1} )\cdot \left(\theta_i \cdot \varphi\right)\circ \kappa_i^{-1} \|_{B^{\Theta_i',t,s}}\\
   &\leq C_i\|\widetilde{\theta}_i\circ \kappa_i^{-1}  \|_{\mathcal C^{r-1}}  \|\left(\theta_i \cdot\varphi\right)\circ \kappa_i^{-1} \|_{B^{\Theta_i,t,s}(\overline{V_i})}.
\end{align*}
Set $ C = \max_{i\in I} C_i\|\widetilde{\theta}_i \circ \kappa_i^{-1}  \|_{\mathcal C^{r-1}}$.
\end{proof}

Below, we state the main result of this section, which establishes a compact embedding of the anisotropic scales $B^{t,s}$.

\begin{theorem}[{Compact embedding of $B^{t,s}$}]\label{thm:compact_embedding}
   Consider two finite families of cone systems in $\mathbb R^d$; $\Theta' = \{\Theta_i'\}_{i\in I}$ and $\Theta =\{\Theta_i\}_{i\in I}$, such that $\Theta_i < \Theta_i'$ for all $i \in I$. Let $s,s',t,t' \in \mathbb R$ be such that $s\leq t < t'$ and $s< s'\leq t'$. Then $B^{\Theta',t',s'}$ compactly embeds in $B^{\Theta, t,s}.$
\end{theorem} 
\begin{proof}
Observe that it is enough to show that there is a compact embedding $B^{\Theta',t',s'}(K)\hookrightarrow B^{\Theta,t,s}(K)$ for a compact domain $K\subset \R^d$ with smooth boundary.
Given $\varepsilon>0$, take $N = N(\varepsilon)$ such that $1/2^{n\min\{t'- t,s'-s\}} <\varepsilon$ for all $n \geq N$. From Proposition~\ref{prop:items} (see also \cite[Corollary~5.19]{Baladi2018-book}) and the definition of the ${B^{\Theta,t,s}}$-norm we obtain that for every $\varphi \in B^{\Theta',t',s'}(K)$,
\begin{align*}
\|\varphi\|_{B^{\Theta,t,s}(K)}&\leq \max_{\substack{n\in \{0,\ldots,N-1\}}}\left\{2^{t n}\|\varphi_{\Theta,n,+}\|_{L_1(\mathcal F)},2^{sn}\|\varphi_{\Theta,n,-}\|_{L_1(\mathcal F)}\right\}+ \frac{1}{2^{N \min\{t'- t,s'-s\}}}\|\varphi\|_{B^{\Theta,t',s'}}\\
&\leq C_\varepsilon \max_{\substack{n\in \{0,\ldots,N\}\\ \sigma \in\{-,+\}}}\{\|\varphi_{\Theta,n,\sigma}\|_{L_1(\mathcal F)}\} + \varepsilon C \|\varphi\|_{B^{\Theta',t',s'}(K)},
\end{align*}
for some constants $C,C_\e >0$.
To prove the theorem, it is enough to show that for every $n\in \mathbb N_0$ and $\sigma\in\{-,+\}$, the map
\begin{align*}
F_{n,\sigma}:B^{\Theta,t,s}(K)&\to L_1(\mathcal F)\\*
\varphi&\mapsto \varphi_{\Theta,n,\sigma}\coloneqq\widehat{\psi}_{\Theta,n,\sigma}\ast\varphi
\end{align*}
is a compact linear operator. Choose a compact domain $\widetilde K\subset\mathbb R^d$ with smooth boundary such that $K\subset\operatorname{int}(\widetilde K)$. By Corollary~\ref{cor:compemb}, we may choose $\ell\in\mathbb N$ such that the inclusion
$$
\iota:B^{\Theta,t,s}(K)\to H^{-\ell}(\widetilde K)
$$
is compact. Choose $\overline\theta\in\mathcal C_c^\infty(\operatorname{int}\widetilde K)$ equal to $1$ on a neighbourhood of $K$, and define
\begin{align*}
G_{n,\sigma}:H^{-\ell}(\widetilde K)&\to L_1(\mathcal F)\\*
\varphi&\mapsto\left(x\mapsto\left\langle\varphi,\overline\theta(\,\cdot\,)\widehat{\psi}_{\Theta,n,\sigma}(x-\cdot)\right\rangle\right).
\end{align*}
Since every element of $B^{\Theta,t,s}(K)$ is supported in $K$ and $\overline\theta=1$ on a neighbourhood of $K$, we have $F_{n,\sigma}=G_{n,\sigma}\circ\iota$. Thus, it suffices to show that $G_{n,\sigma}$ is bounded.
Recall that $\widehat{\psi}_{\Theta,n,\sigma}\in\mathcal S(\mathbb R^d)$ for every fixed $(n,\sigma)\in\mathbb N_0\times\{-,+\}$. For every $\varphi\in H^{-\ell}(\widetilde K)$, we obtain
\begin{align}\label{eq:bound_Hell}
|G_{n,\sigma}\varphi(x)|&=
\left|
\left\langle
\varphi,\overline\theta(\,\cdot\,)
\widehat{\psi}_{\Theta,n,\sigma}(x-\cdot)
\right\rangle
\right|\leq\|\varphi\|_{H^{-\ell}(\widetilde K)}\|\overline\theta(\,\cdot\,)\widehat{\psi}_{\Theta,n,\sigma}(x-\cdot)\|_{H_0^\ell(\widetilde K)}\\
&\leq C\|\varphi\|_{H^{-\ell}(\widetilde K)}\|\overline\theta(\,\cdot\,)\widehat{\psi}_{\Theta,n,\sigma}(x-\cdot)\|_{H^\ell(\mathbb R^d)}.
\end{align}
By the equivalence between the Fourier-defined Sobolev norm and the derivative norm, there exists $C_0>0$ such that
\begin{equation}\label{eq:bound_theta_psi}
\|\overline\theta(\cdot)\widehat{\psi}_{\Theta,n,\sigma}(x-\cdot)\|_{H^\ell(\R^d)}^2
\leq C_0\sum_{|\beta|\leq\ell}\int_{\R^d}\left|D_y^\beta\left(\overline\theta(y)\widehat{\psi}_{\Theta,n,\sigma}(x-y)\right)\right|^2\d y.
\end{equation}
Fix an integer $k>d_u$. Since $\overline\theta$ has compact support and $\widehat{\psi}_{\Theta,n,\sigma}$ is a Schwartz function, the Leibniz rule gives
$$
\left|D_y^\beta\left(\overline\theta(y)\widehat{\psi}_{\Theta,n,\sigma}(x-y)\right)\right|\leq\frac{C_1}{1+\|x\|^k},\qquad |\beta|\leq\ell,
$$
uniformly in $x,y\in\R^d$. All these derivatives are supported, as functions of $y$, in the fixed compact set $\supp\overline\theta$. Thus, increasing $C_1$ if necessary, \eqref{eq:bound_theta_psi} implies
\begin{equation}\label{eq:compact-embedding-kernel-decay}
\|\overline\theta(\cdot)\widehat{\psi}_{\Theta,n,\sigma}(x-\cdot)\|_{H^\ell(\R^d)}
\leq\frac{C_1}{1+\|x\|^k}.
\end{equation}
We next establish the required bound uniformly over $\Gamma\in\mathcal F(\Theta)$. Fix a $d_u$-dimensional subspace $E\subset\Cm$ and let $\pi:\R^d\to E$ be the orthogonal projection. Since $\Gamma\in \mathcal F(\Theta)$ its definition implies that $\pi|_\Gamma$ is injective (see \cite[Remark~4.1]{Baladi-Tsujii2008}). Moreover, there exists $C_2\geq1$, independent of $\Gamma$, such that
\begin{equation}\label{eq:compact-embedding-leaf-projection}
\|x-y\|\leq C_2\|\pi(x)-\pi(y)\|,\qquad x,y\in\Gamma.
\end{equation}
Indeed, for any two distinct points $x,y\in\Gamma$, the vector $x-y$ is orthogonal to some $d_s$-dimensional subspace $L\subset\Cp$. Since $\Cp\cap\Cm=\{0\}$ and $E\subset\Cm$, we have $L\oplus E=\R^d$, so the restriction of $\pi$ to $L^\perp$ is invertible. Compactness of the family of such subspaces $L$ gives a uniform bound on these inverses, proving \eqref{eq:compact-embedding-leaf-projection}.
Hence, $\pi|_\Gamma$ is a $\mathcal C^1$ diffeomorphism onto its image, and the volume Jacobian of its inverse is bounded by $C_2^{d_u}$. Using that $\|\pi(x)\|\leq\|x\|$, the change-of-variables formula gives
\begin{align}\label{eq:compact-embedding-uniform-leaf-bound}
\sup_{\Gamma\in\mathcal F}\int_\Gamma\frac{1}{1+\|x\|^k}\,\d\mu_\Gamma(x)
&\leq C_2^{d_u}\sup_{\Gamma\in\mathcal F}\int_{\pi(\Gamma)}\frac{1}{1+\|u\|^k}\,\d u
\leq C_2^{d_u}\int_E\frac{1}{1+\|u\|^k}\,\d u
=:C_3<\infty,
\end{align}
where finiteness follows from $k>d_u$.
Combining \eqref{eq:bound_Hell}, \eqref{eq:compact-embedding-kernel-decay} and \eqref{eq:compact-embedding-uniform-leaf-bound}, we obtain
\begin{align}\label{eq:compact-embedding-convolution-bound}
\|G_{n,\sigma}\varphi\|_{L_1(\mathcal F)}
&\leq C\|\varphi\|_{H^{-\ell}(\widetilde K)}
\sup_{\Gamma\in\mathcal F}\int_\Gamma\|\overline\theta(\cdot)\widehat{\psi}_{\Theta,n,\sigma}(x-\cdot)\|_{H^\ell(\R^d)}\,\d\mu_\Gamma(x)
\leq CC_1C_3\|\varphi\|_{H^{-\ell}(\widetilde K)}.
\end{align}
Thus $G_{n,\sigma}$ is bounded, and hence $F_{n,\sigma}=G_{n,\sigma}\circ\iota$ is compact.
Finally, let $(\varphi_j)_{j\geq1}$ be a bounded sequence in $B^{\Theta',t',s'}(K)$. Since each $F_{n,\sigma}$ is compact and $\mathbb N_0\times\{-,+\}$ is countable, there exists a subsequence, still denoted by $(\varphi_j)_{j\geq1}$, such that $(\varphi_j)_{\Theta,n,\sigma}$ converges in $L_1(\mathcal F)$ for every fixed pair $(n,\sigma)$. Fix $\varepsilon>0$ and choose $N=N(\varepsilon)$ as above. Applying the first estimate to $\varphi_j-\varphi_m$ gives
$$
\|\varphi_j-\varphi_m\|_{B^{\Theta,t,s}(K)}
\leq C_\varepsilon\max_{\substack{0\leq n\leq N\\ \sigma\in\{-,+\}}}\|(\varphi_j-\varphi_m)_{\Theta,n,\sigma}\|_{L_1(\mathcal F)}
+\varepsilon C\|\varphi_j-\varphi_m\|_{B^{\Theta',t',s'}(K)}.
$$
The first term tends to zero as $j,m\to\infty$, while the second is bounded by a constant times $\varepsilon$. Since $\varepsilon>0$ is arbitrary, the subsequence is Cauchy in $B^{\Theta,t,s}(K)$, which proves the desired compact embedding.
\end{proof}

\begin{remark}
Theorem~\ref{thm:compact_embedding} does not appear to have been previously established in the literature (see e.g.~\cite{Baladi-Tsujii2007, Baladi-Tsujii2008, Baladi2018-book}) even though it is fundamental to our subsequent arguments. This may be due to quasi-compactness not being derived via a standard Lasota--Yorke-type inequality in the cited references. Instead, the approach therein often relies on decomposing the transfer operator $\mathcal{L}_g$ as $\mathcal{L}_g = \mathcal{L}_c + \mathcal{L}_b,$ where $\mathcal{L}_c : \mathcal{B}^{t,s} \to \mathcal{B}^{t,s}$ is a compact operator and $\mathcal{L}_b : \mathcal{B}^{t,s} \to \mathcal{B}^{t,s}$ is a bounded operator with spectral radius strictly smaller than that of $\mathcal{L}_g$, i.e.~$r(\mathcal{L}_b) < r(\mathcal{L}_g)$.
\end{remark}

\section{Properties of \texorpdfstring{$\calP_g$ on $B^{t,s}$ and local results for vanishing weights}{Pg on Bᵗˢ and QEMs for vanishing weights}}
\label{sec:QEMs-vanishing}

We dedicate this section to the study of the deterministic transfer operator $\calP_g$ and the annealed operator $\calP_\e$ acting on $B^{t,s}$ for H\"older weights supported in $V$. Throughout this section, $\calP_\e$ denotes the annealed transfer operator associated with the same weight $g$ as $\calP_g$, namely
$$
\calP_\e f(x)=g(x)\E_\e\left[f\circ T_\omega(x)\cdot\1_{\overline V}\circ T_\omega(x)\right].
$$
In particular, we present a strategy to build quasi-ergodic measures in Lemma~\ref{lem:spectral_gap_QEM} based on the presence of a spectral gap in $\mathcal C^0$ and check that the required assumptions are satisfied.

Let us first recall some standard definitions from operator theory. 

\begin{definition}[Quasi-compact operator, spectral gap]
Let $(\mathcal B,\|\cdot\|) $ be a Banach space and $\mathcal  L:\mathcal B\to \mathcal B$ be a bounded linear transformation. We say that $\mathcal L$ is \emph{quasi-compact} if there exist $\mathcal L$-invariant sub-spaces $F, W\subset \mathcal B$ such that $F \oplus W = \mathcal B$, $F$ is finite-dimensional and $r(\left.\mathcal L\right|_{W}) < r(\mathcal L)$.

We say that $\mathcal L$ has a \emph{spectral gap} if $\mathcal L$ is a quasi-compact operator and there exists only one element in the peripheral spectrum of $\mathcal L$, and such an element is a simple eigenvalue.
\end{definition}

\begin{lemma}\label{lem:spectral_gap_QEM}
    Let $(\Omega,\mathbb P)$ be a probability space and $\{T_\omega:M\to M\}_{\omega\in\Omega}$ be a family of maps. For $Y \subset M$ a compact subset and a non-negative function $\psi:Y \to \mathbb R$, $\psi \in \mathcal C^\alpha(Y)$, $\alpha > 0$, assume that the operator
    \[\calP: f(x) \mapsto \psi(x) \E_\e[f\circ T_\w(x) \cdot \1_{Y}\circ T_\w(x)]\]
    is strong Feller and has a spectral gap in $\mathcal C^0(Y)$. Let $g \in \mathcal C^0(Y)$ be the dominant eigenfunction such that $\calP g = \lambda g$ with $\lambda = r(\calP)$, and let $\mu \in \mathcal M(Y)$ satisfy $\calP^* \mu = \lambda \mu$. Assume that $\mu(g) =1$.
Then the following holds:
\begin{enumerate}
    \item The measure $\nu \in \mathcal M(Y)$, defined by $\nu(h)\coloneqq \mu( h \cdot g)$ for every $h\in\mathcal F_b(Y)$, bounded and measurable, is a quasi-ergodic measure for the process $X^{\log\psi}$ conditioned upon staying in $Y$.
    \item For every $x\in \{g\neq 0\}$ and every $h:Y\to\mathbb R$ bounded and measurable we have that
    $$\lim_{n\to\infty}\E_x^{\log \psi} \left[ \frac{1}{n} \sum_{i=0}^{n-1} h \circ X_i^{\log\psi} \, \bigg |\, \tau >n\right]  = \int_{Y} h \d \nu.$$
    \item Assume that $(\Omega,\mathbb P) = (\Omega_\e,\mathbb P_\e)$ and let $T_\omega$ be defined as in Section~\ref{sec:random_perturb}. Then there exists $m \in \mathcal C^0(Y)$ such that $\mu = m(x) \d x$ and
    \[m (x) = \frac{1}{\lambda} \E_\e\left[\1_{Y} \circ T_\omega^{-1} \frac{\psi\circ T_\omega^{-1}(x)}{|\det D T_\omega | \circ T_\omega^{-1}(x) }m \circ T_\omega^{-1}(x) \right].\]
\end{enumerate}
\end{lemma}

\begin{proof}
The proofs of items (1) and (2) are based on~\cite[Appendix]{BassolsCastroLamb2024} and are included here for the sake of completeness. We first prove (2), which in turn implies (1). Let $x\in \{g\neq 0\}$ and $h:Y \to \mathbb R$ be a bounded and measurable observable. We compute the conditioned Birkhoff average as follows:
\begin{align*}
 \E_x^{\log \psi} \left[ \frac{1}{n} \sum_{i=0}^{n-1} h \circ X_i^{\log \psi} \, \bigg |\, \tau >n\right] 
 &=  \frac{1}{\calP^n (\1_{Y})} \frac{1}{n} \sum_{i=0}^{n-1}\calP^{i}( h \calP^{n-i}(\1_{Y}) ) \\
 &= \frac{\lambda^n}{\calP^n (\1_{Y})} \frac{1}{n} \sum_{i=0}^{n-1}\frac{1}{\lambda^{i}}\calP^{i} \left( h \frac{1}{\lambda^{n-i}}\calP^{n-i}(\1_Y) \right)  \\
 &= \frac{\lambda^n}{\calP^n (\1_{Y})} \frac{1}{n} \left[  
 \sum_{i=0}^{n-1}\frac{1}{\lambda^{i}}\calP^{i}\left( h \left( \frac{1}{\lambda^{n-i}}\calP^{n-i}\1_{Y} -  g \cdot \mu(Y)\right) \right) \right. \\
 &\left.\hspace{4cm} +  \mu(Y)\sum_{i=0}^{n-1}\frac{1}{\lambda^{i}}\calP^{i}\left( h g\right)\right],
\end{align*}
where the last step consists of adding and subtracting $g \cdot \mu(Y)$, bringing the sum inside and rearranging the indexes.
Since
\begin{itemize}
    \item $\lim_{n \to \infty}\left\|{\lambda^{-n}}\calP^n\1_Y -  g\cdot   \mu(Y)\right\|_{\infty} = 0$, converging exponentially fast,
    \item  there exists $C>0$ such that $\left\|{\lambda^{-k}}\calP^{k} f \right\|_{\infty}\leq C \left\| f\right\|_{\infty}$ for every $k\in\mathbb N$ and every bounded and measurable function $f:M\to\R$, and
    \item the operator $\calP$ has a spectral gap and hence is mean ergodic, satisfying~\cite{Eisner2015-book, Castro2023}
    \[\lim_{n\to\infty }\frac{1}{n}\sum_{i=0}^{n-1}\frac{1}{\lambda^{i}}\calP^{i}\left( h g\right)(x) = g(x) \int h g\d\mu,\]
\end{itemize}
we obtain that 
\begin{align*}
    \lim_{n\to\infty} \E_x^{\log \psi} \left[ \frac{1}{n} \sum_{i=0}^{n-1} h \circ X_i^{\log \psi} \, \bigg |\, \tau >n\right] &= \frac{1}{g(x) \mu(Y)} \left(0 + g(x)  \mu(Y)\int_{Y} h g \d\mu \right)  \\
    &= \frac{g(x)}{g(x)} \frac{\mu(Y)}{\mu(Y)} \nu(h) = \nu(h).
\end{align*}
Thus, $\nu$ is a quasi-ergodic measure and the conditioned Birkhoff averages in~\eqref{eq:qem-cond} converge for all $x \in \{g \neq 0\}$.

Finally, we show (3). From the choice of the random perturbation introduced in Section~\ref{sec:random_perturb}, we obtain that $\calP^* \delta_x(\d y) \ll \mathrm{Leb}_{Y}(\d y)$ for every $x\in Y$, where $\delta_x(\d y)$ is the Dirac measure at $x$. Moreover, since
$$\int_{Y} \calP^* \delta_x(\d y)  \mu(\d x) =  \lambda \mu(\d y),$$
we obtain that ${\mu(\d y) \ll \mathrm{Leb}_Y(\d y).}$ In this way, there exists a density $m \in L_1(Y, \mathrm{Leb}_Y)$ such that $\mu(\d x) = m(x) \d x.$ Observe that for each $f \in \mathcal C^0(Y)$ we obtain that
\begin{align*}
   \lambda \int_{Y} f \d \mu &= \int_{Y} \calP f(x) m(x) \d x= \int_{Y}  \psi(x) \E_\e\left[ f \circ T_\omega(x) \cdot \1_{Y} \circ T_\omega(x)\right] m(x) \d x\\
    &=\E_\e\left[\int_{T_\omega(Y)} f(x)  \1_Y(x)  \psi\circ T_\omega^{-1}(x)  \frac{m\circ T_\omega^{-1}(x)}{|\det DT_\omega|\circ T_\omega^{-1}(x)} \d x \right]\\
    &=\int_Y f(x) \E_\e\left[ \frac{\psi\circ T_\omega^{-1}(x)}{{|\det DT_\omega|\circ T_\omega^{-1}(x)}}  m\circ T_\omega^{-1}(x) \cdot \1_Y\circ T_\omega^{-1}(x)  \right]\d x.
\end{align*}
The above equation implies that
\begin{equation}\label{eq:me}
    m(x) = \frac{1}{\lambda} \E_\e\left[ \frac{\psi\circ T_\omega^{-1}(x)}{{|\det DT_\omega|\circ T_\omega^{-1}(x)}}  m\circ T_\omega^{-1}(x) \cdot \1_Y\circ T_\omega^{-1}(x)  \right]. 
\end{equation}
We check that $m \in \mathcal C^0(Y)$. From~\eqref{eq:me} we obtain that for every $x\in Y$   
\begin{align*}
    \lambda m(x) &\leq \sup_{(\omega,x)\in \Omega\times Y} \left|\frac{\psi\circ T_\omega^{-1}(x)}{{|\det DT_\omega|\circ T_\omega^{-1}(x)}}\right| \E_\e[m\circ T_\omega^{-1}(x) \cdot \1_Y\circ T_\omega^{-1}(x)]\\
    &=  \sup_{(\omega,x)\in \Omega\times Y} \left|\frac{\psi\circ T_\omega^{-1}(x)}{{|\det DT_\omega|\circ T_\omega^{-1}(x)}}\right| \int_{Y} \kappa (x,y) m(y) \d y \\
    &\leq \sup_{(\omega,x)\in \Omega\times Y} \left|\frac{\psi\circ T_\omega^{-1}(x)}{{|\det DT_\omega|\circ T_\omega^{-1}(x)}}\right| \|\kappa\|_{\infty} \mu(Y) <\infty 
\end{align*}
for some bounded and measurable function $\kappa: Y\times Y \to \mathbb R$. For a detailed construction of $\kappa$ see~\cite[Proposition~3.5]{BassolsCastroLamb2024}. Hence, $m \in L_\infty(Y,\mathrm{Leb}_Y)$ and the proof is concluded by observing that the linear operator
$$h \in \mathcal F_b(Y) \mapsto \E_\e\left[ \frac{\psi\circ T_\omega^{-1}(x)}{{|\det DT_\omega|\circ T_\omega^{-1}(x)}}  h\circ T_\omega^{-1}(x) \cdot \1_Y\circ T_\omega^{-1}(x)  \right] $$
is strong Feller (see~\cite[Proposition~3.5]{BassolsCastroLamb2024}), which implies that $m$ is continuous.
\end{proof}

The following three results will allow us to apply Lemma~\ref{lem:spectral_gap_QEM}. Theorem~\ref{thm:thm51Baladi} provides a bound on the essential spectral radius of the deterministic transfer operator on $B^{t,s}$, Lemma~\ref{lem:LY-Baladi}, which partly relies on Theorem~\ref{thm:compact_embedding}, ensures that the spectral perturbation arguments of \cite{KellerLiverani1999} are applicable, and Lemma~\ref{lem:strong_feller} shows that the annealed transfer operator is strong Feller.

\begin{theorem}[{see \cite[Theorem~5.1]{Baladi2018-book} or \cite[Theorem~1.1]{Baladi-Tsujii2008}}]\label{thm:thm51Baladi}
    For $r> 1$, let $T: V\to M$ be a $\mathcal C^r$ diffeomorphism onto its image with a hyperbolic basic set $\Lambda \subset V$ and isolating neighbourhood $V$. For any $g \in \mathcal C^{r-1}(\overline{V})$ such that $\mathrm{supp}(g)\Subset V$ and all real numbers $s,t$ such that $t-(r-1) < s<0<t$, the transfer operator $\calP_g \varphi = g(\varphi \circ T)$ extends to a bounded operator on $B^{t,s}(T,V)$ and the essential spectral radius $r_{\mathrm{ess}}(\calP_g \vert_{B^{t,s}})$ satisfies
    \[\begin{split}
    r_{\mathrm{ess}}(\calP_g \vert_{B^{t,s}}) \leq Q^{t,s}(T,g) \coloneqq \exp \sup_{\mu \in \mathrm{Erg}(\Lambda, T)} \bigg\{ h_\mu (T) + \chi_\mu\left( \frac{g}{\det (DT\vert_{E^u})}\right)\\+\max\{t\chi_\mu(DT\vert_{E^s}), |s|\chi_\mu(DT^{-1}\vert_{E^u})\}\bigg\},
    \end{split}\]
    where $\mathrm{Erg}(\Lambda, T)$ denotes the set of $T$-invariant ergodic Borel probability measures on $\Lambda$, $h_\mu(T)$ denotes the metric entropy of $(\mu, T)$, and $\chi_\mu(A) \in \R \cup \{-\infty\}$ denotes the largest Lyapunov exponent of a linear cocycle $A$ over $T\vert_{\Lambda}$, with $(\log \|A\|)^+ \in L_1(\d \mu)$.
\end{theorem}

\begin{remark}\label{rmk:deterministic-spectral-data}In addition to the bound on the essential spectral radius above, the results in~\cite[Remark~1.6 and Lemmas~3.2 and~3.5]{Baladi-Tsujii2008} imply that
$$
r(\calP_g)=Q^{0,0}(T,g)
=\exp P_{\mathrm{top}}\left(T,\Lambda,\log g-\log|\det DT|_{E^u}|\right).
$$
Moreover, \cite[Remark~1.3]{Baladi-Tsujii2008} gives $Q^{t,s}(T,g)<Q^{0,0}(T,g)=r(\calP_g).$ Since $T:\Lambda\to\Lambda$ is topologically mixing, \cite[Theorem~5.1]{Gouezel-Liverani2008} gives that $r(\calP_g)$ is a simple eigenvalue and the unique eigenvalue of modulus $r(\calP_g)$. Finally, the construction in \cite[Lemmas~6.1--6.2 and Theorem~6.4]{Gouezel-Liverani2008} identifies the probability measure obtained from the leading left and right eigenvectors with the unique equilibrium state. Thus, if
$\calP_g\gamma=r(\calP_g)\gamma$ and $\calP_g^*\mu=r(\calP_g)\mu$ are normalised by $\mu(\gamma)=1$, then $
\nu^g(f)=\mu(f\gamma)$ for every continuous $f$.
 \end{remark}

\begin{lemma}\label{lem:LY-Baladi}
Consider two finite families of cone systems $\Theta=\{\Theta_i\}_{i\in I}$ and $\Theta'=\{\Theta_i'\}_{i\in I}$ satisfying $\Theta_i'<\Theta_i$ for all $i\in I$. Choose $r>1$ such that $T\in\mathcal C^r$ and $g\in\mathcal C^{r-1}$, and choose $s,t$ such that $t-(r-1)<s<0<t$. Let $0<\eta<\min\{1,t\}$ and let $\rho_0=Q^{t,s}(T,g)$ be as in Theorem~\ref{thm:thm51Baladi}. By Remark~\ref{rmk:deterministic-spectral-data}, $\rho_0<r(\calP_g),$
so we may fix $\widetilde{\rho}$ such that $\rho_0<\widetilde{\rho}<r(\calP_g).$
Then, for $\e>0$ sufficiently small, there exist $C,C_0>0$, independent of $\e$ and $m$, such that
\begin{align*}
&\|\mathcal P_\e^m\|_{B^{\Theta',t-\eta,s-\eta}}\leq CC_0^m,\\
&\|\calP_\e^m\varphi\|_{B^{\Theta,t,s}}
\leq C\widetilde{\rho}^{\,m}\|\varphi\|_{B^{\Theta,t,s}}
+CC_0^m\|\varphi\|_{B^{\Theta',t-\eta,s-\eta}},\\
&\|\calP_\e\varphi-\calP_g\varphi\|_{B^{\Theta',t-\eta,s-\eta}}
\leq C\e^\eta\|\varphi\|_{B^{\Theta,t,s}}.
\end{align*}
In particular, setting $\lambda=r(\mathcal P_g)$, the operators $\frac{1}{\lambda}\calP_g$ and $\frac{1}{\lambda}\calP_\e$ satisfy the assumptions of Theorem~\ref{thm:keller} with
$$
\mathcal B=B^{\Theta,t,s},\quad
\|\cdot\|=\|\cdot\|_{B^{\Theta,t,s}},\quad
|\cdot|=\|\cdot\|_{B^{\Theta',t-\eta,s-\eta}},\quad\text{and}\quad \tau(\e)=C\e^\eta.$$
\end{lemma}

\begin{proof}
Since $0<\eta<t$, the weak exponents remain in the admissible range $$(t-\eta)-(r-1)<s-\eta<0<t-\eta.$$
Hence multiplication by a $\mathcal C^{r-1}$ function is bounded on $B^{\Theta',t-\eta,s-\eta}$, and the composition estimates used in the proof of \cite[Theorem~5.22 and Remark~5.23]{Baladi2018-book} apply to this scale. The cone-hyperbolicity conditions are open in the $\mathcal C^1$ topology, and the maps $T_\omega$ remain in a sufficiently small $\mathcal C^r$ neighbourhood of $T$. Therefore, the constants in these estimates may be chosen uniformly for arbitrary finite compositions of the maps $T_\omega$. Averaging over the noise then gives $\|\mathcal P_\e^m\|_{B^{\Theta',t-\eta,s-\eta}}\leq CC_0^m.$

The strong Lasota--Yorke estimate in the proof of \cite[Theorem~5.22 and Remark~5.23]{Baladi2018-book} gives
$$
\|\calP_\e^m\varphi\|_{B^{\Theta,t,s}}
\leq C\widetilde{\rho}^{\,m}\|\varphi\|_{B^{\Theta,t,s}}
+CC_0^m\|\varphi\|_{B^{\Theta',t-1,s-1}}.
$$
Since $0<\eta<1$, $B^{\Theta',t-\eta,s-\eta}\hookrightarrow B^{\Theta',t-1,s-1}$
continuously, and the second estimate follows.

It remains to prove the perturbation estimate. On the common core $\mathcal C_c^u(V)$, extend $\varphi$ by zero outside $V$. The local composition estimates used in the proof of \cite[Theorem~5.22]{Baladi2018-book} give, uniformly in $\omega\in\Omega_\e$, 
$$
\|\varphi\circ T_\omega-\varphi\circ T\|_{B^{\Theta',t,s}}
\leq C\|\varphi\|_{B^{\Theta,t,s}}\quad\text{and}\quad\|\varphi\circ T_\omega-\varphi\circ T\|_{B^{\Theta',t-1,s-1}}
\leq C\e\|\varphi\|_{B^{\Theta,t,s}}.$$
For every $h$ belonging to both spaces, the definition of the anisotropic norms gives the interpolation estimate
$$
\|h\|_{B^{\Theta',t-\eta,s-\eta}}
\leq
\|h\|_{B^{\Theta',t,s}}^{1-\eta}
\|h\|_{B^{\Theta',t-1,s-1}}^\eta.
$$
Hence, $\|\varphi\circ T_\omega-\varphi\circ T\|_{B^{\Theta',t-\eta,s-\eta}}
\leq C\e^\eta\|\varphi\|_{B^{\Theta,t,s}}.$
Since multiplication by $g\in\mathcal C^{r-1}$ is bounded on the admissible space $B^{\Theta',t-\eta,s-\eta}$, averaging with respect to $\omega$ gives $$\|\calP_\e\varphi-\calP_g\varphi\|_{B^{\Theta',t-\eta,s-\eta}}
\leq C\e^\eta\|\varphi\|_{B^{\Theta,t,s}}.$$
By density, the estimate extends to $B^{\Theta,t,s}$. Finally, Theorem~\ref{thm:compact_embedding} gives the compact embedding $B^{\Theta,t,s}\hookrightarrow B^{\Theta',t-\eta,s-\eta},$
so the assumptions of Theorem~\ref{thm:keller} are satisfied by the normalised operators described above.
\end{proof}

We shall use the bounds of Lemma~\ref{lem:LY-Baladi} later in the proof of Theorem~\ref{thm:local} for vanishing weights. The final condition to be satisfied in order to apply Lemma~\ref{lem:spectral_gap_QEM} is for the operator to be strong Feller.

\begin{lemma}[{\cite[Proposition~3.5]{BassolsCastroLamb2024}}]\label{lem:strong_feller}
 For any compact set with non-empty interior $Y \subset M$ and any non-negative continuous function $\varphi:Y\to \mathbb R$, the operator
 \begin{align*}
     \calP_{\e,Y} : \mathcal F_b(Y) &\to \mathcal F_b(Y)\\
     h&\mapsto \varphi(x)\E_x\left[h\circ X_1^\e \cdot \1_{Y}\circ X_1^{\e}\right],
 \end{align*}
where $\mathcal F_b(Y) \coloneqq \{h:Y\to \R;\ h \text{ is bounded and measurable}\}$, is strong Feller, i.e.~for every $h\in \mathcal F_b(Y)$ we have that $\calP_{\e, Y} h $ is a continuous function. The same strong Feller property holds for the corresponding annealed transfer operator
$$
\mathcal L_{\e,Y}f(y)
\coloneqq
\E_\e\left[
\1_Y\circ T_\omega^{-1}(y)
\frac{\varphi\circ T_\omega^{-1}(y)}
{|\det DT_\omega|\circ T_\omega^{-1}(y)}
f\circ T_\omega^{-1}(y)
\right],
$$
which is the predual of $\calP_{\e,Y}$ with respect to Lebesgue measure.
\end{lemma}

We are ready to prove the local existence and uniqueness of quasi-ergodic measures as well as conditioned stochastic stability, i.e.~Theorem~\ref{thm:local}, for H\"older weights vanishing on the boundary of the isolating neighbourhood $V$.

\begin{proof}[Proof of Theorem~\ref{thm:qem-local} and Theorem~\ref{thm:local} for vanishing H\"older weights]
Let us consider three finite families of cone systems $\Theta = \{\Theta_i\}_{i \in I}$, $\Theta' = \{\Theta'_i\}_{i \in I}$ and $\Theta'' = \{\Theta''_i\}_{i \in I}$ satisfying $\Theta''_i<\Theta'_i<\Theta_i$ for all $i\in I$. Choose $s',t'$ sufficiently close to $s,t$ so that $B^{\Theta,t,s}$ compactly embeds in $B^{\Theta',t',s'}$, which in turn continuously embeds in the weak space appearing in Lemma~\ref{lem:LY-Baladi}. We denote the norm of this weak space by $|\cdot|$.

Let $\mu \in (B^{\Theta,t,s})^*$ and $\gamma \in B^{\Theta,t,s}$ be the leading left and right eigenvectors from Remark~\ref{rmk:deterministic-spectral-data}, normalised so that $\mu(\gamma)=1$. Thus, $\calP_g^*\mu=\lambda\mu$, $\calP_g\gamma=\lambda\gamma$, $\lambda=r(\calP_g)$, and $\mu(\phi\gamma)=\nu^g(\phi)$ for every continuous $\phi$.

For $\delta > 0$ small enough such that $B_\delta(\lambda)\cap\sigma(\calP_g)=\{\lambda\}$, define the spectral projection
\begin{equation}\label{eq:projection}
    \Pi_\e^{(\lambda,\delta)}\coloneqq\frac{1}{2\pi i}\int_{\partial B_\delta(\lambda)}(z-\calP_\e)^{-1}\d z.
\end{equation}
From item~\ref{it:keller1} of Theorem~\ref{thm:keller},
$$
|\Pi_\e^{(\lambda,\delta)}\gamma-\gamma|
\leq K_1\tau(\e)^\eta\|\gamma\|_{B^{\Theta,t,s}}\to0.
$$
For $\e>0$ small enough, set $g_\e\coloneqq\Pi_\e^{(\lambda,\delta)}\gamma$. By item~\ref{it:keller3} of Theorem~\ref{thm:keller}, $\lambda_\e$ is simple, so we may choose $\mu_\e\in(B^{\Theta,t,s})^*$ such that $\Pi_\e^{(\lambda,\delta)}f=\mu_\e(f)g_\e$ for every $f\in B^{\Theta,t,s}$. In particular, $\calP_\e g_\e=\lambda_\e g_\e$, $\calP_\e^*\mu_\e=\lambda_\e\mu_\e$ and $\mu_\e(g_\e)=1$, where $\lambda_\e=r(\calP_\e)$.

From item~\ref{it:keller2} of Theorem~\ref{thm:keller}, there exists $K>0$ such that
\begin{equation}\label{eq:proj-ineq}
\|\Pi_\e^{(\lambda,\delta)}f\|_{B^{\Theta,t,s}}
\leq K|\Pi_\e^{(\lambda,\delta)}f|
\end{equation}
for every $f\in B^{\Theta,t,s}$. Applying this with $f=\gamma$, we obtain
$$
\|g_\e\|_{B^{\Theta,t,s}}
\leq K|g_\e|,
$$
and hence $\{g_\e\}_{\e>0}$ is uniformly bounded in $B^{\Theta,t,s}$.

We show that for every $\phi \in \mathcal C^{\infty}(\overline{V})$, the limit
$
\lim_{\e\to0}\mu_\e(\phi\cdot g_\e)=\mu(\phi\cdot\gamma)
$
holds.

Recall that $B^{\Theta,t,s}$ compactly embeds in $B^{\Theta',t',s'}$ by Theorem~\ref{thm:compact_embedding}. Thus, every sequence $\e_n\to0$ admits a subsequence such that $g_{\e_n}\to g_0$ in $B^{\Theta',t',s'}$. Since $B^{\Theta',t',s'}$ continuously embeds in the weak space and $|g_\e-\gamma|\to0$, necessarily $g_0=\gamma$. Since this holds for every subsequence, we conclude that $g_\e\to\gamma$ in $B^{\Theta',t',s'}$ as $\e\to0$.

Moreover, item~\ref{it:keller3} of Theorem~\ref{thm:keller} gives
\begin{equation}\label{eq:growth-rate}
\lim_{\e \to 0} r(\calP_\e) = \lim_{\e \to 0}\lambda_\e = \lambda = r(\calP_g).
\end{equation}

Applying Lemma~\ref{lem:LY-Baladi} with the cone systems $\Theta'$ and $\Theta''$ and exponents $t',s'$, choosing the corresponding weak space within the admissible regularity range, Theorem~\ref{thm:keller} also applies with strong space $B^{\Theta',t',s'}$. The corresponding Riesz projections agree with $\Pi_\e^{(\lambda,\delta)}$ on $B^{\Theta,t,s}$. Moreover, items~\ref{it:keller1} and~\ref{it:keller2} imply that, for some $C>0$ independent of $\e$,
$\|\Pi_\e^{(\lambda,\delta)}f\|_{B^{\Theta',t',s'}}
\leq C\|f\|_{B^{\Theta',t',s'}}$ for every $f\in B^{\Theta',t',s'}$ and $\e>0$ small enough. Since $\lambda_\e$ is simple and $g_\e\to\gamma$ in $B^{\Theta',t',s'}$, the identity $\Pi_\e^{(\lambda,\delta)}f=\mu_\e(f)g_\e$ defines an extension of $\mu_\e$ to $(B^{\Theta',t',s'})^*$ and $\sup_{\e>0}\|\mu_\e\|_{(B^{\Theta',t',s'})^*}<\infty.$
The same argument for $\e=0$ shows that $\mu$ extends to $(B^{\Theta',t',s'})^*$.

On the other hand, item~\ref{it:keller1} for the original strong and weak spaces gives, for every $f\in B^{\Theta,t,s}$,
$$
|\mu_\e(f)g_\e-\mu(f)\gamma|
\leq K_1\tau(\e)^\eta\|f\|_{B^{\Theta,t,s}}.
$$
Therefore,
$$
|\mu_\e(f)-\mu(f)|\,|g_\e|
\leq K_1\tau(\e)^\eta\|f\|_{B^{\Theta,t,s}}
+|\mu(f)|\,|\gamma-g_\e|.
$$
Since $|g_\e-\gamma|\to0$ and $|\gamma|>0$, it follows that $\mu_\e(f)\to\mu(f)$ for every $f\in B^{\Theta,t,s}$. Since $B^{\Theta,t,s}$ is dense in $B^{\Theta',t',s'}$ and the functionals $\mu_\e$ are uniformly bounded on $B^{\Theta',t',s'}$, we conclude that
$$
\mu_\e\xlongrightarrow[\mathrm{w}^*]{}\mu
\quad\text{in }(B^{\Theta',t',s'})^*.
$$

Finally, let $\phi \in \mathcal C^\infty(\overline{V})$. By the multiplication bounds on the anisotropic spaces, $\phi g_\e\to\phi\gamma$ in $B^{\Theta',t',s'}$. Hence,
$$|\mu_\e(\phi g_\e)-\mu(\phi\gamma)|
\leq \|\mu_\e\|_{(B^{\Theta',t',s'})^*}\|\phi(g_\e-\gamma)\|_{B^{\Theta',t',s'}}
+|(\mu_\e-\mu)(\phi\gamma)|,$$
and therefore $\mu_\e(\phi g_\e)\longrightarrow\mu(\phi\gamma)=\nu^g(\phi).$

To conclude the proof of the theorem, we show that $\nu_\e(\d x) = g_\e(x) \mu_\e(\d x)$ is a quasi-ergodic measure of $X_\e^{\log g}$ on $\overline{V}$. The operator $\calP_{\e}$ is strong Feller so $\calP^2_{\e}$ is compact (recall Lemma~\ref{lem:strong_feller}) and therefore $\calP_\e: \mathcal C^0(\overline{V})\to \mathcal C^0(\overline{V})$ is quasi-compact~\cite[Equation~(8), Theorem~4]{Yosida2}. We argue that $\calP_\e: \mathcal C^0(\overline{V})\to \mathcal C^0(\overline{V})$ has a spectral gap. To see this, observe that since $\calP_\e(\mathcal C^0(\overline{V})) \subset \mathcal C^0_0 (\overline{V})$, we obtain that $\left.\calP_\e\right|_{\mathcal C^0_0 (\overline{V})}$ is a quasi-compact operator. 
Choose $u$ such that $t<u<r-1.$
Since the weight $g$ is compactly supported in $\overline{V}$, the space $\mathcal C_0^u(\overline{V})$ is invariant under $\calP_\e$. Moreover, $\mathcal C_0^u(\overline{V})$ is dense in both $\mathcal C_0^0(\overline V)$ and $B^{t,s}(T,V)$. Therefore Proposition~\ref{prop:SpecBaladi} applies with
$$
\mathcal B_0=\mathcal C_0^u(\overline{V}),\qquad
\mathcal B_1=\mathcal C_0^0(\overline V),\qquad
\mathcal B_2=B^{t,s}(T,V).
$$
It follows that the spectra of $\calP_\e$ on $\mathcal C_0^0(\overline V)$ and $B^{t,s}(T,V)$ coincide outside the corresponding essential spectral radii. Since $\calP_\e$ has a spectral gap on $B^{t,s}(T,V)$, the same holds on $\mathcal C_0^0(\overline V)$.
Moreover, observe that $g_\e\in\mathcal C_0^0(\overline V)$ and $\mu_\e \in \mathcal M (\overline{V})$. A posteriori, we conclude that $\calP_\e:\mathcal C^0(\overline{V})\to \mathcal C^0(\overline{V})$ also has a spectral gap as desired. Finally, invoking Lemma~\ref{lem:spectral_gap_QEM}, it follows that $\nu_\e (\d x)= g_\e(x) \mu_\e (\d x)$ is a quasi-ergodic probability measure of $X_\e^{\log g}$ on $\overline{V}$. Since $\nu_\e$ and $\nu^g$ are probability measures, the convergence established above for smooth observables extends to every $\phi\in\mathcal C^0(\overline V)$ by uniform approximation. Indeed, if $\phi_n\in\mathcal C^\infty(\overline V)$ satisfies $\|\phi_n-\phi\|_\infty\to0$, then
$$
|\nu_\e(\phi)-\nu^g(\phi)|\leq 2\|\phi-\phi_n\|_\infty+|\nu_\e(\phi_n)-\nu^g(\phi_n)|.
$$
Letting first $\e\to0$ and then $n\to\infty$, we obtain $\nu_\e(\phi)\to\nu^g(\phi)$ for every $\phi\in\mathcal C^0(\overline V)$, and hence $\nu_\e\xlongrightarrow[\mathrm{w}^*]{}\nu^g$ in $\mathcal M(\overline V)$.
\end{proof}

\begin{remark}\label{rmk:spectral_gap}
    We emphasise that in the last argument in the proof of Theorem~\ref{thm:qem-local} and Theorem~\ref{thm:local} for vanishing H\"older weights we have shown that $\calP_\e$ has a spectral gap in $\mathcal C^0 (\overline{V})$. This is crucial to apply Lemma~\ref{lem:spectral_gap_QEM} and is used in the following section.
\end{remark}

\section{Local quasi-ergodic measures for non-vanishing weights}\label{sec:QEMs-non-vanishing}

Thus far, we proved the local existence and uniqueness of quasi-ergodic measures for $\e > 0$ small and, importantly, for non-negative weights vanishing on the boundary of the isolating neighbourhood $V \supset \Lambda$. In this section, we extend these results to the case of non-negative weights which do not (necessarily) vanish on $\partial V \coloneqq \overline{V} \setminus V$. 

Given a $\mathcal C^{r-1}$ H\"older weight $e^\phi:\overline{V}\to(0,\infty)$, consider a suitable neighbourhood $\U\supset\Lambda$, which we shall specify later in this section. Since $\overline{\U}\subset V$, choose $\chi\in\mathcal C_c^\infty(V)$ such that $0\leq\chi\leq1$ and $\chi\equiv1$ on a neighbourhood of $\overline{\U}$, and define
\begin{equation}\label{eq:phi}
\varphi\coloneqq\chi e^\phi.
\end{equation}
Then the zero extension of $\varphi$ belongs to $\mathcal C^{r-1}(M)$, is compactly supported in $V$, and satisfies $\varphi=e^\phi$ on a neighbourhood of $\overline{\U}$. Thus, $\varphi$ is a weight for which Section~\ref{sec:QEMs-vanishing} applies on $V$. We apply the results of that section with $g$ replaced by $\varphi$. The notation $X_\e^{\log\varphi}$ is understood according to Remark~\ref{rmk:positive}, including at points where $\varphi=0$.

The strategy we present here consists in constructing a suitable neighbourhood $\U$ of $\Lambda$ satisfying $\Lambda \subset \U \subset \overline{\U} \subset V$ and such that the quasi-ergodic measures of:
\begin{itemize}
    \item $X_\e^{\log \varphi }$ absorbed outside of $\overline{V}$ and with weight $\log \varphi$,
    \item  $X_\e^{\phi}$ absorbed outside of $\overline{V}$ and with weight $\phi$, and
    \item  $\widetilde{X}_\e^{\phi}$ absorbed outside of $\overline{\U}$ and with weight $\phi$
\end{itemize}
all coincide on $\overline{\U}$. So far, we have shown the existence, uniqueness and stochastic stability of the quasi-ergodic measure for $X_\e^{\log \varphi}$, which we proved is built from the left and right dominant eigenfunctions of $\mathcal P_\e$. In Section~\ref{sec:correspondence-qems}, we show that the quasi-ergodic measures for $X_\e^{\phi}$ and $\widetilde{X}_\e^{\phi}$ are built analogously with eigenfunctions that can be obtained from those of $\mathcal P_\e$.

\subsection{Constructing \texorpdfstring{$\U$}{U}}\label{sec:construct-U}
Recall that $V$ is an isolating neighbourhood of $\Lambda$ and write
$$W_V^s(\Lambda)\coloneqq\bigcap_{n\geq0}T^{-n}(\overline V),\qquad W_V^u(\Lambda)\coloneqq\bigcap_{n\geq0}T^n(\overline V).$$
These sets are compact and satisfy $W_V^s(\Lambda)\cap W_V^u(\Lambda)=\Lambda$. We construct an arbitrarily small neighbourhood $\U$ of $\Lambda$ such that points leaving $\U$ cannot land sufficiently close to $W_V^s(\Lambda)$, uniformly under sufficiently small random perturbations.

\begin{lemma}\label{lem:set_U}
For every $\rho>0$ such that $\overline{B_\rho(\Lambda)}\subset V$, there exists $N\geq1$ such that the open set
\begin{equation}\label{eq:local-U-construction}
\U\coloneqq\left\{x\in M:\max_{0\leq j<N}2^{j/N}\dist(T^j(x),\Lambda)<\rho\right\}
\end{equation}
satisfies $\Lambda\subset\U\subset\overline\U\subset V$ and
\begin{equation}\label{eq:local-U-strict-inclusion}
T(\overline\U)\cap W_V^s(\Lambda)\subset\U.
\end{equation}
\end{lemma}
\begin{proof}
Since $T(W_V^s(\Lambda))\subset W_V^s(\Lambda)$, the compact sets $T^n(W_V^s(\Lambda))$, $n\geq0$, form a decreasing sequence. Moreover,
\begin{equation}\label{eq:local-U-stable-intersection}
\bigcap_{n\geq0}T^n(W_V^s(\Lambda))=\bigcap_{j\in\Z}T^j(\overline V)=\Lambda.
\end{equation}
Indeed, a point lies in the first intersection precisely when all its forward and backward iterates lie in $\overline V$. By compactness and \eqref{eq:local-U-stable-intersection}, there exists $N\geq1$ such that
\begin{equation}\label{eq:local-U-stable-tail}
T^N(W_V^s(\Lambda))\subset B_{\rho/2}(\Lambda).
\end{equation}
Define $\U$ by \eqref{eq:local-U-construction}. This set is open and contains $\Lambda$, since $T(\Lambda)=\Lambda$. The condition corresponding to $j=0$ gives $\U\subset B_\rho(\Lambda)$, and therefore $\overline\U\subset\overline{B_\rho(\Lambda)}\subset V$.

Let $x\in\overline\U$ and suppose that $T(x)\in W_V^s(\Lambda)$. Since $x\in V$ and all the forward iterates of $T(x)$ lie in $\overline V$, we also have $x\in W_V^s(\Lambda)$. Using \eqref{eq:local-U-construction} and \eqref{eq:local-U-stable-tail}, we obtain
\begin{align}\label{eq:local-U-bound}
\max_{0\leq j<N}2^{j/N}\dist(T^j(T(x)),\Lambda)
&=2^{-1/N}\max_{1\leq j\leq N}2^{j/N}\dist(T^j(x),\Lambda)\nonumber\\
&\leq2^{-1/N}\max\left\{\rho,2\dist(T^N(x),\Lambda)\right\}\leq2^{-1/N}\rho<\rho.
\end{align}
Hence $T(x)\in\U$ by \eqref{eq:local-U-construction}, proving \eqref{eq:local-U-strict-inclusion}.
\end{proof}

The use of $\overline\U$ in \eqref{eq:local-U-strict-inclusion} gives the separation needed to obtain the following stochastic version.

\begin{lemma}\label{lem:set_U_stoch}
For $\U$ constructed in Lemma~\ref{lem:set_U}, there exist $\delta>0$ and $\e_0>0$ such that, for every $0<\e<\e_0$, $x\in\overline\U$ and $\omega\in\Omega_\e$,
\begin{equation}\label{eq:local-U-random-nonreturn}
T_\omega(x)\notin B_\delta(W_V^s(\Lambda))\setminus\U.
\end{equation}
\end{lemma}
\begin{proof}
Suppose otherwise. Then there exist $0<\e_n<1/n$, $x_n\in\overline\U$ and $\omega^{(n)}\in\Omega_{\e_n}$ such that
$$
T_{\omega^{(n)}}(x_n)\in B_{1/n}(W_V^s(\Lambda))\setminus\U.
$$
Passing to a subsequence, compactness gives $x_n\to x\in\overline\U$. From the definition of the perturbations $T_\omega$ in Section~\ref{sec:random_perturb}, $
T_{\omega^{(n)}}(x_n)\longrightarrow T(x).$
Since $W_V^s(\Lambda)$ and $M\setminus\U$ are closed, it follows that $T(x)\in W_V^s(\Lambda)\setminus\U,$
contradicting \eqref{eq:local-U-strict-inclusion}.
\end{proof}

Finally, since $W_V^s(\Lambda)$ and $W_V^u(\Lambda)$ are compact and their intersection is $\Lambda\subset\U$, we may decrease $\delta$ so that
\begin{equation}\label{eq:local-U-stable-unstable-separation}
B_\delta(W_V^s(\Lambda))\cap B_\delta(W_V^u(\Lambda))\subset\U.
\end{equation}
Indeed, otherwise there would exist $x_n\in M\setminus\U$ whose distances to both sets tend to zero. A convergent subsequence would then have its limit in $W_V^s(\Lambda)\cap W_V^u(\Lambda)=\Lambda$, contradicting $\Lambda\subset\U$.

\begin{notation}\label{not:delta}
We fix $\U$ as constructed in Lemma~\ref{lem:set_U} and choose $\delta>0$ and $\e_0>0$ such that \eqref{eq:local-U-random-nonreturn} and \eqref{eq:local-U-stable-unstable-separation} hold.
\end{notation}

\subsection{Correspondence of quasi-ergodic measures}\label{sec:correspondence-qems}
Consider the three operators
\begin{align*}
    \overline{\mathcal P}_\e f \coloneqq \varphi(x)\E_x\left[f \circ X_1^\e \cdot \1_{\overline{V}}\circ X^\e_1 \right]\\
    \mathcal P_\e f \coloneqq e^{\phi(x)}\E_x\left[f \circ X_1^\e \cdot \1_{\overline{V}}\circ X^\e_1 \right] \\
    \mathcal P_{\e, \U} f \coloneqq e^{\phi(x)}\E_x\left[f \circ \widetilde{X}_1^\e \cdot \1_{\overline{\U}}\circ \widetilde{X}^\e_1 \right]
\end{align*}
where $\widetilde{X}^\e$ is the process $X^\e$ absorbed outside of $\overline{\U}$, i.e.~if $\widetilde{X}^{\e}_n \not\in \overline{\U}$ for some $n \in \N$ then $\widetilde{X}^{\e}_n \in \partial$ and $\widetilde{X}^{\e}_{n+1} \in \partial$, and $\varphi$ is defined in \eqref{eq:phi}. Moreover, recall from Remark~\ref{rmk:abuse}, we say that $T_\omega(\partial) = \partial$ for all $\w \in \Omega_\e$; in this way, if $T_\omega^i(x) \not\in \overline{V}$ for some $i\in\mathbb N$, we have that $T^{i+j}_\omega (x) \in \partial$ for every $j\in\mathbb N.$ Finally, note that in this section, the operator $\mathcal P_\e$ of Section~\ref{sec:QEMs-vanishing} is denoted by $\overline{\mathcal P}_\e$. 

We are able to control the support of the peripheral eigenfunctions of $\calP_\e$ as follows.

\begin{lemma}\label{lem:support_g}
    For every $\e <\e_0$, if a continuous function $g_\varepsilon:\overline{V}\to\mathbb R$ satisfies $\overline{\mathcal P}_\e g_\e = \lambda_\e g_\e$ (resp. $\mathcal P_\e g_\e = \lambda_\e g_\e)$ and $\lambda_\e = r(\overline{\calP}_\e)$ (resp. $r({\calP}_\e)$) , then $\{g_\e > 0\} \subseteq \U \cup B_\delta(W^s_V(\Lambda))$.
\end{lemma}
\begin{proof}
    We show the result only for the operator $\overline{\mathcal P}_\e$ and note that the same proof applies when considering the operator $\mathcal P_\e$. Observe that there exists $N\in\mathbb N$ such that $T^{N}(x) \not\in \overline{V}$ for every $x \not\in \U \cup B_\delta(W^s_V(\Lambda))$. Since $\overline{V}$ is compact, there exists $\e>0$ small enough such that $T_\w^{N}(x) \not\in \overline{V}$, for every $x\in \overline{V}\setminus \left(\mathcal U \cup B_{\delta}(W^s_V(\Lambda)) \right)$ and $ \w \in \Omega_\e$. Hence,
    \[ g_\e(x) = \frac{1}{\lambda_\e^N}  \E_\e\left[\prod_{i = 0}^{N-1} \varphi \circ T_\w^i(x) \cdot g_\e\circ T_\w^N(x) \cdot \1_{\overline{V}}\circ T^N_\omega(x)\right] = 0\]
    and the claim follows.
\end{proof}

This allows us to obtain eigenfunctions for $\calP_{\e, \U}$ from those of $\overline{\calP}_\e$ and ${\calP}_\e$ by simply restricting them to $\overline{\U}$.

\begin{lemma} \label{lem:restrict_eval}
    For every $\e <\e_0$ and $g_\e\in \ker(\overline{\calP}_\e - \lambda_\e)$, with $\lambda_\e = r(\overline{\calP}_\e)$ we have that $\left.g_\e\right|_{\overline{\U}} \neq 0$ and $\calP_{\e, \U}(\left.g_\e\right|_{\overline{\U}}) = \lambda_\e \left.g_\e\right|_{\overline{\U}}$. The same statement holds for $\mathcal{P}_\e$ instead of $\overline{\calP}_\e$.
\end{lemma}
\begin{proof}
   We show the result only for the operator $\overline{\mathcal P}_\e$ and note that the same proof applies when considering the operator $\mathcal P_\e$. Given $x \in \overline{\U}$, from Lemma~\ref{lem:set_U_stoch} we know that $T_\w(x) \not \in B_\delta(W^s_V(\Lambda)) \setminus \U$ for every $\w \in \Omega_\e$, which implies that $\1_{\U}\circ T_\w(x) = \1_{{\U} \cup B_\delta(W^s_V(\Lambda))}\circ T_\w(x)$ on $\{g_\e\circ T_\w(x)\neq0\}$. Thus, for $x \in \overline{\U}$, and using Lemmas~\ref{lem:set_U_stoch} and \ref{lem:support_g} we obtain
    \begin{align*}
    \calP_{\e, \U} (\left.g_\e\right|_{\overline{\U}})(x)&= e^{\phi(x)}\E_\e\left[g_\e\circ T_\w(x) \cdot \1_{\overline{\U}} \circ T_\w (x) \right] = e^{\phi(x)}\E_\e\left[g_\e\circ T_\w(x) \cdot \1_{\U} \circ T_\w (x) \right] \\& = \varphi(x)\E_\e\left[g_\e\circ T_\w(x) \cdot \1_{\U \cup B_\delta(W^s_V(\Lambda))} \circ T_\w (x) \right] \\
    &= \varphi(x) \E_\e\left[g_\e \circ T_\w(x) \cdot \1_{\overline{V}} \circ T_\w(x)\right] = \overline{\calP}_\e g_\e(x) = \lambda_\e \left.g_\e\right|_{\overline{\U}}(x).
    \end{align*}
    
    We show that $\left.g_\e\right|_{\overline{\U}}\neq 0.$  Fix $0<\e<\e_0,$ and assume for a contradiction that $\left.g_\e\right|_{\overline{\U}} =0$. From Lemma~\ref{lem:support_g}, there exists $\delta>0$ such that $\{g_\e >0\}\subset \U\cup B_\delta(W^s_V(\Lambda)).$ Observe that for every $x\in B_\delta(W^s_V(\Lambda))$ there exists $N>0$ such that either $T^i_\omega(x)\in \overline{\mathcal U}$ for some $i\in\{1,\ldots, N-1\}$ or $T^{N-1}_\omega (x) =\partial$. From Lemma~\ref{lem:set_U_stoch} we have that  $T^N_\omega(x) \not\in  B_\delta(W^s_V(\Lambda))\setminus {\U}$  for every $x\in B_\delta(W^s_V(\Lambda))\setminus \overline{\U}$ and $\omega\in \Omega_\e$. It follows that
    $$g_\e(x) =\frac{1}{\lambda_\e^{N} } \E_\e\left[e^{S_N\phi(x)} g_\e\circ T_\omega^N(x)  \right] =\frac{1}{\lambda_\e^{N} } \E_\e\left[e^{S_N\phi(x)} \left.g_\e\right|_{\overline{\mathcal U}}\circ T_\omega^N(x)  \right]= 0,$$ 
    since $g_\e(\partial) = 0$ and $T_\omega(\partial) = \partial$ for every $\omega\in \Omega_\e.$
\end{proof}

We may also go in the other direction and obtain eigenfunctions for $\overline{\calP}_\e$ from those of $\calP_{\e,\U}$.

\begin{lemma}\label{lem:lift_eval} For each eigenfunction $g$ satisfying $\mathcal P_{\e, \U} g= \lambda g$, $\lambda \neq 0$, we can induce an eigenfunction of $\overline{\calP}_\e$ and of $\calP_\e$.
\end{lemma}
\begin{proof}
    We prove the statement for $\overline{\calP}_\e$ and note that the same proof for $\calP_\e$ follows changing $\varphi$ for $e^{\phi}$ below. For each $x \in \overline{V}, \w \in \Omega_\e$ let us define $\sigma(\w, x) \coloneqq \min\{n \geq 0:\, T^n_\w(x) \in \overline{\U} \cup \partial\}$, which is a uniformly bounded stopping time. Consider the function
    \[\overline{g}(x) \coloneqq \E_\e\left[ \tilde{g}(\w,x) \right] \coloneqq \E_\e\left[ \frac{1}{\lambda^{\sigma(\w, x)}} \prod_{i = 0}^{\sigma(\w, x) -1} \varphi\circ T^i_\w(x) \cdot g\circ T^{\sigma(\w, x)}_\w(x) \right].\]
    Notice that if $x \in \overline{\U}$, then $\sigma(\w, x) =0$ and $\overline{g}(x) = g(x)$. We check that $\overline{\mathcal P}_\e \overline{g} = \lambda \overline{g}$. If $x \in \overline{V} \setminus \overline{\U}$, then $\sigma(\w, x) \geq 1$ and $\sigma \circ\Theta(\w, x)= \sigma(\w, x) -1$. Using the Markov property, we obtain
    \begin{align*}
        \overline{\mathcal P}_\e \overline{g}(x) &= \varphi(x)\E_\e^{\w}\left[\1_{\overline{V}}\circ T_\w(x) \cdot 
        \E^{\nu}_{\e}\left[ \tilde{g}(\nu, T_\w(x)) \right]\right]\\
        &=\E_\e\left[\1_{\overline{V}}\circ T_\w(x) \frac{\lambda}{\lambda^{\sigma (\w, x)}} \prod_{i = 0}^{\sigma(\w, x) -1} \varphi\circ T^i_\w(x) \cdot g\circ T^{\sigma(\w, x)}_\w(x)\right] = \lambda \overline{g}(x),
    \end{align*}
    where we have used the superscript on the expectation symbol to denote the noise variable with respect to which the average is taken.
    If $x \in \overline{\U}$, given $\omega\in\Omega_\e$ there are two possible cases: either $T_\w(x) \in \overline{\U}$ or $T_\w(x) \in \overline{V} \setminus \overline{\U}$. Therefore, 
    \begin{align*}
        \overline{\calP}_\e \overline{g}(x)&=\varphi(x)\E_\e\left[\1_{\overline{\U}} \circ T_\w(x) \cdot \overline{g} \circ T_\w(x) \right] + \varphi(x)\E_\e\left[\1_{\overline{V}\setminus \overline{\U}} \circ T_\w(x)\cdot \overline{g} \circ T_\w(x) \right]\\
        &= \mathcal P_{\e,\U}g(x) + \varphi(x)\E_\e\left[\1_{\overline{V}\setminus \overline{\U}} \circ T_\w(x)\cdot \overline{g} \circ T_\w(x) \right]\\
        &= \lambda \bar{g}(x) + \varphi(x)\E_\e\left[\1_{\overline{V}\setminus \overline{\U}} \circ T_\w(x)\cdot \overline{g} \circ T_\w(x) \right].
    \end{align*}
   To prove the claim it is sufficient to show that $\E_x\left[\1_{\overline{V}\setminus \overline{\U}} \circ T_\w(x)\cdot \overline{g} \circ T_\w(x) \right] = 0$ for every $x\in \overline{\U}$, which is equivalent to showing that $\overline{g} \circ T_\omega(x) =0$ if $T_\omega(x) \not \in \overline{\U}$ for each $x\in \overline{\U}$. From Lemma~\ref{lem:set_U_stoch}, since $T_\w(x) \not\in \U$ we have that $T_\w (x) \not \in B_\delta(W^s_V(\Lambda))$ so there exists $N> 0$ sufficiently large such that for any $\nu \in \Omega_\e$, with $\e>0$ small enough, it holds that $T^N_\nu \circ T_\w(x) \in \partial$. In particular, $g\circ T^{\sigma(\nu,T_\w(x))}_\nu(T_\w(x)) = 0$ for all $\nu \in \Omega_\e$ so $\overline{g}\circ T_\w(x) = 0$.
\end{proof}

With Lemma~\ref{lem:spectral_gap_QEM} in mind, we show that $\calP_\e$ and $\calP_{\e, \U}$ have a spectral gap in $\mathcal C^0$ and that their spectral radii coincide with that of $\overline{\calP}_\e$.

\begin{lemma}\label{lem:sg_PeU}
    The operator $\mathcal P_{\e, \U}$ has a spectral gap in $\mathcal C^0(\overline{\U})$ and $r(\calP_{\e,\U}) = r(\overline{\calP}_\e)$.
\end{lemma}
\begin{proof}
    Observe that for all $x \in \overline{\U}$ and non-negative $f \in \mathcal C^{0}(\overline{V}),$ we have $\overline{\calP}_\e(f \1_{\overline\U})(x) \geq \calP_{\e, \U}(f \1_{\overline{\U}})(x)$ so that $r(\calP_{\e, \U}) \leq r(\overline{\calP}_\e)$. On the other hand, Lemma~\ref{lem:restrict_eval} provides $r(\calP_{\e, \U}) \geq r(\overline{\calP}_\e)$, so $r(\calP_{\e, \U}) = r(\overline{\calP}_\e)$.

    The operator $\calP_{\e, \U}$ is strong Feller and thus $\calP^2_{\e,\U}$ is compact (see~\cite[Proposition~3.5]{BassolsCastroLamb2024}). Moreover, $\frac{1}{\lambda_\e}\overline{\calP}_\e$ is power-bounded since it has a spectral gap in $\mathcal C^0(\overline{V})$ (see Remark~\ref{rmk:spectral_gap}), and so is $\frac{1}{\lambda_\e}\calP_{\e, \U}$. Finally, since there is a one-to-one correspondence between the eigenvalues of $\mathcal P_{\e,\U}$ and $\overline{\mathcal P}_\e$ (recall Lemmas \ref{lem:restrict_eval} and \ref{lem:lift_eval}), this provides the presence of a spectral gap in $\mathcal C^0(\overline{\U})$ for $\calP_{\e, \U}$. 
\end{proof}

\begin{lemma}\label{lem:PeSG}
    The operator $\mathcal P_\e$ has a spectral gap in $\mathcal C^0(\overline{V})$ and $r(\calP_{\e}) = r(\calP_{\e, \U}) = r(\overline{\calP}_\e).$ 
\end{lemma}
\begin{proof}
     $\mathcal P_\e$ is a strong Feller operator and so quasi-compact. Using a similar construction to that in Lemma~\ref{lem:lift_eval}, we may induce eigenvectors for $\calP_\e$ from those of $\calP_{\e,\U}$. Moreover, observe that $|\calP_{\e, \U} 1| \leq |\calP_{\e}1|$ so that $r(\calP_{\e, \U}) \leq r(\calP_{\e})$, and following Lemmas~\ref{lem:support_g} and \ref{lem:restrict_eval} we have that $r(\calP_{\e}) \leq r(\calP_{\e, \U})$. It follows that $r(\calP_{\e}) = r(\calP_{\e, \U})$.
     
     We show that $\calP_\e$ does not have any Jordan blocks on the peripheral spectrum. Assume that $f_\e$ satisfies $\calP_\e f_\e = \lambda_\e f_\e + g_\e$, with $g_\e$ an eigenfunction of $\calP_\e$ associated with the eigenvalue $\lambda_\e = r(\calP_\e)$. Following the same steps as in the proof of Lemma~\ref{lem:support_g} we have that $\calP_{\e,\U} (f_\e\1_{\overline\U} ) = \lambda_\e f_\e\1_{\overline\U}+ g_\e\1_{\overline\U}$. However, $g_\e\1_{\overline\U} \neq 0$ and $\calP_{\e, \U}$ does not have Jordan blocks as $\frac{1}{\lambda_\e}\mathcal P_{\e, \U}$ is power-bounded, thus reaching a contradiction. 
\end{proof}

Bringing all these properties together, we conclude with the main result of this section. Corollary~\ref{coro:all-qems-same} provides a clear link between the quasi-ergodic measures of the three processes considered and extends the existence and uniqueness of quasi-ergodic measures for non-vanishing weights.

\begin{corollary}\label{coro:all-qems-same}
    Each of the three processes:
    \begin{itemize}
    \item $X_\e^{\log \varphi },$ absorbed outside of $\overline{V}$ and with weight $\varphi$,
    \item  $X_\e^{\phi},$ absorbed outside of $\overline{V}$ and with weight $e^\phi$, and
    \item  $\widetilde{X}_\e^{\phi},$ absorbed outside of $\overline{\U}$ and with weight $e^\phi$
\end{itemize}
admits the same quasi-ergodic measure $\nu_\e$, which is supported on $\overline{\U}$.
\end{corollary}
\begin{proof}
    The transfer operators $\overline{\calP}_\e, \calP_{\e},$ and $\calP_{\e, \U}$ are all strong Feller (see Lemma~\ref{lem:strong_feller}), have a spectral gap in $\mathcal C^0$ (see Remark~\ref{rmk:spectral_gap} and Lemmas~\ref{lem:PeSG} and \ref{lem:sg_PeU}) and their left and right eigenfunctions coincide in $\U$ (see Lemmas~\ref{lem:restrict_eval} and \ref{lem:lift_eval}). From Lemma~\ref{lem:support_g} we obtain that dominant eigenfunctions of $\overline{\calP}_\e$ and $\calP_{\e}$ are supported on $\overline{\U}\cup B_\delta(W^s_V(\Lambda))$ (recall the definition of $\delta$ from Notation \ref{not:delta}). Analogously, we also obtain that the dominant eigenmeasures of $\overline{\calP}_\e^*$ and $\calP_{\e}^*$ are supported on $\U \cup B_\delta(W_V^u(\Lambda)).$ From Lemma~\ref{lem:spectral_gap_QEM} we conclude that the quasi-ergodic measures of $X_\e^{\log \varphi }, X_\e^{\phi},$ and $\widetilde{X}_\e^{\phi}$ are supported on $\overline{\U}$ and therefore coincide. 
\end{proof}

\section{Proof of the theorems}\label{sec:proofs}
\subsection{Local conditioned stochastic stability}\label{sec:proof_local}

We have developed the necessary ingredients in Sections~\ref{sec:BanachSpaces} and~\ref{sec:QEMs-non-vanishing} to show the existence and conditioned stochastic stability of a quasi-ergodic measure $\nu_\e^\phi$ for the process $X_\e^\phi$ conditioned upon not leaving $\overline{V}$. Recall that this measure satisfies $\nu_\e^\phi=  g_\e(x) \mu_\e(\d x)$, where $g_\e \in \ker(\calP_\e - \lambda_\e) \cap \mathcal C^0_+(\overline{V})$ and $\mu_\e(\dx) = m_\e(x)\dx \in \ker(\calP_\e^* - \lambda_\e)$ with $m_\e \in \mathcal C^0_+(\overline{V})$ (see Lemma~\ref{lem:spectral_gap_QEM} (3)). Here, the transfer operator considered is 
\begin{align*}
    \calP_\e: \mathcal C^0(\overline{V}) &\to \mathcal C^0(\overline{V}),\\
   f(x) &\mapsto e^{\phi(x)} \E_\e\left[f\circ T_\w(x) \cdot \1_{\overline{V}}\circ T_\w(x)\right]
\end{align*}
and $\lambda_\e = r(\calP_\e)$. Following from Corollary~\ref{coro:all-qems-same}, it remains to show that the quasi-ergodic measure $\nu_\e^\phi$ satisfies $\supp \nu_\e^\phi \supset \Lambda$ to conclude the proof of Theorem~\ref{thm:qem-local}. We dedicate the rest of Section~\ref{sec:proof_local} to this statement. Recall from Remark~\ref{rmk:abuse} that if $T^n_\omega(x) \not\in \overline{V}$ we say that $T^{n+m}_\omega(x) \not\in \overline{V}$ for every $m\in \mathbb N.$
\begin{proof}[Proof of Theorem~\ref{thm:qem-local}]
It remains to show that $\Lambda\subset\supp\nu_\e^\phi$. Let $g_\e,m_\e\in\mathcal C^0_+(\overline V)$ be the functions introduced at the beginning of Section~\ref{sec:proof_local}, satisfying
\begin{equation}\label{eq:qem-support-eigenvectors}
\calP_\e g_\e=\lambda_\e g_\e,\qquad \calP_\e^*\mu_\e=\lambda_\e\mu_\e,\qquad \mu_\e(\dx)=m_\e(x)\dx,\qquad \nu_\e^\phi(\dx)=g_\e(x)m_\e(x)\dx,
\end{equation}
where $\lambda_\e=r(\calP_\e)>0$ and $\mu_\e(g_\e)=1$. By Corollary~\ref{coro:all-qems-same}, $\supp\nu_\e^\phi\subset\overline{\U}\subset V$. Since $g_\e m_\e$ is continuous, it follows that
\begin{equation}\label{eq:qem-support-normalisation}
\int_V g_\e(x)m_\e(x)\dx=1,\quad\text{and}\quad \{g_\e m_\e>0\}\subset V.
\end{equation}
For the second inclusion, observe that a strictly positive value of $g_\e m_\e$ at a point of $\partial V$ would give positive $\nu_\e^\phi$-measure to a neighbourhood of that point disjoint from $\overline{\U}$, contradicting Corollary~\ref{coro:all-qems-same}.

We define the following set
\begin{equation}\label{eq:qem-support-accessible}
A_\e^+\coloneqq\left\{T_\omega^n(x):\,x\in\Lambda,\ \omega\in\Omega_\e,\ n\in\N_0,\ T_\omega^j(x)\in V\text{ for }0\leq j\leq n\right\}.
\end{equation}
The proof is divided into four steps.

\begin{step}[1]\label{step:qem-support-zeros}
If $g_\e(x)=0$ for some $x\in\Lambda$, then $g_\e$ vanishes on $\overline{A_\e^+}$. Moreover, for every fixed sufficiently small $\e>0$, there exists $\delta=\delta(\e)>0$ such that $B_\delta(\Lambda)\subset A_\e^+$.
\end{step}
\begin{stepproof}[Proof of Step~\ref{step:qem-support-zeros}]
By the eigenfunction identity in \eqref{eq:qem-support-eigenvectors},
\begin{equation}\label{eq:qem-support-ge}
\lambda_\e g_\e(x)=e^{\phi(x)}\int_{[-\e,\e]^m}g_\e(T_\omega(x))\1_{\overline V}(T_\omega(x))\,\d\upsilon(\omega),\qquad x\in\overline V.
\end{equation}
Moreover, from \eqref{eq:me}, with  $Y=\overline V$ and $\psi=e^\phi$, we obtain
$$\lambda_\e m_\e(y)=\int_{[-\e,\e]^m}\frac{e^{\phi(T_\omega^{-1}(y))}}{|\det DT_\omega(T_\omega^{-1}(y))|}m_\e(T_\omega^{-1}(y))\1_{\overline V}(T_\omega^{-1}(y))\,\d\upsilon(\omega),\qquad y\in\overline V.$$
The integrands in \eqref{eq:qem-support-ge} and \eqref{eq:me} are non-negative. Since $\upsilon$ has full support, continuity implies that
\begin{equation}\label{eq:qem-support-zero-propagation}
\begin{aligned}
g_\e(x)=0,\quad T_\omega(x)\in V&\ \text{implies}\ g_\e(T_\omega(x))=0,\\
\ \text{and}\ m_\e(y)=0,\quad T_\omega^{-1}(y)\in V&\ \text{implies}\ m_\e(T_\omega^{-1}(y))=0.
\end{aligned}
\end{equation}

Suppose for a contradiction that $g_\e(x)=0$ for some $x\in\Lambda$. By the assumptions of Section~\ref{sec:random_perturb}, the map $\omega\mapsto T_\omega(x)$ is a submersion at $\omega=0$ and $T_0=T$. Its image therefore contains an open neighbourhood of $T(x)$ contained in $V$, on which $g_\e$ vanishes by \eqref{eq:qem-support-zero-propagation}. Intersecting this neighbourhood with $\Lambda$ gives a non-empty relatively open subset of $\Lambda$. Applying \eqref{eq:qem-support-zero-propagation} with $\omega=0$, we conclude that $g_\e$ vanishes on all its forward images under $T$. Since $T|_\Lambda$ is topologically transitive, the union of these images is dense in $\Lambda$. Thus $g_\e|_\Lambda=0$. Repeatedly applying \eqref{eq:qem-support-zero-propagation} along the paths in \eqref{eq:qem-support-accessible}, and then using continuity, implies
\begin{equation}\label{eq:qem-support-vanishing}
g_\e|_\Lambda=0\quad\text{and}\quad g_\e|_{\overline{A_\e^+}}=0.
\end{equation}

Finally, fix $\e>0$ sufficiently small. For each $p\in\Lambda$, the map $\omega\mapsto T_\omega(T^{-1}(p))$ is a submersion at $\omega=0$, with image containing an open neighbourhood of $p$. Choosing this neighbourhood inside $V$, each of its points is reached in one step from $T^{-1}(p)\in\Lambda$ and therefore belongs to $A_\e^+$ by \eqref{eq:qem-support-accessible}. Compactness of $\Lambda$ gives $\delta=\delta(\e)>0$ such that
\begin{equation}\label{eq:qem-support-neighbourhood}
B_\delta(\Lambda)\subset A_\e^+\subset V.
\end{equation}
\end{stepproof}

\begin{step}[2]\label{step:qem-support-unstable}
For $\e>0$ sufficiently small, for every $\omega\in\Omega_\e$ and $p\in\Lambda_\omega$, every sufficiently small local unstable manifold $W^u_{\mathrm{loc}}(\omega,p)$ is contained in $\overline{A_\e^+}$.
\end{step}
\begin{stepproof}[Proof of Step~\ref{step:qem-support-unstable}]
Define the set
$$\Lambda_\omega\coloneqq\left\{x\in\overline V:\,T_\omega^n(x)\in\overline V\text{ for every }n\in\Z\right\}.$$
By the cocycle identity, it is readily verified that $
T_\w(\Lambda_\w)=\Lambda_{\theta\w}.$
Choose local unstable manifolds $W^u_{\mathrm{loc}}(\omega,p)$, $p\in\Lambda_\omega$, of a uniform sufficiently small size, so that all their backward iterates remain in $V$.

By the structural stability of hyperbolic sets under small random perturbations (see, e.g., \cite[Section~4.2]{KiferLiu2006} and \cite[Proposition~1.3]{Liu1998}), the sets $\Lambda_\omega$ and their corresponding local unstable manifolds depend continuously on the noise sequence. In particular, if $\omega^{(k)}\to\omega$,
$p\in\Lambda_\omega$ and $z\in W^u_{\mathrm{loc}}(\omega,p)$, then there
exist $p_k\in\Lambda_{\omega^{(k)}}$ and
$z_k\in W^u_{\mathrm{loc}}(\omega^{(k)},p_k)$ such that $p_k\to p$ and $z_k\to z.$ 

Fix $\omega\in\Omega_\e$, $p\in\Lambda_\omega$ and $z\in W^u_{\mathrm{loc}}(\omega,p)$. For each $k\geq1$, define
\begin{equation}\label{eq:qem-support-truncated-noise}
\omega_i^{(k)}\coloneqq
\begin{cases}
0,&i<-k,\\
\omega_i,&i\geq-k.
\end{cases}
\end{equation}
Then $\omega^{(k)}\to\omega$ in the product topology. By the preceding continuity property, choose $p_k\in\Lambda_{\omega^{(k)}}$ and $z_k\in W^u_{\mathrm{loc}}(\omega^{(k)},p_k)$ such that $z_k\to z$.

For each fixed $k$, the backward orbit of $z_k$ remains in $V$. Since $\omega_i^{(k)}=0$ for every $i<-k$, we have
$$
T_{\omega^{(k)}}^{-(k+j)}(z_k)
=
T^{-j}\!\left(T_{\omega^{(k)}}^{-k}(z_k)\right)
\qquad\text{for every }j\geq0.
$$
Thus, after $k$ backward iterates, the orbit evolves exactly according to the deterministic map $T^{-1}$. Its accumulation set in that direction is a non-empty compact $T$-invariant subset of $\overline V$. Since $V$ is an isolating neighbourhood, this accumulation set is contained in $\Lambda$. Consequently,
$$\lim_{n\to\infty}\dist\bigl(T_{\omega^{(k)}}^{-n}(z_k),\Lambda\bigr)=0.$$
From \eqref{eq:qem-support-neighbourhood}, a sufficiently far backward iterate of $z_k$ lies in  $A_\e^+$. Concatenating a path from $\Lambda$ to this iterate with the finite forward segment returning to $z_k$, all of whose intermediate points lie in $V$, gives $z_k\in A_\e^+$ by \eqref{eq:qem-support-accessible}. Letting $k\to\infty$, we obtain
\begin{equation}\label{eq:qem-support-unstable-inclusion}
W^u_{\mathrm{loc}}(\omega,p)\subset\overline{A_\e^+},\qquad \omega\in\Omega_\e,\quad p\in\Lambda_\omega.
\end{equation}
\end{stepproof}

\begin{step}[3]\label{step:qem-support-positive-g}
We have $g_\e(x)>0$ for every $x\in\Lambda$.
\end{step}
\begin{stepproof}[Proof of Step~\ref{step:qem-support-positive-g}]
We first show that every backward orbit $(x_{-j})_{j\geq1}$ contained in $V$, with $$T_{\omega_{-j}}(x_{-j})=x_{-j+1}\ \text{for each }j\geq 1\ \text{for some }\omega\in\Omega_\e,$$ has a sufficiently far backward iterate lying on a local unstable manifold $W^u_{\mathrm{loc}}(\theta^{-j}\omega,p)$ for some $p\in\Lambda_{\theta^{-j}\omega}$. Choose a sufficiently small neighbourhood $U$ of $\Lambda$, with $\overline U\subset V$, such that the local stable and unstable manifold theorem applies uniformly on $U$. Since $V$ is isolating, compactness gives $N\geq1$ such that
\begin{equation}\label{eq:qem-support-finite-isolation}
\bigcap_{j=-N}^{N}T^{-j}(\overline V)\subset U.
\end{equation}
It follows that, for any backward orbit $(x_{-j})_{j\geq0}$ contained in $V$, we have $x_{-j}\in U$ for every $j\geq N$, since $x_{-j}$ is the middle point of the orbit segment $\{x_{-(j+N)},\ldots,x_{-j},\ldots,x_{-(j-N)}\}$, which has $2N$ steps and is contained in $\overline V$.
By the local unstable manifold theorem and the definition of
$\Lambda_\omega$, the backward orbit $(x_{-j})_{j\geq N}$ is
asymptotic to the backward orbit $(T_\omega^{-j}(p))_{j\geq N}$
of some point $p\in\Lambda_\omega$.
Hence, after completing the parameters to a two-sided sequence $\omega$, there exist $j\geq N$ and $p\in\Lambda_{\theta^{-j}\omega}$ such that
\begin{equation}\label{eq:qem-support-backward-plaque}
x_{-j}\in W^u_{\mathrm{loc}}(\theta^{-j}\omega,p).
\end{equation}

Suppose now, for a contradiction, that $g_\e(x)=0$ for some $x\in\Lambda$. By \eqref{eq:qem-support-normalisation}, there exists $x_0\in V$ such that $g_\e(x_0)m_\e(x_0)>0$. Applying equation~\eqref{eq:me} at $x_0$, we obtain a parameter $\omega_{-1}\in[-\e,\e]^m$ for which
$$x_{-1}\coloneqq T_{\omega_{-1}}^{-1}(x_0)\in\overline V,\qquad m_\e(x_{-1})>0.$$
Moreover, $g_\e(x_{-1})>0$. Indeed, otherwise \eqref{eq:qem-support-zero-propagation}, applied at $x_{-1}$ with $T_{\omega_{-1}}(x_{-1})=x_0\in V$, would imply $g_\e(x_0)=0$. Hence $g_\e(x_{-1})m_\e(x_{-1})>0$, and \eqref{eq:qem-support-normalisation} implies $x_{-1}\in V$.

Repeating this argument, we obtain some $\omega\in\Omega_\e$ and a backward orbit satisfying
\begin{equation}\label{eq:qem-support-positive-backward-orbit}
x_{-j}\in V,\qquad T_{\omega_{-j}}(x_{-j})=x_{-j+1},\qquad g_\e(x_{-j})m_\e(x_{-j})>0,\qquad j\geq1.
\end{equation} By \eqref{eq:qem-support-backward-plaque}, a sufficiently far backward iterate belongs to a local unstable manifold. Equations~\eqref{eq:qem-support-vanishing} and \eqref{eq:qem-support-unstable-inclusion} imply that $g_\e$ vanishes at this iterate, contradicting \eqref{eq:qem-support-positive-backward-orbit}. Therefore $g_\e(x)>0$ for every $x\in\Lambda$.
\end{stepproof}

\begin{step}[4]\label{step:qem-support-positive-m}
We have $m_\e(x)>0$ for every $x\in\Lambda$.
\end{step}
\begin{stepproof}[Proof of Step~\ref{step:qem-support-positive-m}]
Consider the annealed transfer operator
\begin{equation}\label{eq:qem-support-inverse-operator}
\mathcal L_\e f(y)
\coloneqq
\int_{[-\e,\e]^m}
\frac{e^{\phi(T_\omega^{-1}(y))}}
{|\det DT_\omega(T_\omega^{-1}(y))|}
f(T_\omega^{-1}(y))
\1_{\overline V}(T_\omega^{-1}(y))
\,\d\upsilon(\omega),
\qquad y\in\overline V.
\end{equation}
Equation~\eqref{eq:me} gives $\mathcal L_\e m_\e=\lambda_\e m_\e$. Moreover, changing variables $y=T_\omega(x)$ in \eqref{eq:qem-support-inverse-operator}, for every $f\in\mathcal C^0(\overline V)$ we obtain
\begin{align}\label{eq:qem-support-inverse-duality}
\int_{\overline V}\mathcal L_\e f(y)g_\e(y)\,\d y &=
\int_{\overline V}
f(x)e^{\phi(x)}
\int_{[-\e,\e]^m}
g_\e(T_\omega(x))
\1_{\overline V}(T_\omega(x))
\,\d\upsilon(\omega)\,\dx
\nonumber=
\lambda_\e\int_{\overline V}f(x)g_\e(x)\,\dx,
\end{align}
where the last equality follows from \eqref{eq:qem-support-ge}. Hence $
\mathcal L_\e^*(g_\e(x)\dx)=\lambda_\e g_\e(x)\dx.$

The family $T_\omega^{-1}$ is a small random perturbation of $T^{-1}$. Indeed, $T_0^{-1}=T^{-1}$ and, differentiating $T_\omega(T_\omega^{-1}(y))=y$
with respect to $\omega$, we obtain
\begin{equation}\label{eq:qem-support-inverse-nondegeneracy}
D_\omega T_\omega^{-1}(y)
=
-\bigl(DT_\omega(T_\omega^{-1}(y))\bigr)^{-1}
D_\omega F_\e(\omega,T_\omega^{-1}(y)).
\end{equation}
Since $D_\omega F_\e(\omega,T_\omega^{-1}(y))$ is surjective, \eqref{eq:qem-support-inverse-nondegeneracy} implies that $\omega\mapsto T_\omega^{-1}(y)$ is a submersion.

The arguments in Steps~\ref{step:qem-support-zeros}--\ref{step:qem-support-positive-g} therefore apply to \eqref{eq:qem-support-inverse-operator}, with $m_\e$ as the eigenfunction and $g_\e(x)\dx$ as the eigenmeasure. Repeating those arguments for the inverse random dynamics yields $m_\e(x)>0$ for every $x\in\Lambda$.
\end{stepproof} 

By Steps~\ref{step:qem-support-positive-g} and \ref{step:qem-support-positive-m}, the continuous function $g_\e m_\e$ is strictly positive at every point of $\Lambda$. Therefore every neighbourhood of each such point has positive $\nu_\e^\phi$-measure, by \eqref{eq:qem-support-eigenvectors}. Hence $\Lambda\subset\supp\nu_\e^\phi$, completing the proof.
\end{proof}

\begin{proof}[Proof of Theorem~\ref{thm:local}]
    By Corollary~\ref{coro:all-qems-same}, the quasi-ergodic measure of $X_\e^\phi$ on $\overline V$ coincides with that of $X_\e^{\log\varphi}$, where $\varphi=\chi e^\phi$ is defined in \eqref{eq:phi}. Applying Section~\ref{sec:QEMs-vanishing} with $g$ replaced by $\varphi$ gives the desired convergence, since $\varphi=g$ on $\Lambda$ and hence the corresponding equilibrium state is $\nu^g$.
\end{proof}

\subsection{Global conditioned stochastic stability} \label{sec:global-proof}

The results obtained thus far provide a strategy for approximating equilibrium states on uniformly hyperbolic sets $\Lambda$ via quasi-ergodic measures. In turn, these may be obtained by combining the dominant eigenfunctions of the annealed transfer operator $\calP_\e$ and its dual. To compute these operators, however, it is necessary to have information on the set $\Lambda$ and an isolating neighbourhood $V\supset \Lambda$ (recall the ``conditioning'' term $\1_V \circ T_\w(x)$ in the definition of $\calP_\e$). As argued in Section~\ref{sec:statements}, this assumption is often not valid in model-free, data-driven applications. Instead, we may only be able to identify an attracting or trapping region $\mathcal A$ in the state space to which trajectories converge after a sufficiently long time. In this context, the conditioning step can only be considered upon not entering this region $\mathcal A$.

Moreover, by construction and the results of Lemma~\ref{lem:spectral_gap_QEM}, it follows that the quasi-ergodic limit in \eqref{eq:qem-cond} holds for $x \in \{g_\e > 0\}$ which is only a subset of $V$. It is only natural to ask which other points in $M$ also satisfy this limit for a given quasi-ergodic measure and how this problem can be addressed in the presence of several hyperbolic invariant sets.

We focus on these two aspects in this final section, where we build on the results obtained thus far to address the \emph{global} description of a system. Here, we work under the following assumption slightly extending~\ref{hyp:global}:
\begin{hypothesis}[HG2]\label{hyp:global-global} We say that $(T,g,\mathcal A)$ satisfies Hypothesis~\ref{hyp:global-global} if:
\begin{enumerate}[label = (\roman*)]
    \item $(T,g,\mathcal A)$ satisfies \ref{hyp:global}, 
    \item the topological pressures are pairwise distinct, i.e.~the topological pressure on each hyperbolic basic set $\Lambda_i\subset M\setminus \mathcal A$ is different from the rest, and
    \item $T:\Lambda_i\to\Lambda_i$ is topologically mixing for every basic set $\Lambda_i\subset M\setminus\mathcal A$. 
\end{enumerate}
\end{hypothesis}
In the remainder of this section, we relabel only the basic sets contained in $M\setminus\mathcal A$ as $\Lambda_1,\ldots,\Lambda_k$, and write $\Lambda=\bigcup_{i=1}^k\Lambda_i$. The condition $T(\overline{\mathcal A})\subset\mathcal A$ implies that each basic set is either contained in $\mathcal A$ or disjoint from $\overline{\mathcal A}$. The basic sets contained in $\mathcal A$ are excluded from the
graph and the ordering constructed below. In the following, we construct a dynamical filtration of the space $M$ in the spirit of \cite{Conley1978} and \cite{Shub1987}. Let us recall some standard definitions.

\begin{definition}[Adapted filtration]
    A \emph{filtration adapted} to $T$ is a nested family of $\mathcal C^r$, compact, codimension 0 submanifolds $\mathbf{M} = \{M_i\}_{i =0}^{\ell}$, i.e.~with
    $\emptyset = M_0 \subset M_1 \subset \dots \subset M_\ell = M,$
    and such that $T(M_i) \subset \mathrm{int}\left( M_i\right), i = 0,\dots, \ell$.
\end{definition}

For the conditioned dynamics, we use a relative filtration $\mathcal A=M_0\subset M_1\subset\cdots\subset M_\ell=M$. The positive levels $M_i$, $i\geq1$, are compact
codimension-zero submanifolds satisfying $T(M_i)\subset\operatorname{int}(M_i)$, and
$\overline{\mathcal A}\subset\operatorname{int}(M_1)$. The initial level $M_0=\mathcal A$ represents the absorbing region and need not be compact. Observe that this filtration is stable under small perturbations since $T(\overline{\mathcal A})\subset\mathcal A$.

For a filtration $\mathbf{M}$ adapted to $T$, we denote by $K_i(\mathbf{M})$ the maximal (compact) $T$-invariant subset of $M_i \setminus M_{i-1}$ for each $i = 1, \dots, \ell$, namely:
\[K_i(\mathbf{M}) = \bigcap_{n \in \mathbb{Z}} T^n(M_i \setminus M_{i-1}).\]

\begin{definition}[Filtration ordering]\label{def:order}
    We define a \emph{preorder} $\gg$ on the $\Lambda_i$'s as follows: $\Lambda_i \gg \Lambda_j$ if and only if $(W^u(\Lambda_i) \setminus \Lambda_i)\cap(W^s(\Lambda_j) \setminus \Lambda_j) \neq \emptyset$, i.e.~there are points that are both in the unstable set of $\Lambda_i$ and in the stable set of $\Lambda_j$. If there exists a sequence $\Lambda_{i_0} \gg \dots \gg \Lambda_{i_{r-1}} = \Lambda_{i_0}$, we say that the preorder has an $r$-cycle. In the absence of cycles, we write $\Lambda_i > \Lambda_j$ if there exists a sequence such that $\Lambda_i \gg \dots \gg \Lambda_j$.
\end{definition}

Denote by $G = (V, E)$ the graph generated by the vertices $V = \{\Lambda_i\}_{i = 1}^{k}$ and edges $E$ given by the preorder $\gg$.
To construct a total order on $V$, we consider the following ordering rules:
\begin{enumerate}[label = (R\arabic*)]
    \item The order inferred by $\gg$ is always preserved, i.e.~if $\Lambda_i \gg \Lambda_j$ then $\Lambda_i > \Lambda_j$. 
    \item This rule applies only to the blocks $\mathscr G_{i_s}=(V_{i_s},E_{i_s})$ produced by the recursive construction below. If $s<q$, every vertex of $V_{i_s}$ is placed above every vertex of $V_{i_q}$. Thus the blocks are ordered by strictly decreasing maximal topological pressure.
   \item Within each subgraph, topological pressure is used only as a tie-breaker. More precisely, suppose that $R$ is the set of vertices which have not yet been ordered. Among the vertices $\Lambda\in R$ for which there is no $\widetilde{\Lambda}\in R$ satisfying $\widetilde{\Lambda}>\Lambda$, we choose the one with largest topological pressure, place it next in the ordering, and repeat the procedure with the remaining vertices.
\end{enumerate}

We apply (R1), (R2) and (R3) to the following steps:
\begin{enumerate}[label = (\arabic*)]
    \item Begin by choosing the basic set of maximal topological pressure, which we assume to be $\Lambda_{i_0}$ without loss of generality. Denote by $\mathscr{G}_{i_0}$ the graph consisting of vertices
$$
V_{i_0}\coloneqq\{\Lambda_i\in V:\Lambda_i>\Lambda_{i_0}\}\cup\{\Lambda_{i_0}\}
$$
and their corresponding edges given by the relation $\gg$. We order the vertices of $\mathscr G_{i_0}$ according to (R1), using (R3) to choose between the maximal elements whenever there is more than one possible choice.
    \item Consider the graph $G_1 = G \setminus \mathscr{G}_{i_0} = (V_1, E_1)$.
    \item Choose the basic set in the graph $G_1$ of maximal topological pressure, which we denote by $\Lambda_{i_1}$. Denote by $\mathscr{G}_{i_1}$ the graph consisting of the vertices $$V_{i_1} \coloneqq
\{\Lambda_i \in V \setminus V_{i_0}\,:\, \Lambda_i > \Lambda_{i_1}\}
\cup\{\Lambda_{i_1}\}$$ and their corresponding edges. Order the vertices of $\mathscr G_{i_1}$ according to (R1), using (R3) whenever there is more than one maximal element, and compare $\mathscr G_{i_1}$ with $\mathscr G_{i_0}$ via (R2).
    \item 
    Consider the graph $G_2=G_1\setminus\mathscr G_{i_1}$ and repeat Step~(3) on the remaining graph.
    \item Repeat until the remaining graph is empty. Writing $G_0=G$, let $t$ be the total number of blocks selected, so that $G_t=G_{t-1}\setminus\mathscr G_{i_{t-1}}=\emptyset,$ $1\leq t\leq k.$
\end{enumerate}
Observe that this recursion is finite as there are finitely many basic sets. Once the ordering is established, we relabel the sets $\Lambda_i$ such that $\Lambda_i > \Lambda_j$ if $i > j$ while also relabelling the indexes $i_s, s \in \{0,\dots,t-1\},\quad 1\leq t\leq k$, so that they refer to the $t$ basic sets of ``maximal'' topological pressure identified in step (3).

Following \cite{Shub1987}, we call this a \emph{filtration ordering}. Let us provide an example of the ordering construction proposed above. 

\begin{example}\label{ex:filtration-ordering}
Consider ten basic sets $\Lambda_1,\ldots,\Lambda_{10}$ with heteroclinic relations
$$
\Lambda_5\gg\Lambda_2\gg\Lambda_4,\quad
\Lambda_1\gg\Lambda_3\gg\Lambda_4,\quad
\Lambda_6\gg\Lambda_3,
$$
and
$$
\Lambda_4\gg\Lambda_8,\quad
\Lambda_7\gg\Lambda_8\gg\Lambda_{10},\quad
\Lambda_9\gg\Lambda_{10},
$$
with no other relations. Suppose that the corresponding topological pressures are
$$
\begin{array}{c|cccccccccc}
\Lambda_i
&\Lambda_1&\Lambda_2&\Lambda_3&\Lambda_4&\Lambda_5
&\Lambda_6&\Lambda_7&\Lambda_8&\Lambda_9&\Lambda_{10}\\ \hline
P_{\mathrm{top}}(T,\Lambda_i,\psi)
&3.2&6.4&4.5&10&5.8&2.1&7.1&9&1.2&8
\end{array}
$$
The graph and the resulting ordering are illustrated in Figure~\ref{fig:filtration-ordering-example}.

\begin{figure}[h!]
\centering
\begin{subfigure}[b]{\textwidth}
\centering
\resizebox{\textwidth}{!}{%
\begin{tikzpicture}[
    >=stealth,
    every node/.style={font=\small},
    vertex/.style={
        draw,
        rounded corners=1.5pt,
        minimum width=1.35cm,
        minimum height=.72cm,
        align=center,
        inner sep=2pt
    },
    selected/.style={
        vertex,
        double,
        double distance=1pt
    }
]
\node[vertex] (L5) at (1,4) {$\Lambda_5$\\[-1mm]\scriptsize $P=5.8$};
\node[vertex]   (L1) at (4,4) {$\Lambda_1$\\[-1mm]\scriptsize $P=3.2$};
\node[vertex]   (L6) at (7,4) {$\Lambda_6$\\[-1mm]\scriptsize $P=2.1$};

\node[vertex]   (L2) at (1,2.5) {$\Lambda_2$\\[-1mm]\scriptsize $P=6.4$};
\node[vertex]   (L3) at (5.5,2.5) {$\Lambda_3$\\[-1mm]\scriptsize $P=4.5$};

\node[selected] (L4) at (3.25,1) {$\Lambda_4$\\[-1mm]\scriptsize $P=10$};

\node[vertex]   (L7) at (10.7,4) {$\Lambda_7$\\[-1mm]\scriptsize $P=7.1$};
\node[selected] (L8) at (10.7,1) {$\Lambda_8$\\[-1mm]\scriptsize $P=9$};

\node[vertex]   (L9) at (15,4) {$\Lambda_9$\\[-1mm]\scriptsize $P=1.2$};
\node[selected] (L10) at (15,-1) {$\Lambda_{10}$\\[-1mm]\scriptsize $P=8$};

\draw[->] (L5) -- (L2);
\draw[->] (L2) -- (L4);
\draw[->] (L1) -- (L3);
\draw[->] (L6) -- (L3);
\draw[->] (L3) -- (L4);

\draw[->] (L4) -- (L8);
\draw[->] (L7) -- (L8);

\draw[->] (L8) -- (L10);
\draw[->] (L9) -- (L10);

\draw[dashed,rounded corners=2pt]
    (0.1,.25) rectangle (8,4.8);
\node[anchor=west] at (0.25,4.6)
    {$\mathscr G_{i_0}$};

\draw[dashed,rounded corners=2pt]
    (9.55,.25) rectangle (11.85,4.8);
\node[anchor=west] at (9.7,4.6)
    {$\mathscr G_{i_1}$};

\draw[dashed,rounded corners=2pt]
    (13.7,-1.75) rectangle (16.3,4.8);
\node[anchor=west] at (13.85,4.6)
    {$\mathscr G_{i_2}$};
\end{tikzpicture}
}
\caption{Heteroclinic graph and the blocks selected by the algorithm. An arrow $\Lambda_i\to\Lambda_j$ denotes the preorder $\Lambda_i\gg\Lambda_j$. Double boxes indicate the maximal-pressure basic set selected at each stage.}
\end{subfigure}

\vspace{.5cm}

\begin{subfigure}[b]{\textwidth}
\centering
\resizebox{\textwidth}{!}{%
\begin{tikzpicture}[
    >=stealth,
    every node/.style={font=\small},
    vertex/.style={
        draw,
        rounded corners=1.5pt,
        minimum width=1.35cm,
        minimum height=.72cm,
        align=center,
        inner sep=2pt
    },
    selected/.style={
        vertex,
        double,
        double distance=1pt
    }
]

\node[vertex]   (O5)  at (0,-4) {$\Lambda_5$};
\node[vertex]   (O2)  at (1.7,-4) {$\Lambda_2$};
\node[vertex]   (O1)  at (3.4,-4) {$\Lambda_1$};
\node[vertex]   (O6)  at (5.1,-4) {$\Lambda_6$};
\node[vertex]   (O3)  at (6.8,-4) {$\Lambda_3$};
\node[selected] (O4)  at (8.5,-4) {$\Lambda_4$};

\node[vertex]   (O7)  at (10.7,-4) {$\Lambda_7$};
\node[selected] (O8)  at (12.4,-4) {$\Lambda_8$};

\node[vertex]   (O9)  at (14.6,-4) {$\Lambda_9$};
\node[selected] (O10) at (16.3,-4) {$\Lambda_{10}$};

\draw[->] (O5) -- (O2);
\draw[->] (O2) -- (O1);
\draw[->] (O1) -- (O6);
\draw[->] (O6) -- (O3);
\draw[->] (O3) -- (O4);

\draw[->] (O4) -- node[above,font=\scriptsize] {(R2)} (O7);
\draw[->] (O7) -- (O8);

\draw[->] (O8) -- node[above,font=\scriptsize] {(R2)} (O9);
\draw[->] (O9) -- (O10);

\draw[dashed,rounded corners=2pt]
    (-.75,-4.65) rectangle (9.25,-3.35);

\draw[dashed,rounded corners=2pt]
    (9.95,-4.65) rectangle (13.15,-3.35);

\draw[dashed,rounded corners=2pt]
    (13.85,-4.65) rectangle (17.05,-3.35);
\end{tikzpicture}%
}
\caption{Total order obtained after applying (R1)--(R3). An arrow $\Lambda_i \to \Lambda_j$ denotes the total order $\Lambda_i > \Lambda_j$.}
\end{subfigure}
\caption{Example of the filtration ordering. Dashed boxes show the subgraphs $\mathscr G_{i_0}$, $\mathscr G_{i_1}$ and $\mathscr G_{i_2}$.}
\label{fig:filtration-ordering-example}
\end{figure}

We now apply the algorithm explicitly. Since $\Lambda_4$ has maximal topological pressure, Step~(1) selects $\Lambda_{i_0}=\Lambda_4$. Thus,
$$
V_{i_0}=\{\Lambda_1,\Lambda_2,\Lambda_3,\Lambda_4,\Lambda_5,\Lambda_6\},
$$
so $\mathscr G_{i_0}$ is the subgraph generated by these vertices. Initially, its maximal elements are $\Lambda_1$, $\Lambda_5$, and $\Lambda_6$. Rule~(R3) selects $\Lambda_5$, which has the largest pressure among them. After removing $\Lambda_5$, the vertex $\Lambda_2$ becomes maximal, and (R3) selects it next. Notice that, although $\Lambda_2$ has larger pressure than $\Lambda_5$, it could not be selected first because $\Lambda_5>\Lambda_2$. Continuing recursively gives
$$
\Lambda_5>\Lambda_2>\Lambda_1>\Lambda_6>\Lambda_3>\Lambda_4.
$$
After removing $\mathscr G_{i_0}$, the remaining vertices are $\Lambda_7,\Lambda_8,\Lambda_9,\Lambda_{10}$, with $\Lambda_7\gg\Lambda_8\gg\Lambda_{10}$ and $\Lambda_9\gg\Lambda_{10}$. Amongst these, the maximal-pressure vertex is $\Lambda_8$, hence $\Lambda_{i_1}=\Lambda_8$ and $V_{i_1}=\{\Lambda_7,\Lambda_8\}$, giving $\Lambda_7>\Lambda_8$. Removing this block leaves $\Lambda_9$ and $\Lambda_{10}$, so $\Lambda_{i_2}=\Lambda_{10}$, $V_{i_2}=\{\Lambda_9,\Lambda_{10}\}$, and $\Lambda_9>\Lambda_{10}$.
Since
$$
P_{\mathrm{top}}(T,\Lambda_4,\psi)>P_{\mathrm{top}}(T,\Lambda_8,\psi)>P_{\mathrm{top}}(T,\Lambda_{10},\psi),
$$
rule~(R2) gives $\mathscr G_{i_0}>\mathscr G_{i_1}>\mathscr G_{i_2}$. Therefore, the resulting total order is
\begin{equation}\label{eq:filtration-order-example}
\Lambda_5>\Lambda_2>\Lambda_1>\Lambda_6>\Lambda_3>\Lambda_4>\Lambda_7>\Lambda_8>\Lambda_9>\Lambda_{10}.
\end{equation}
This also illustrates the role of (R3), showing that pressure is used only to choose among the maximal vertices available at each stage.
\end{example}

Recall that in Hypothesis~\ref{hyp:global}~\ref{item:gh1iii} we require that $\Lambda$ has no cycles. In particular, this assumption allows for the order construction above to be well-defined and sets us up to apply the following result.

\begin{proposition}[{\cite[Theorem~2.3]{Shub1987}}] \label{prop:shub} Let $T:M \to M$ be a homeomorphism and let $\Lambda = \Lambda_1 \sqcup \dots \sqcup \Lambda_k$ be the union of $k$ closed invariant sets containing all $\alpha$ and $\omega$ limit sets of $T$. Then $\Lambda_i = K_i(\mathbf{M}), i = 1, \dots, k,$ for some filtration $\mathbf{M}$ adapted to $T$ if and only if $\Lambda$ has no cycles and the ordering by indices is a filtration ordering. 
\end{proposition}



Applying Proposition~\ref{prop:shub}, we obtain an adapted
filtration
$$
\emptyset=F_0\subset\cdots\subset F_q
\subset F_{q+1}\subset\cdots\subset F_{q+k}=M,
$$
with the $q$ basic sets contained in $\mathcal A$ first,
followed by $\Lambda_1,\ldots,\Lambda_k$.
This ordering is compatible with the heteroclinic relations
because $W^u(\Lambda_j)\subset\mathcal A$ whenever
$\Lambda_j\subset\mathcal A$.

Every point of $\overline{\mathcal A}$ eventually enters
$\operatorname{int}(F_q)$. By compactness and forward
invariance, choose $N\geq0$ we define
$T^N(\overline{\mathcal A})\subset\operatorname{int}(F_q)$.
Define
$$
M_0=\mathcal A,\quad\text{and}\quad M_i=T^{-N}(F_{q+i}),\quad\text{for each}\ i=1,\ldots,k.
$$
We obtain
$\mathcal A=M_0\subset M_1\subset\cdots\subset M_k=M$.
Writing $C_i=\overline{M_i\setminus M_{i-1}}$, we have
$$
\Lambda_i\subset\operatorname{int}(C_i)\quad\text{and}\quad\bigcap_{n\in\mathbb Z}T^n(C_i)=\Lambda_i.
$$

Observe that the inclusion $T(M_i)\subset\operatorname{int}(M_i)$ 
and $T(\overline{\mathcal A})\subset\mathcal A$ persist under sufficiently small random perturbations. 

It is immediate to see that for every $x\in W^s(\mathcal A)$, there exist $\e>0$ and $N\in\mathbb N$ such that $\tau(\omega,x)\leq N$ for every $\omega\in\Omega_\e$ (see Theorem~\ref{thm:global-global}, item~\ref{item:global-global-global-A}).

Let $x \in W^s(\Lambda_j) \cap M_j$ for some $j \in \{1,\dots, k\}$ and assume that $\Lambda_j \in V_{i_s}$ for some $s \in \{0, \dots, t-1\}$. Then we show that the quasi-ergodic measure corresponding to the conditioned statistics for an initial condition $x$ when conditioned upon not entering any trapping region $\mathcal A \subset M_{i_{s}-1}$ coincides with the quasi-ergodic measure for the process starting in $M_{i_s} \setminus M_{i_s -1}$ and conditioned upon remaining there. We formalise this statement in the following theorem.

\begin{theorem}[Global conditioned stochastic stability II, $\partial = \mathcal A$]\label{thm:global-global} Assume that $(T,g,\mathcal A)$ satisfies Hypothesis~\ref{hyp:global-global}. Let $\phi:M\to\mathbb R$ be such that $e^\phi$ is continuous on $M$ and of class $\mathcal C^{r-1}$ on an open neighbourhood of every retained basic set $\Lambda_i$. Assume that $e^\phi=g$ on $\Lambda=\bigcup_{i=1}^k\Lambda_i$. Consider the filtration ordering presented above and let $i_0>\dots>i_{t-1}$ denote the indices identified after the final relabelling step. The following statements hold true:
\begin{enumerate}[label = (\roman*)]
    \item \label{item:global-global-M_j} There exists $\delta>0$ such that for every $x \in \cup_{i_s \leq j < i_{s-1}} W^s(\Lambda_j) \cap M_{i_{s-1} -1}$ for some $s \in \{0, \dots, t-1\}$, with $i_{-1} = k+1$, we have that for every $0<\e <\delta$ and $\varphi: M \to \R$ bounded and measurable, 
    \begin{equation*}
            \lim_{n \to \infty}\E_x^\phi\left[\frac{1}{n}\sum_{i= 0}^{n-1}\varphi \circ X_i^\phi \, \bigg|\, \tau >n \right] = \int \varphi \,\d\nu_{\e, i_s}^\phi,
    \end{equation*}
    where $\nu_{\e, i_s}^\phi$ is the unique quasi-ergodic measure of $X_\e^\phi$ on $M_{i_s} \setminus M_{i_s - 1}$ when conditioned upon staying within $M_{i_s} \setminus M_{i_s - 1}$ such that $\Lambda_{i_s} \subset \supp \nu^\phi_{\e,i_s}$ (i.e.~of Theorem~\ref{thm:qem-local}); and Theorem~\ref{thm:global-existence} holds.

    \item \label{item:global-global-global} Given $i_s \leq j < i_{s-1}$, for some $s \in \{0, \dots, t-1\}$, and $x\in W^s(\Lambda_j)$ there exists $\delta =\delta(x)$ such that for every $0<\e<\delta$ and every bounded and measurable function $\varphi: M\to \R$, we have that
    \begin{equation*}
            \lim_{n \to \infty}\E_x^\phi\left[\frac{1}{n}\sum_{i= 0}^{n-1}\varphi \circ X_i^\phi \, \bigg|\, \tau >n \right] = \int \varphi \,\d\nu_{\e, i_s}^\phi.
    \end{equation*}

    \item \label{item:global-global-global-stoch-stab} The quasi-ergodic measures $\nu^\phi_{\e, i_s}$ of item~\ref{item:global-global-M_j} converge in weak-$^*$ as $\e \to 0$ to the unique equilibrium state on $\Lambda_{i_s}$ associated with the potential $\phi - \log |\left.\det DT\right|_{E^u}|$. Moreover, if $s=0$ we only need to assume Hypothesis~\ref{hyp:global} for this item to hold, which implies Theorem \ref{thm:global}.
    
\item \label{item:global-global-global-A} If $x$ eventually enters $\mathcal A$ under the deterministic dynamics, then there exist $\e>0$ and $N\in\mathbb N$ such that $
\tau(\omega,x)\leq N$ for every $\omega\in\Omega_\e$.
\end{enumerate}

\end{theorem}

\begin{proof}[Proof of Theorem~\ref{thm:global-global}]
We prove item~\ref{item:global-global-M_j}. Assume first that $x \in W^s(\Lambda_j)$ for $j \geq i_0$ so that $s = 0$ and $M_{i_{-1}-1} = M_k = M$. For $\e > 0$ small enough, consider the global transfer operator acting on suitable functions $f:M\setminus\mathcal A\to\R$ given by
$$
\widehat{\mathcal P}_\e f(x)\coloneqq e^{\phi(x)}\E_\e\left[f\circ T_\w(x)\cdot\1_{M\setminus\mathcal A}\circ T_\w(x)\right].
$$
Given a subset $X\subset M$, we denote by $\mathcal P_{\e,X}$ the transfer operator restricted to $X$, i.e.
$$
\calP_{\e,X}f(x)\coloneqq e^{\phi(x)}\E_\e\left[f\circ T_\w(x)\cdot\1_X\circ T_\w(x)\right]
$$
for suitable $f:X\to\R$. Both operators are strong Feller by definition of the perturbation (cf.~Section~\ref{sec:random_perturb}). We write $C_i=\overline{M_i\setminus M_{i-1}}$ for every $i\in\{1,\dots,k\}$.

Under Hypothesis~\ref{hyp:global}, Lemma~\ref{lem:PeSG} gives a spectral gap for the dominant operator $\calP_{\e,C_{i_0}}$. For each $i$, choose an isolating neighbourhood $V_i$ of $\Lambda_i$ with $\overline{V_i}\subset\operatorname{int}(C_i)$ such that $e^\phi$ is of class $\mathcal C^{r-1}$ on a neighbourhood of $\overline{V_i}$. Choose a cutoff function $\chi_i\in\mathcal C_c^\infty(V_i)$ with $0\leq\chi_i\leq1$ and $\chi_i=1$ near $\Lambda_i$. Define $
g_i\coloneqq\chi_i e^\phi$ and $\calP_{g_i}f(x)\coloneqq g_i(x)f(T(x)).$ The spectral-radius formula in Remark~\ref{rmk:deterministic-spectral-data} gives
$$
\lambda_i^0\coloneqq r\!\left(\calP_{g_i}\big|_{B^{t,s}(T,V_i)}\right)=\exp P_{\mathrm{top}}\left(T,\Lambda_i,\log g-\log|\det DT|_{E^u}|\right)>0.
$$
Moreover, Theorem~\ref{thm:thm51Baladi} and Remark~\ref{rmk:deterministic-spectral-data} give
$$
r_{\mathrm{ess}}\!\left(\calP_{g_i}\big|_{B^{t,s}(T,V_i)}\right)\leq Q^{t,s}(T,g_i)<\lambda_i^0.
$$
Define the localised annealed operators by
$$
\calP_{\e,g_i}f(x)\coloneqq g_i(x)\E_\e\left[f(T_\w(x))\right]\quad \text{and}\quad\calP_{0,g_i}=\calP_{g_i},
$$
where functions are extended by zero outside $V_i$. Using the estimates of Lemma~\ref{lem:LY-Baladi}, apply \cite[Theorem~1 and Corollary~1]{KellerLiverani1999} to the family $\bigl(\frac{1}{\lambda_i^0}\calP_{\e,g_i}\bigr)_{\e\geq0}$ acting on $B^{t,s}(T,V_i)$.

From Proposition~\ref{prop:SpecBaladi} and the same proofs of Lemmas~\ref{lem:restrict_eval} and~\ref{lem:lift_eval}, with outer domain $C_i$ and  weight $g_i$, we obtain
$$
r(\calP_{\e,C_i})=r\!\left(\calP_{\e,g_i}\big|_{B^{t,s}(T,V_i)}\right)\longrightarrow\lambda_i^0>0\quad\text{as }\e\to0.
$$

Moreover, each operator $\calP_{\e,C_i}^2$ is compact, since $\calP_{\e,C_i}$ is strong Feller property, as in the proof of Lemma~\ref{lem:sg_PeU}. Thus, for fixed sufficiently small $\e>0$, the finite-dimensional peripheral spectral subspaces and the spectral decomposition give constants $K_i>0$ and $q_i\in\mathbb N_0$ such that
$$
\|\calP_{\e,C_i}^{\,n}\|_{\mathcal C^0(C_i)}\leq K_i n^{q_i}r(\calP_{\e,C_i})^n\quad\text{for every }n\geq1.
$$
The factor $n^{q_i}$ accounts for possible nontrivial Jordan blocks associated with the peripheral eigenvalues.

Let $\delta>0$ be sufficiently small and set
$$
\lambda_\e\coloneqq r(\calP_{\e,C_{i_0}})\quad\text{and}\quad\rho_\e\coloneqq\max_{1\leq j<i_0}r(\calP_{\e,C_j}).
$$
We use the convention $\rho_\e=0$ if $i_0=1$, in which case the lower surviving region is empty and $\widehat{\calP}_\e^n\1_{\mathcal A}=0$ for every $n\geq1$. By Hypothesis~\ref{hyp:global} and the convergence of the local spectral radii as $\e\to0$, after decreasing $\delta$ if necessary we have $\rho_\e<\lambda_\e$ for every $0<\e<\delta$.

Since the filtration is adapted to $T_\omega$ for $\e>0$ sufficiently small, an orbit starting in $M_{i_0-1}$ can only move through the finitely many layers $C_1,\ldots,C_{i_0-1}$. Using the polynomial power bounds for $\calP_{\e,C_j}$ on each layer and arguing inductively over the filtration, there exist $K>0$ and $q\in\mathbb N_0$ such that $
0\leq\widehat{\calP}_\e^n\1_{M_{i_0-1}}(x)\leq K n^q\rho_\e^n$ for every $x\in M_{i_0-1}$ and $n\geq1$. Therefore
\begin{equation}\label{eq:limitdecay}
0\leq\frac{1}{\lambda_\e^n}\widehat{\calP}_\e^n\1_{M_{i_0-1}}(x)\leq K n^q\left(\frac{\rho_\e}{\lambda_\e}\right)^n\xrightarrow[n\to\infty]{}0.
\end{equation}
Arguing in the same way over the whole finite filtration, and using that $\lambda_\e$ is the largest of the local spectral radii, there exist $K'>0$ and $q'\in\mathbb N_0$ such that $\|\widehat{\calP}_\e^n\|\leq K'n^{q'}\lambda_\e^n$. Hence $r(\widehat{\calP}_\e)\leq\lambda_\e$. On the other hand, 
every path staying in $C_{i_0}$ also survives in $M\setminus\mathcal A$. Positivity of the operators considered implies that  $\|\calP_{\e,C_{i_0}}^n\|\leq\|\widehat{\calP}_\e^n\|$ for every $n\geq1$, and hence $\lambda_\e\leq r(\widehat{\calP}_\e)$. 
Therefore $r(\widehat{\calP}_\e)=\lambda_\e=r(\calP_{\e,C_{i_0}})$. In particular, if $g_{\e,i_0}\in\mathcal C^0(M)$ satisfies $\widehat{\mathcal P}_\e g_{\e,i_0}=\lambda_\e g_{\e,i_0}$ with $\lambda_\e=r(\widehat{\mathcal P}_\e)$, then
$$
g_{\e,i_0}=\1_{M\setminus M_{i_0-1}}g_{\e,i_0}+\1_{M_{i_0-1}}g_{\e,i_0}=\1_{M\setminus M_{i_0-1}}g_{\e,i_0}.
$$
Since $\mathbf M$ is a dynamical filtration, it follows that if $x\in C_{i_0}$ is such that $g_{\e,i_0}(x)>0$ and $g_{\e,i_0}(T_\w(x))>0$, then $T_\w(x)\in C_{i_0}$. Therefore, $\mathcal P_{\e,C_{i_0}}g_{\e,i_0}=\lambda_\e\1_{C_{i_0}}g_{\e,i_0}$. The same steps can be performed for $\mu_\e\in(\mathcal C^0(M))^*$ satisfying $\widehat{\mathcal P}_\e^*\mu_\e=\lambda_\e\mu_\e$, from which we obtain that $\mu_\e=\1_{M_{i_0}}\mu_\e$ and $\calP_{\e,C_{i_0}}^*\mu_\e=\lambda_\e\1_{C_{i_0}}\mu_\e$.

An analogue of Lemma~\ref{lem:lift_eval} is also applicable here. We may induce eigenfunctions of $\widehat{\calP}_\e$ from eigenfunctions of $\calP_{\e,C_{i_0}}$, preserving their eigenvalue and setting up a link between both spectra. It follows that $\widehat{\mathcal P}_\e$ has a spectral gap in $\mathcal C^0(M)$ (cf.~Lemma~\ref{lem:sg_PeU}).

Finally, by construction of the filtration $\mathbf M$, for every $\delta>\e>0$ and $x\in W^s(\Lambda_j)$, $j\geq i_0$, there exists $N=N(x,\e)$ such that $T_\w^N(x)\in\Lambda_{i_0}$ for some $\w\in\Omega_\e$. This yields $g_{\e,i_0}(x)>0$ and thus $\bigcup_{j\geq i_0}W^s(\Lambda_j)\subseteq\{g_{\e,i_0}>0\}$. Applying Lemma~\ref{lem:spectral_gap_QEM}, we conclude that for $x\in W^s(\Lambda_j)$, $j\geq i_0$, the conditioned Birkhoff averages converge
$$
\lim_{n\to\infty}\E_x^\phi\left[\frac1n\sum_{i=0}^{n-1}h\circ X_i^\phi\,\bigg|\,\tau>n\right]=\int h\,\d\nu_{\e,i_0}^\phi,
$$
where $\nu_{\e,i_0}^\phi(\dx)=g_{\e,i_0}(x)\mu_\e(\dx)$ is the unique quasi-ergodic measure for the process $X_\e^\phi$ conditioned upon remaining within $M_{i_0}\setminus M_{i_0-1}$ of the local Theorem~\ref{thm:qem-local}. In fact, this holds for every $x\in\{g_{\e,i_0}>0\}$. Observe that this conclusion only uses Hypothesis~\ref{hyp:global} since we only need a spectral gap of $\mathcal P_{\e,C_{i_0}}$. The proof also implies items \ref{item:gex1}-\ref{item:gex3} of Theorem~\ref{thm:global-existence}.

For $x\in W^s(\Lambda_j)$, $i_s\leq j<i_{s-1}$ and $s=1,\ldots,t-1$, we argue in a similar fashion. The same argument above applies simply replacing $M$ by $M_{i_{s-1}-1}$ and $i_0$ by $i_s$. In this setting, one needs to assume Hypothesis~\ref{hyp:global-global} when repeating the computations in \eqref{eq:limitdecay} with $\lambda_\e$ replaced by $r(\mathcal P_{\e,C_{i_s}})$. Observe that since $\{g_{\e,i_s}>0\}\subset M_{i_{s-1}-1}$, we obtain the quasi-ergodic limit for
$$
x\in\left(\bigcup_{i_s\leq j<i_{s-1}}W^s(\Lambda_j)\right)\cap M_{i_{s-1}-1}.
$$

We prove item~\ref{item:global-global-global}. Let $s\in\{0,\ldots,t-1\}$, $i_s\leq j<i_{s-1}$, and $x\in W^s(\Lambda_j)$. From Hypothesis~\ref{hyp:global-global}, for every $\e>0$ small enough there exists a filtration $\mathcal A=M_0\subsetneq M_1\subsetneq M_2\subsetneq\ldots\subsetneq M_k$ of $M$ adapted to $T_\omega$ for every $\w\in\Omega_\e$. Observe that $x\in M_{j_0}$ for some $j_0\geq j$. From item~\ref{item:global-global-M_j}, if $x\in M_j$ the proof is done. Assume then that $x\notin M_j$ (in particular $x\notin\Lambda_j$). It is sufficient to construct a new filtration $\widetilde{\mathbf M}$ of $M$ adapted to $T_\omega$ for every $\omega\in\Omega_\delta$ respecting the filtration ordering, i.e. $\Lambda_i\subset\widetilde M_i$ for every $i\in\{1,\dots,k\}$, and such that $x\in\widetilde M_j$, where $\delta=\delta(x)$ is a small constant that will be chosen appropriately during the proof.

Since $x\in W^s(\Lambda_j)$ and $\Lambda_j\subset\operatorname{int}(M_j)$, there exists $N=N(x)\geq1$ such that $T^N(x)\in\operatorname{int}(M_j)$. Define $\widetilde M_i\coloneqq T^{-N}(M_i)$, $i=0,\ldots,k$. Since $T$ is a diffeomorphism, $\widetilde{\mathbf M}=\{\widetilde M_i\}_{i=0}^k$ is again a nested filtration of $M$. Moreover, since each $\Lambda_i$ is $T$-invariant,
$$
\Lambda_i=T^{-N}(\Lambda_i)\subset T^{-N}(M_i)=\widetilde M_i,
$$
and $x\in\operatorname{int}(\widetilde M_j)$. The filtration is adapted to $T$, since
$$
T(\widetilde M_i)=T^{-N}(T(M_i))\subset T^{-N}(\operatorname{int}(M_i))=\operatorname{int}(\widetilde M_i).
$$
Furthermore,
$$
K_i(\widetilde{\mathbf M})=\bigcap_{m\in\mathbb Z}T^m\bigl(T^{-N}(M_i\setminus M_{i-1})\bigr)=T^{-N}\bigl(K_i(\mathbf M)\bigr)=\Lambda_i,
$$
so the filtration ordering is unchanged. Since there are only finitely many sets $\widetilde M_i$ and each inclusion $T(\widetilde M_i)\subset\operatorname{int}(\widetilde M_i)$ is strict, the $\mathcal C^0$-closeness of $T_\omega$ to $T$ implies that there exists $\delta=\delta(x)>0$ such that $T_\omega(\widetilde M_i)\subset\operatorname{int}(\widetilde M_i)$ for every $i=0,\ldots,k$, every $\omega\in\Omega_\e$, and every $0<\e<\delta$. Hence $\widetilde{\mathbf M}$ is adapted to all sufficiently small random perturbations and contains $x$ at level $j$.

Observe that item~\ref{item:global-global-global-stoch-stab} is a direct consequence of the local conditioned stochastic stability shown in Theorem~\ref{thm:local}. Indeed, $\nu_{\e,i_s}^\phi$ is the unique quasi-ergodic measure of $X_\e^\phi$ on the closure of the isolating neighbourhood $V=\mathrm{int}(M_{i_s}\setminus M_{i_s-1})$ such that $\Lambda_{i_s}\subset\supp\nu_{\e,i_s}^\phi$ (as established in item~\ref{item:global-global-M_j}). If $s=0$, Hypothesis~\ref{hyp:global} is sufficient to imply the existence of $\nu_{\e,i_0}$, as mentioned in the proof of item~\ref{item:global-global-M_j} above. Theorem~\ref{thm:local} implies the second part of item~\ref{item:global-global-global-stoch-stab} and therefore Theorem~\ref{thm:global}.

Finally, item~\ref{item:global-global-global-A} follows from the fact that, if $T^N(x)\in\mathcal A$ for some $N\in\mathbb N$, then, since $\mathcal A$ is open, by continuity there exists $\e=\e(x)>0$ such that $T_\omega^N(x)\in\mathcal A$ for every $\omega\in\Omega_\e$. Hence $\tau(\omega,x)\leq N$ for every $\omega\in\Omega_\e$.
\end{proof}

\section{Lean formalisation}
\label{sec:lean-formalisation}

Theorems~\ref{thm:qem-local}, \ref{thm:local}, \ref{thm:global-existence}, \ref{thm:global} and~\ref{thm:global-global} have been formalised in Lean~4.35.0 using Mathlib. The source code is available on the GitHub repository \url{https://github.com/Manzatto/conditioned-stochastic-stability-hyperbolic}. The statements and proofs are collected in \texttt{MainTheorems.lean}, in the namespace \texttt{ConditionedStochasticStabilityHyperbolic}, with the following declarations:
\begin{center}
\small
\setlength{\tabcolsep}{3pt}
\begin{tabular}{@{}lcccccc@{}}
\hline
Theorem
& \ref{thm:qem-local}
& \ref{thm:local}
& \ref{thm:global-existence}
& \ref{thm:global}
& \ref{thm:compact_embedding}
& \ref{thm:global-global} \\
Lean
& \texttt{Theorem2\_7}
& \texttt{Theorem2\_8}
& \texttt{Theorem2\_10}
& \texttt{Theorem2\_11}
& \texttt{Theorem3\_15}
& \texttt{Theorem6\_5} \\
\hline
\end{tabular}
\end{center}

In this repository, we have also formalised every result cited throughout the paper, if leveraged in any of the proofs above. This includes the construction of the anisotropic Banach spaces in Section~\ref{sec:BanachSpaces} and the estimates for transfer operators underlying Theorem~\ref{thm:thm51Baladi}, following
Baladi--Tsujii~\cite{Baladi-Tsujii2008} and Baladi~\cite{Baladi2018-book}; the spectral and equilibrium state identifications of Baladi--Tsujii \cite{Baladi-Tsujii2008} and Gou\"ezel--Liverani \cite{Gouezel-Liverani2008} recalled in Remark~\ref{rmk:deterministic-spectral-data}; and the spectral perturbation results of Keller--Liverani \cite{KellerLiverani1999} in Appendix~\ref{sec:appendix}.
The correspondence between the paper and their formal statements is documented in \texttt{docs/StatementAudit.md}.

The development contains no \texttt{sorry} placeholders or additional axioms. An audit of the project declarations, including private declarations, verifies that their axiom dependencies are confined to \texttt{propext}, \texttt{Classical.choice} and \texttt{Quot.sound}. The compiler and dependency versions are fixed by \texttt{lean-toolchain} and \texttt{lake-manifest.json}. Compilation and axiom dependencies can be checked from the repository root by running:
\begin{verbatim}
lake exe cache get
lake build
lake env lean Audit.lean
\end{verbatim}

\section*{Acknowledgements}
    We thank Jeroen S.W. Lamb for very valuable comments and insightful discussions throughout the preparation of this manuscript. BBC thanks José F. Alves and Martin Rasmussen for a very careful and detailed review of an earlier version of the paper.
    The authors gratefully acknowledge support from the EPSRC Centre for Doctoral Training in Mathematics of Random Systems: Analysis, Modelling and Simulation (EP/S023925/1). BBC is supported by a Chapman Fellowship at Imperial College London. MMC’s research has been supported by an Imperial College President’s PhD scholarship, the Australian Research Council Laureate Fellowship (FL230100088), and the São Paulo Research Foundation (FAPESP, grant no. 2025/26997-9). 

\addtocontents{toc}{\protect\setcounter{tocdepth}{1}} 

\subsection*{Use of AI declaration.} The Lean~4 formalisation of this paper has been done by Claude Opus 5.5 (Anthropic) and ChatGPT Sol 6.1 (OpenAI).

\subsection*{Data availability statement} Data sharing is not applicable to this article as no datasets were generated or analysed during the current study.

\subsection*{Conflict of interest statement} The authors have no competing interests to declare that are relevant to the content of this article.

\addtocontents{toc}{\protect\setcounter{tocdepth}{2}}

\bibliographystyle{abbrvnat}
\bibliography{refs}

\newpage
\appendix
\markboth{APPENDIX}{APPENDIX}
\section{Some useful results from functional analysis}\label{sec:appendix}

For the sake of completeness and the reader's convenience, we recall two well-known spectral theorems which we use in Section~\ref{sec:QEMs-vanishing}. The first one is a particular version of the well-established perturbation theory developed in \cite{KellerLiverani1999}. The second is a simple abstract result that provides conditions under which an operator defined on two different spaces shares the same eigenvectors for eigenvalues which are not in the essential spectrum (see also \cite[Appendix~A.2]{Baladi2018-book}).

\begin{theorem}[{\cite[Corollary~1]{KellerLiverani1999}}]\label{thm:keller} Let $(\mathcal{B},\|\cdot\|)$ be a Banach space and let $|\cdot|\leq \|\cdot\|$ be a second norm defined on $\mathcal B$. 
Consider a family of bounded linear operators ${(P_\e)_{\e \geq 0}:(\mathcal{B},\|\cdot\|)} \to (\mathcal{B},\|\cdot\|)$ such that the following holds:  
\begin{enumerate}[label = (\roman*)]
    \item There are constants $C_1, M>0$ such that for every $\e\geq 0$ small enough, $$|P_\e^n| \leq C_1 M^n,$$ 
    for all $n\in\mathbb N$.
    \item There exist $C_2,C_3>0$ and $\alpha \in(0,1)$, $\alpha<M,$ such that for every $\e\geq 0$
    $$\|P_\e^n f\| \leq C_2 \alpha^n \|f\| + C_3 M^n |f|$$
    for every $n\in\mathbb N$ and for each $f\in \mathcal{B}$.
    \item The closed unit ball of $(\mathcal{B},\|\cdot\|)$ is $|\cdot|$-compact in  the completion of $\mathcal B$ with respect to the norm $|\cdot|$.
    \item There exists a monotone upper-semicontinuous function $\tau:[0,\infty)\to[0,\infty),$ such that $\tau(\e) >0$ if $\e >0$ and
    $$ \sup\{\left|P_0 f - P_\e f\right|\colon f \in \mathcal{B},\|f\|\leq 1 \} \leq \tau(\e) \xrightarrow[]{\e\to 0} 0.$$
\end{enumerate}
Moreover, let $\lambda$ be an isolated simple eigenvalue of $P_0$ with $\lambda>\alpha$. Let $\delta_0 >0$ and let $r\in(\alpha,\lambda)$ be such that $B_{\delta_0}(\lambda)\cap \sigma(P_0) =\{\lambda\}$. For $\delta < \delta_0$, consider the projection
$$\Pi_\e^{(\lambda,\delta)} \coloneqq \frac{1}{2\pi i} \int_{\partial B_{\delta}(\lambda)} (z- P_\e)^{-1} \d z.$$
Then the following holds:
\begin{enumerate}[label = (\arabic*)]
     \item \label{it:keller1} there exist constants $K_1= K_1(\delta,r) >0$ and $\eta>0$, such that
    $$\left|\Pi_\e^{(\lambda,\delta)} f - \Pi_0^{(\lambda,\delta)}f \right| \leq K_1 \tau(\e)^\eta \|f\|\ \text{for all }f\in \mathcal{B}.$$
    \item \label{it:keller2} There are constants $\delta = \delta(r)>0$ and $K_2 = K_2(\delta,r)>0$ such that $\|\Pi_{\e}^{(\lambda,\delta)} f\| \leq K_2 |\Pi_\e^{(\lambda,\delta)}f|$ for all $f\in \mathcal{B}$, $\delta\in(0,\delta_0]$ and $\e$ small enough.
    \item\label{it:keller3} For each $\e>0$ small enough $P_\e$ has a unique eigenvalue $\lambda_\e \in B_\delta(\lambda)\cap \sigma(P_\e)$. Moreover, $\lambda_\e \xrightarrow[]{\e\to 0} \lambda$ and $\lambda_\e$ is a simple eigenvalue of $P_\e$.
    \item\label{it:keller4} If $P_0$ has a spectral gap then each $P_\e$ also has a spectral gap for $\e>0$ small enough.  
\end{enumerate}
\end{theorem}

\begin{proposition}[{\cite[Appendix~A]{Baladi-Tsujii2008}}]\label{prop:SpecBaladi}
    
Let $\mathcal{B}$ be a Hausdorff topological linear space and let $(\mathcal{B}_1, \|\cdot\|_1)$ and $(\mathcal{B}_2, \|\cdot\|_2)$ be Banach spaces that are continuously embedded in $\mathcal{B}$. Suppose that there is a subspace $\mathcal{B}_0 \subset \mathcal{B}_1 \cap \mathcal{B}_2$ that is dense both in the Banach spaces $(\mathcal{B}_1, \|\cdot\|_1)$ and $(\mathcal{B}_2, \|\cdot\|_2)$. Let $\mathcal{L} : \mathcal{B} \rightarrow \mathcal{B}$ be a continuous linear map, which preserves the subspaces $\mathcal{B}_0$, $\mathcal{B}_1$, and $\mathcal{B}_2$. Suppose that the restrictions of $\mathcal{L}$ to $\mathcal{B}_1$ and $\mathcal{B}_2$ are bounded operators whose essential spectral radii are both strictly smaller than some number $\rho > 0$. Then the eigenvalues of $\mathcal{L}|_{\mathcal{B}_1}$ and $\mathcal{L}|_{\mathcal{B}_2}$ in $\{ z \in \mathbb{C} \mid |z| > \rho \}$ coincide. Furthermore, the corresponding generalised eigenspaces coincide and are contained in $\mathcal{B}_1 \cap \mathcal{B}_2$.
\end{proposition}

\end{document}